\documentclass[11pt]{article}
\usepackage[T1]{fontenc}
\usepackage{lmodern}
\usepackage[margin=1.05in]{geometry}
\usepackage{amsmath,amssymb,amsthm,mathtools}
\usepackage{booktabs,array,microtype}
\usepackage[hidelinks]{hyperref}
\hypersetup{pdftitle={Invariant subsets and invariable generation of Ewens random permutations and cyclations},
  pdfauthor={Charles Burnette},
  pdfsubject={Invariant subsets and invariable generation under the Ewens sampling formula},
  pdfkeywords={Ewens sampling formula; invariant subsets; invariable generation; random sumsets; random cyclations},bookmarksdepth=3}
\usepackage{enumitem}
\setlist[enumerate]{itemsep=3pt,topsep=5pt}

\newcommand{\Pp}{\mathbb P}
\newcommand{\Pref}{\mathbb P_{\!*}}
\newcommand{\E}{\mathbb E}
\newcommand{\TV}{d_{\mathrm{TV}}}
\newcommand{\Law}{\mathcal L}
\newcommand{\ind}{\mathbf 1}

\newtheorem{theorem}{Theorem}[section]
\newtheorem{lemma}[theorem]{Lemma}
\newtheorem{proposition}[theorem]{Proposition}
\newtheorem{corollary}[theorem]{Corollary}
\newtheorem*{mainA}{Theorem A}
\newtheorem*{mainB}{Theorem B}
\theoremstyle{definition}

\theoremstyle{remark}
\newtheorem{remark}[theorem]{Remark}
\numberwithin{equation}{section}
\allowdisplaybreaks[1]
\title{Invariant subsets and invariable generation\\
of Ewens random permutations and cyclations}
\author{Charles Burnette\\[3pt]
\small Department of Mathematics\\
\small Xavier University of Louisiana\\
\small New Orleans, LA 70125-1098\\
\small \texttt{cburnet2@xula.edu}}
\date{}
\begin{document}
\maketitle
\begin{abstract}
A permutation fixes a set of size \(k\) if and only if some of its
cycle lengths sum to \(k\). We study these sums for Ewens permutations
of parameter \(\theta\) and for random cyclations, whose loop-length profile
has the Ewens law at \(\theta=1/2\). At
\(\theta_c=1/(2\log 2)\), two independent Poisson cycle profiles
almost surely have infinitely many common positive sums, even after
all lengths through any fixed cutoff are deleted. At this parameter,
two independent Ewens permutations invariably generate \(S_n\) with
probability of order \((\log n)^{-1/2}\) as \(n\to\infty\), whereas three do so with a
positive limiting probability. For every \(\theta>0\), the least
number of samples giving a positive limit inferior agrees with the
corresponding Poisson sumset threshold. We determine the limiting
generation probabilities in both strict ranges and give the precise
scope of earlier bounds for the Ewens law.

We also determine the order of the probability of fixing a \(k\)-set at
every fixed \(\theta>0\), uniformly for \(1\le k\le n/2\). Below
\(\theta_c\), conditioning on one prescribed \(k\)-set being invariant
gives an approximation to the number of invariant \(k\)-sets by
exchanges between two Poisson profiles. The error decreases as a power
of \(k\) in total variation. The difference \(n-2k\) accounts for
invariant sets near the complement of the prescribed set.
\end{abstract}

\section{Introduction}

We say that a set $A$ is invariant under a function $f$ if $f(A) = A.$ A permutation fixes a set if and only if that set is a union of its
cycles. If several permutations each fix a set of size \(k\), we can
conjugate them so that they all fix the same \(k\)-set. For \(0<k<n\)
these conjugates generate an intransitive subgroup of the symmetric
group \(S_n\). We ask how often independent permutations have a common invariant size
and how likely a single permutation is to fix a set of a prescribed size.

Let \(\pi_n\) be an Ewens permutation of \([n]=\{1,\ldots,n\}\) with
parameter \(\theta>0\), and write \(K_n\) for its number of cycles.
A permutation with \(c\) cycles has probability
\(\theta^c/[\theta(\theta+1)\cdots(\theta+n-1)]\). Write \(I_{n,k}\)
for the number of \(k\)-sets fixed by \(\pi_n\), and put
\(i_\theta(n,k)=\Pp(I_{n,k}>0)\).
The empty set and the whole set give \(I_{n,0}=I_{n,n}=1\).
Complementation gives \(I_{n,k}=I_{n,n-k}\), so we take \(n\ge2\),
\(1\le k\le n/2\), and write \(d=n-2k\) for the difference between
the two side sizes. All logarithms are natural unless a base is indicated.

\subsection*{Main results}
\addcontentsline{toc}{subsection}{Main results}
A family of permutations \emph{invariably generates} \(S_n\) if every
choice of their conjugates generates \(S_n\). Let \(G_{n,m}(\theta)\)
be the probability that \(m\) independent Ewens permutations invariably
generate \(S_n\). The corresponding infinite model has independent
counts \(Z_j\sim\operatorname{Pois}(\theta/j)\) of copies of each
length \(j\ge1\). Its \emph{sumset} is the subset of \(\mathbb Z_{\ge0}\)
consisting of the totals of all finite selections of copies, including
the empty selection. Each attainable total occurs once in the sumset.
Write
\(p_{0,m}(\theta)\) for the probability that \(m\) independent such
sumsets intersect only at zero. Let \(s_\theta\) and \(m_\theta\) be the
least positive integers \(m\) for which \(p_{0,m}(\theta)>0\) and
\(\liminf_{n\to\infty}G_{n,m}(\theta)>0\), respectively, with value
\(\infty\) if no such \(m\) exists. Brito, Fowler, Junge, and Levy
\cite[Eq.~(1.3)]{brito} use \(m_\alpha\), with parameter \(\alpha\),
for the threshold defined by an infimum over all \(n\ge2\). We denote
their all-degree Ewens threshold by \(m_\theta^{\mathrm{all}}\).

\begin{mainA}[Common sums and invariable generation]
Put\pdfbookmark[3]{Theorem A. Common sums and invariable generation}{main-theorem-A} \(\theta_c=1/(2\log 2)\) and \(\lambda_0=2\theta_c=1/\log 2\).
For every \(\theta>0\), \(s_\theta=m_\theta\).
The common value is infinite for \(\theta\ge\lambda_0\).
At the first critical parameter, the threshold and the large-degree
two-sample probability satisfy
\[
     s_{\theta_c}=m_{\theta_c}=3,\qquad
     G_{n,2}(\theta_c)\asymp (\log n)^{-1/2}
     \quad(n\to\infty).
\]
For \(\theta\ge\theta_c\), two independent Poisson sumsets almost surely
have infinitely many common positive elements, simultaneously after the
deletion, from both profiles, of all lengths through any fixed cutoff.
For each integer \(m\ge2\) and fixed \(\theta>0\) in either strict range,
\[
 \lim_{n\to\infty}G_{n,m}(\theta)=
 \begin{cases}
  (1-2^{-m})p_{0,m}(\theta)>0,&m(1-\theta\log2)>1,\\
  0,&m(1-\theta\log2)<1 .
 \end{cases}
\]
When \(m(1-\theta\log2)<1\), the \(m\)-fold Poisson intersection
is almost surely infinite simultaneously after every fixed lower cutoff.
\end{mainA}

For \(0<\theta<\lambda_0\), put \(x_\theta=(1-\theta\log 2)^{-1}\).
When \(x_\theta\) is not an integer, Theorem~A gives
\(m_\theta=\lceil x_\theta\rceil\), the value proposed by Brito,
Fowler, Junge, and Levy \cite{brito}. At \(x_\theta=2\), the
two-sample probability tends to zero at the rate given in Theorem~A,
answering the first Ewens case of \cite[Question~1]{brito}. At each
higher jump \(x_\theta=m\ge3\), Corollary~\ref{cor:threshold-equivalence}
reduces the choice between \(m\) and \(m+1\) to deciding whether
\(p_{0,m}(\theta)>0\). With two profiles, any nonempty exchange that preserves both
totals gives a common positive sum. For three profiles, consider the
lists \((1,2),(3),(4)\). Exchanging the first two lists leaves all
three totals unchanged, yet no positive sum occurs in all three.

The distinction between the two degree conventions already matters
below \(\theta_c\). There the large-degree threshold is two, while no
pair invariably generates \(S_6\). Proposition~\ref{prop:all-degree-threshold}
gives \(m_\theta^{\mathrm{all}}=\max\{3,m_\theta\}\).
Thus degree six alone supplies the lower bound of three under the
all-degree convention at \(\theta_c\). The assertion about
\(G_{n,2}(\theta_c)\) concerns arbitrarily large degrees.

An additional obstruction comes from transitive subgroups. Two
\(4\)-cycles in \(S_4\), for instance, are odd and have no proper invariant
size, yet one may choose their conjugates to be equal. For
\(\theta>\lambda_0\) in even degree, an invariant half places a sample
in a proper transitive subgroup with probability bounded away from
zero, contrary to \cite[Theorem~1.3]{brito}. Let \(F_n\) be the event
that some transitive subgroup of \(S_n\) contains \(\pi_n\) and does
not contain \(A_n\). At \(\lambda_0\), \(\Pp(F_n)\) tends to zero
and has order \((\log n)^{-1/2}\) in even degrees
(Proposition~\ref{prop:endpoint-overgroups}). The class of cycle
laws also matters. Under the printed strong logarithmic condition
of \cite{brito}, changing the fixed-point count alone can give every
sample a fixed point. Section~\ref{sec:earlier-generation} treats
these assertions in the notation of the original paper.

At \(\theta=1\), Pemantle, Peres, and Rivin \cite{pemantle-peres-rivin}
proved positive probability of invariable transitivity for four samples.
Eberhard, Ford, and Green \cite{efg-invariable} established the
four-sample assertion for \(S_n\) and proved that three samples fail
with probability tending to one. Their estimates for common sizes
and signs also give the ingredients for the uniform four-sample
limit \(15p_{0,4}(1)/16\). Theorem~A gives the limit throughout the
strict Ewens ranges. See Eberhard \cite[Section~1]{eberhard-ewens}
for ordinary generation under the Ewens law.

Now fix a \(k\)-set \(A\). Conditional on its invariance, the two
restrictions of \(\pi_n\) are independent Ewens permutations.
A competitor near \(A\) exchanges equal totals across the two sides.
A competitor \(B\) near \(A^c\) retains \(a\) labels of \(A\).
It therefore uses \(k-a\) of the \(k+d\) labels of \(A^c\), leaving
\(a+d\) there unused. For two independent Poisson profiles, let
\(U_a,U'_a\) count selections of total length \(a\), distinguishing
copies of the same length. Using the same profiles in both sums, put
\[
 W_\infty=\sum_{a\ge0}U_aU'_a,\qquad
 W_d=W_\infty+\sum_{a\ge0}U_aU'_{a+d}\quad(0\le d<\infty).
\]

Put
\begin{equation}\label{eq:rate-function}
 D(\theta)=2\theta-\lambda_0+
          \lambda_0\log\bigl(\lambda_0/(2\theta)\bigr).
\end{equation}

\begin{mainB}[Invariant sets of a prescribed size]
For each fixed\pdfbookmark[3]{Theorem B. Invariant sets of a prescribed size}{main-theorem-B} \(0<\theta<\theta_c\) there is an \(\eta=\eta(\theta)>0\) such
that, uniformly for \(1\le k\le n/2\) with \(d=n-2k\),
\begin{align*}
 \TV\bigl(\Law(I_{n,k}\mid A\text{ is invariant}),\Law(W_d)\bigr)
     &=O_\theta(k^{-\eta}),\\
 \frac{i_\theta(n,k)}{\E I_{n,k}}&=c_d+O_\theta(k^{-\eta}),
 \qquad c_d=\E(1/W_d).
\end{align*}
In this range each \(W_d\), including \(W_\infty\), is almost surely finite and at
least one. The constants \(c_d\) are uniformly bounded away from zero,
\(c_0=c_\infty/2\), and \(c_d\to c_\infty\) as \(d\to\infty\).
At \(\theta=\theta_c\),
\[
 \frac{i_{\theta_c}(n,k)}{\E I_{n,k}}\asymp(1+\log k)^{-1/2}.
\]
For every fixed \(\theta>0\), the second column of
Table~\ref{tab:ewens-regimes} gives the order of \(i_\theta(n,k)\)
uniformly for \(1\le k\le n/2\).
\end{mainB}

To see how the two theorems are related, let \(\Sigma_k\) be the set of
attainable totals in one Poisson profile through length \(k\). In two
independent profiles let \(V_k\) count the pairs of selections of
equal total, including the empty pair, using copies of individual
length at most \(k\). Writing \(H_k=\sum_{j=1}^k1/j\), we have
\begin{equation}\label{eq:intro-reciprocal-comparison}
 \frac{i_\theta(n,k)}{\E I_{n,k}}
 \asymp_\theta e^{-\theta H_k}\E_\theta|\Sigma_k|
 =\E_{\theta,\theta}(1/V_k),\qquad n\ge2k .
\end{equation}
Indeed, Proposition~\ref{cor:finite-support-comparison} compares
\(i_\theta(n,k)\) with \(\E_\theta|\Sigma_k|/k\), and
\eqref{s3:first} gives \(\E I_{n,k}\asymp_\theta k^{\theta-1}\).
Since \(e^{-\theta H_k}\asymp_\theta k^{-\theta}\), the comparison follows.
The equality is the reciprocal exchange identity of
Lemma~\ref{lem:reciprocal-exchange}. The shifted exchanges in \(W_d\)
give the dependence on \(d\) in Theorem~B. At equal halves each invariant
set has an invariant complement, and \(W_0=2W_\infty\).

\begin{table}[htbp]
\centering\small
\caption{Invariant \(k\)-sets, generation thresholds, and transitive
overgroups under the Ewens law. The second column gives the order up
to factors depending only on \(\theta\), with \(L=1+\log k\).
The third column gives the least \(m\) for which \(p_{0,m}(\theta)>0\).
The last column gives the probability of a transitive overgroup missing \(A_n\).}
\label{tab:ewens-regimes}
\renewcommand{\arraystretch}{1.22}
\setlength{\tabcolsep}{5pt}
\begin{tabular}{@{}llll@{}}
\toprule
Parameter & \(i_\theta(n,k)\) & \(m_\theta=s_\theta\) & \(\Pp(F_n)\)\\
\midrule
\(0<\theta<\theta_c\) & \(k^{\theta-1}\) & \(2\) & \(o(1)\)\\
\(\theta=\theta_c\) & \(k^{\theta_c-1}L^{-1/2}\) & \(3\) & \(o(1)\)\\
\(\theta_c<\theta<\lambda_0\)
 & \(k^{\theta-1-D(\theta)}L^{-3/2}\)
 & \(\lceil x_\theta\rceil\) if \(x_\theta\notin\mathbb Z\) & \(o(1)\)\\
 & & \(x_\theta\) or \(x_\theta+1\) if \(x_\theta\in\mathbb Z\) & \\
\(\theta=\lambda_0\) & \(L^{-1/2}\) & \(\infty\) & \(o(1)\)\\
\(\theta>\lambda_0\) & \(1\) & \(\infty\)
 & \(\asymp_\theta 1\) for even \(n\)\\
\bottomrule
\end{tabular}
\end{table}

Adding a profile can only shrink the common intersection, so
\(p_{0,m}(\theta)>0\) precisely for the finite integers \(m\ge m_\theta\).
At \(\lambda_0\), Proposition~\ref{prop:endpoint-overgroups} gives the
sharper even-degree estimate \(\Pp(F_n)\asymp(\log n)^{-1/2}\).

A single profile has \(\theta\log 2+o(1)\) copies on average in
\((2^{v-1},2^v]\) as \(v\to\infty\). A merged pair has
\(2\theta\log 2+o(1)\). Each copy doubles the number of selections,
and each block doubles the scale of the lengths. The mean number of
copies per block reaches one at \(\theta_c\) for a pair and at
\(\lambda_0\) for one profile. These counts also explain the two
logarithmic powers in the second column. Between \(\theta_c\) and
\(\lambda_0\), the lower bound uses a centered walk that stays
nonnegative for \(b\) dyadic blocks and returns to zero, an event of
order \(b^{-3/2}\). At \(\lambda_0\) the endpoint is free, and the
probability of staying nonnegative has order \(b^{-1/2}\).
Here \(b\) is of order \(\log k\).

Conditioning on the invariance of \(A\) size-biases the law of \(I_{n,k}\).
Thus \(i_\theta(n,k)/\E I_{n,k}\) is exactly the conditional mean of
\(1/I_{n,k}\). To condition instead on \(I_{n,k}>0\), divide the
prescribed-set mass at each multiplicity \(j\) by \(j\) and normalize
(Corollary~\ref{s3:exist-cor}).

For \(1/2\le\theta<\theta_c\), \(W_\infty\) is almost surely finite
and has infinite mean. Indeed, the Poisson generating function gives
\[
 \sum_{j\ge0}\E U_j\,z^j
 =\exp\!\left(\theta\sum_{\ell\ge1}\frac{z^\ell}{\ell}\right)
 =(1-z)^{-\theta}.
\]
Thus \(\E U_j=h_j(\theta)=[z^j](1-z)^{-\theta}
\sim j^{\theta-1}/\Gamma(\theta)\) as \(j\to\infty\), and independence gives
\(\E W_\infty=\sum_{j\ge0}h_j(\theta)^2\). This series diverges
precisely for \(\theta\ge1/2\).

For the uniform permutation, Eberhard, Ford, and Green \cite{efg} give
the sharp order \(i_1(n,k)\asymp k^{-\delta}(1+\log k)^{-3/2}\),
where \(\delta=D(1)=1-(1+\log\log 2)/\log 2\).
The intermediate range in Theorem~B extends this order to fixed Ewens
parameters on both sides of one. For \(\theta=1\), Green and
Sawhney \cite{gs} take \(n\to\infty\) first and obtain
an asymptotic equivalent with a smooth positive factor periodic in
\(\log_2 k\). They conjecture that this factor is nonconstant. For fixed \(k\) at \(\theta=1\), Diaconis, Fulman, and Guralnick
\cite[Theorem~4.1]{dfg} express the limiting count as a polynomial
in independent Poisson cycle counts. The same polynomial applies at
general \(\theta\), with means \(\theta/j\), by \eqref{s2:feller}.
In Arratia's continuum model the closure of the finite subset sums
of a Poisson process of intensity \(\theta\,dx/x\) contains \(1\) with
positive probability precisely when \(\theta>\lambda_0\). For the
supercritical implication see the discussion of Tenenbaum's result
following Problem~22.1 in \cite[Section~22]{ArratiaScale}.
Theorem~B treats exact attainment of an integer and gives the rate
at equality.

Pippenger's random cyclation \cite{pippenger} is obtained by
pairing the endpoints of \(n\) labeled intervals uniformly at random.
Its loop lengths determine a conjugacy class in \(S_n\), distributed
as the class of an Ewens permutation at \(\theta=1/2\). Conditional on the absence of cyclation loops of lengths
one and two, the graph obtained by contracting the intervals is a
uniform labeled simple \(2\)-regular graph. Write \(I_{n,k}^{[3]}\) for its number of \(k\)-vertex unions
of components. In the following corollary \(c_d\) is evaluated at
\(\theta=1/2\), and \(c_d^{[3]}\) uses the same two Poisson profiles
with all copies of lengths one and two removed
(Corollary~\ref{s3:simple-cor}).

\begin{corollary}[Cyclations and simple \(2\)-regular graphs]
\label{cor:cyclation-consequences}
The probability that the classes of two independent random
\(n\)-cyclations invariably generate \(S_n\) tends to
\(3p_{0,2}(1/2)/4>0\). For \(k\to\infty\), uniformly in \(n\ge2k\),
the probabilities of a \(k\)-interval union of loops and a
\(k\)-vertex union of graph components satisfy
\[
 \begin{aligned}
 \sqrt{\frac{\pi k(n-k)}{n}}\,i_{1/2}(n,k)
      &=c_d+O(k^{-\eta}),\\
 \sqrt{\frac{\pi k(n-k)}{n}}\,\Pp(I_{n,k}^{[3]}>0)
      &=e^{-3/4}c_d^{[3]}+O(k^{-\eta}),
 \end{aligned}
 \qquad d=n-2k,
\]
where \(\eta=\delta/(1+\delta)=0.079250\ldots\) and \(\delta=D(1/2)=D(1)\).
In particular, for \(n=2k\), multiplying the two probabilities by
\(\sqrt{k}\) gives limits \(c_\infty/\sqrt{2\pi}\) and
\(e^{-3/4}c_\infty^{[3]}/\sqrt{2\pi}\), respectively.
\end{corollary}

The generation limit follows from Theorem~A. For the two
probabilities, apply Theorem~B and Corollary~\ref{s3:simple-cor},
using \(h_m(1/2)=(\pi m)^{-1/2}(1+O(m^{-1}))\).
The restriction to lengths at least three multiplies \(h_m(1/2)\)
by \(e^{-3/4}+O(m^{-1})\), the probability of avoiding the two
shortest lengths in degree \(m\). Here \(k,n-k,n\) are all at
least \(k\), and \(\eta<1\).

For two simple \(2\)-regular graphs of equal size, D'Ascenzo and
Michielan \cite[Section~3.1 and Appendix~A]{das} prove tightness of
the number of distinct common sizes. Their estimates, together with truncation, give
convergence in law of the number of equal-size exchanges.
Corollary~\ref{s3:simple-cor} gives a rate in total variation and also
allows unequal sizes. For uniform permutations, Qian \cite{qian}
estimates the joint probability of fixing a \(k\)-set and having no
cycle shorter than \(r\). We will treat higher moments of \(I_{n,k}\) and
the mean number of exchanges between two simple graphs in a follow-up paper.

\subsection*{Outline of the proofs}
\addcontentsline{toc}{subsection}{Outline of the proofs}
Merge two independent Poisson profiles and mark each copy by its
side. Given the lengths, the first side is a uniformly chosen subset
of the merged copies. Every subset with the same total gives an equal
exchange. Thus the conditional mean reciprocal number of exchanges
is the number of attainable totals divided by the number of subsets.
Green and Sawhney \cite[Eq.~(1.5), Definition~6.1, and Lemma~6.2]{gs}
use this normalized support. Averaging over the merged profile
changes its intensity from \(2\theta/j\) to \(\theta/j\) and gives the
equality in \eqref{eq:intro-reciprocal-comparison}.

Section~\ref{sec:prelim} gives the estimates for Ewens cycle counts.
Section~\ref{sec:conditional} proves the subcritical approximation
and the critical common-sum theorem. At the first jump, a centered
walk records the excess number of copies over dyadic blocks.
In Section~\ref{sec:generation} we complete short profiles having
no positive sum common to all samples by adding long cycles and
prescribing an odd sign. One sample whose transitive overgroups all
contain \(A_n\) then suffices. Under the original product law the
exceptional overgroup events are independent, giving a two-sample
error of \(O((\log n)^{-4\theta_c/5})\), where
\(4\theta_c/5>1/2\). For the single-permutation estimates we compare
\(i_\theta(n,k)\) with the expected number of attainable totals
through length \(k\), divided by \(k\), and estimate this expectation
in each parameter range. Section~\ref{sec:earlier-generation} returns
to the assertions of \cite{brito}.

\section{Cycle profiles under the Ewens sampling formula}
\label{sec:prelim}

A uniform cyclation on \(n\) labeled intervals is obtained by choosing a
uniform perfect matching of their \(2n\) endpoints. Following interval
edges and matched pairs alternately traces disjoint loops. Pippenger
\cite{pippenger} gives their length profile. To recover it, traverse a
loop on a specified \(j\)-set of intervals. Its traversal gives a cyclic
order of the labels and a sign at each label recording its entrance
endpoint. The matched pairs join each exit to the next entrance. The two
traversals of the same loop reverse the order and flip every sign. They
remain distinct for \(j=1,2\), when reversal fixes the cyclic order.
There are therefore \(2^{j-1}(j-1)!\) loops on the specified set.
For one interval its two endpoints must be paired. For two intervals
there are two ways to pair their endpoints across the intervals.

For the rising factorial we write
\[
 x^{\overline m}=\prod_{i=0}^{m-1}(x+i),\qquad x^{\overline 0}=1.
\]

Choose a traversal of each loop independently and uniformly, and read its cyclic order as a cycle on the interval labels. The loop lengths alone determine the conjugacy class of the resulting permutation in \(S_n\). For a specified permutation with \(c\) cycles, its directed cyclic lists admit \(2^n\) sign assignments. They give distinct pairs of a matching and a choice of traversals. Each of these pairs has probability \(2^{-c}\) times the uniform matching probability. The resulting probability of the specified permutation is
\[
 \frac{2^n\,2^{-c}}{(2n-1)!!}
   =\frac{(1/2)^c}{(1/2)^{\overline n}},
\]
which is its probability under the Ewens law at \(\theta=1/2\). See also the endpoint exploration of \cite[Section~1.2]{abewens} and the matching double cosets of \cite[Theorem~4.4]{ds}.

For the rest of the paper, \(\theta>0\) is a fixed Ewens parameter. Statements about the matching itself use \(\theta=1/2\). Put
\[
 h_m(\theta)=[z^m](1-z)^{-\theta}=\frac{\theta^{\overline m}}{m!},
 \qquad h_m=h_m(\theta),
 \qquad H_m=\sum_{j=1}^m\frac1j.
\]
Take \(h_0=1\), \(h_m=0\) for \(m<0\), and \(H_0=0\). Recall that \(H_m=\log m+\gamma+O(1/m)\) for \(m\ge1\), where \(\gamma\) is Euler's constant. At \(\theta=1/2\), \(h_m=4^{-m}\binom{2m}{m}\). Write \(K_{n,j}\) for the number of cycles of length \(j\) and \(K_n=\sum_j K_{n,j}\). We use \(0^0=1\) in products of marks.
In \(O_\theta\), \(\ll_\theta\), and \(\asymp_\theta\), the implicit constants may depend on \(\theta\). Other subscripts indicate further allowed dependence. An unsubscripted implied constant is absolute unless a dependence is specified locally.

The marked series and the resulting profile law are
\begin{align}
 \E\prod_{j=1}^n w_j^{K_{n,j}}
 &=\frac{[z^n]}{h_n}\exp\!\left(\theta\sum_{j=1}^n w_jz^j/j\right),
 \label{eq:marked}\\
 \Pp(K_{n,j}=a_j,\ 1\le j\le n)
 &=\frac1{h_n}\prod_{j=1}^n\frac{(\theta/j)^{a_j}}{a_j!},
 \qquad a_j\in\mathbb Z_{\ge0},\quad \sum_j ja_j=n.
 \label{eq:profile}
\end{align}
For general \(\theta\), these identities follow by assigning the weight \(\theta^{K_n}\) and counting permutations with the prescribed profile. The one-loop weights \(\theta/j\) put the profile in the framework of logarithmic combinatorial structures of Arratia, Barbour, and Tavaré \cite{abt2000,abt}. If \(Z_j\), \(j\ge1\), are independent with \(Z_j\sim\operatorname{Pois}(\theta/j)\), then
\begin{equation}\label{eq:condition}
 (K_{n,j})_{1\le j\le n}\ \overset d=\
 (Z_j)_{1\le j\le n}\ \bigg|\ \sum_{j=1}^n jZ_j=n .
\end{equation}
The exponential factors in the Poisson masses cancel on conditioning.

We use the coefficient estimates
\begin{align}
 h_i &\asymp_\theta (i+1)^{\theta-1}, \label{s2:gamma}\\
 \frac{h_{n-s}}{h_n}
     &=1+(1-\theta)\frac{s}{n}+O_\theta(s^2/n^2)
         \qquad(n\ge1,\ 0\le s\le n/2).
 \label{s2:ratio}
\end{align}
The first estimate follows from the gamma ratio. For the second, use the
consecutive ratios \(h_{i-1}/h_i=1+(1-\theta)/(i-1+\theta)\).
For \(1\le s\le n/2\)
the indices in the product satisfy \(i>n/2\). Thus for \(n\) sufficiently
large in terms of \(\theta\) the logarithm of the ratio is
\((1-\theta)\sum_{n-s<i\le n} i^{-1}+O_\theta(s/n^2)\).
The sum is \(s/n+O(s^2/n^2)\), and \(s/n^2\le s^2/n^2\) for \(s\ge1\).
Exponentiation gives \eqref{s2:ratio}. Increasing the implied constant
covers bounded \(n\), and \(s=0\) gives equality.

For the Feller coupling \cite{abt,feller}, take independent Bernoulli variables \(\xi_i\), \(i\ge1\), of means \(\theta/(i-1+\theta)\), so \(\xi_1=1\) almost surely. The gaps between ones in \(\xi_1,\ldots,\xi_n,1\), with the last one appended at position \(n+1\), give the \(K_{n,j}\). In the infinite word the gap counts give the independent \(Z_j\). On this coupling,
\begin{equation}\label{s2:feller}
 \E\sum_{j\le b}|K_{n,j}-Z_j|\ll_\theta b/n
       \quad(1\le b\le n/2),\qquad
 \E\sum_{j\le n}|K_{n,j}-Z_j|\ll_\theta1 .
\end{equation}
Let \(L_n\) be the position of the last one among \(\xi_1,\ldots,\xi_n\).
Only the finite terminal gap and the infinite gap starting at \(L_n\)
can differ among gaps starting at \(i\le n\). The former has length
\(n+1-L_n\) and the latter has length at least \(n+1-L_n\).
Thus either gap can contribute a discrepancy through length \(b\)
only if \(L_n\in\{n-b+1,\ldots,n\}\). For \(b\le n/2\), a
union bound gives
\[
 \mathbb P\bigl(L_n\ge n-b+1\bigr)
 \le \sum_{i=n-b+1}^n\frac{\theta}{i-1+\theta}
 =O_\theta(b/n),
\]
since \(i-1+\theta\ge n/2\) throughout the sum.
For \(i>n\), independence gives
\(\Pp(\xi_i=\xi_{i+j}=1)\le\theta^2/(i-1+\theta)^2\).
A gap of length \(j\) starting at \(i\) requires ones at positions \(i\) and \(i+j\). Thus the expected number of infinite gaps of length at most \(b\) starting at \(i>n\) is at most \(\sum_{i>n}\sum_{j\le b}\Pp(\xi_i=\xi_{i+j}=1)=O_\theta(b/n)\). For \(b=n\) the same double sum is \(O_\theta(1)\). Regardless of \(L_n\), the two boundary gaps starting at positions \(i\le n\) contribute at most two to \(\sum_{j\le n}|K_{n,j}-Z_j|\). In particular, the only finite gap that need not occur in the infinite word is the terminal one, so \(\sum_{j\le b}K_{n,j}\le 1+\sum_{j\le b}Z_j\) for every \(b\).

\begin{lemma}[Comparison of short profiles]\label{lem:short-profile-comparison}
Write \(g_{s,L}=[z^s]\exp(\theta\sum_{j>L} z^j/j)\) for the
Ewens coefficient with cycle lengths up to \(L\) excluded.
For some \(c_\theta>0\), uniformly for integers \(1\le L\le c_\theta s\)
and all sufficiently large \(s\),
\[
 g_{s,L}\gg_\theta s^{\theta-1}e^{-\theta H_L}.
\]
In particular, for all sufficiently large \(m\) and
\(1\le L\le c_\theta m/2\), the Ewens probability of a short-cycle
profile through \(L\) with total length at most \(m/2\) is at least
a fixed positive multiple of its Poisson probability.
\end{lemma}

\begin{proof}
Take \(b=\lfloor s/3\rfloor\), \(Q=\lfloor s/(6(1+\theta))\rfloor\),
and \(0<c_\theta<1/(6(1+\theta))\). Then \(L\le Q\) for all
sufficiently large \(s\). For each \(0\le v\le b\), put a cycle
of length \(s-v\) beside an arbitrary profile of total length
\(v\) with lengths in \((L,Q]\). Since \(s-v\ge s-b\ge 2s/3>Q\),
the cycle of length \(s-v\) is unique and has weight at least \(\theta/s\). Multiplication by
\(e^{-\theta(H_Q-H_L)}\) turns the weights on \((L,Q]\) into
independent Poisson probabilities. Summing the distinct profiles
therefore gives
\[
 g_{s,L}\ge \frac{\theta}{s} e^{\theta(H_Q-H_L)}
 \mathbb P\!\left(\sum_{L<j\le Q}j Z_j\le b\right),
 \qquad Z_j\sim\operatorname{Pois}(\theta/j)\ \text{independently}.
\]
The sum has mean \(\theta(Q-L)\le s/6\) and
\(b+1>s/3\), so the displayed probability is at least \(1/2\)
by Markov's inequality. The coefficient bound follows since
\(e^{\theta H_Q}\asymp_\theta s^\theta\).

For a short profile of total length \(u\) in an Ewens permutation
of size \(m\), the ratio of its Ewens probability to its Poisson
probability is \(e^{\theta H_L}g_{m-u,L}/h_m\), by the cycle
profile formula. Now \(m/2\le m-u\le m\) and \(L\le c_\theta(m-u)\).
As soon as \(m/2\) exceeds the lower threshold on \(s\), the coefficient
bound applies to every such \(s=m-u\). Using
\(h_m\asymp_\theta m^{\theta-1}\) gives the claimed comparison.
\end{proof}

\section{Exchange counts}\label{sec:conditional}

Conditioning on a prescribed invariant set separates the permutation
into two independent parts. The reciprocity in
\eqref{s3:existence} then reduces the probability of an invariant
set to the multiplicity of the possible exchanges. The Poisson
version of the same argument will also give a criterion for common sums.

\subsection{Conditioning on a prescribed set}

Fix a \(k\)-set \(A\), and write \(\mathbb P_A\) and \(\mathbb E_A\) for
probability and expectation conditional on \(A\) being invariant. The
Ewens weight \(\theta^{K_n}\) factors over the two restrictions, so
under \(\mathbb P_A\) the permutations on \(A\) and \(A^c\) are
independent with Ewens laws of sizes \(k\) and \(k+d\), where \(d=n-2k\). In particular, \(I_{n,k}\ge1\)
under \(\mathbb P_A\). For an invariant \(k\)-set \(B\), put
\(s=|A\setminus B|\) and \(a=|A\cap B|=k-s\). One can obtain \(B\) by
exchanging \(s\) labels on each side. One can also retain \(a\) labels
of \(A\) and omit \(a+d\) labels of \(A^c\), leaving \(k-a\) labels
of \(A^c\) in \(B\). Each selection is a union of cycles. We use
the exchange description for \(s\le k/2\) and the omission
description for \(s>k/2\) when \(d\le k\).

Let \((Z_j)_{j\ge1}\) and \((Z'_j)_{j\ge1}\) be independent families of
independent Poisson variables of means \(\theta/j\). Set
\[
 U_s=[u^s]\prod_{j=1}^{s}(1+u^j)^{Z_j},\qquad
 U'_s=[u^s]\prod_{j=1}^{s}(1+u^j)^{Z'_j},\qquad U_0=U'_0=1.
\]
The coefficient \(U_s\) counts subsets of cycle copies of total
length \(s\), with copies of the same length treated as distinct.
The two ends of the finite overlap range thus suggest the following
counts in the Poisson model.
\begin{equation}\label{s3:W}
 W_\infty=\sum_{s\ge0}U_sU'_s,\qquad
 W_d=W_\infty+\sum_{a\ge0}U_aU'_{a+d}\quad(d\ge0).
\end{equation}
Both sums in \eqref{s3:W} use the same realization of
\((Z_j)_{j\ge1}\) and \((Z'_j)_{j\ge1}\). A short cycle may therefore
participate in both an unshifted exchange and a shifted exchange. A common positive sum \(s\) contributes
\(U_sU'_s\) exchanges to \(W_\infty\). The \(s=0\) term counts the empty exchange at
\(A\). The \(a=0\) term counts an empty left selection and a right
selection of length \(d\). At \(d=0\) the two sums coincide, giving
\(W_0=2W_\infty\), just as \(B\) and \(B^c\) have the same size in
the finite problem. We allow infinite values and put \(1/\infty=0\).
Whenever \(W_d\) is almost surely finite we write
\(c_d=\E(1/W_d)\in(0,1]\), since \(W_d\ge1\).

If prescribed cells of sizes \(b_1,\ldots,b_\ell\) partition the
labels, the probability that they are all invariant is
\(\bigl(\prod_i b_i!h_{b_i}\bigr)/(n!h_n)\).
This follows by multiplying the Ewens
weights summed over the permutations of each cell. The multinomial
coefficient cancels the factorials when the cells vary, giving
\begin{align}
 \mathbb E I_{n,k}&=\frac{h_k h_{n-k}}{h_n},\label{s3:first}\\
 \mathbb E I_{n,k}^2&=\frac{1}{h_n}\sum_{a=0}^k h_a h_{k-a}^2h_{a+d}.
 \label{s3:second}
\end{align}
For the second identity, two \(k\)-sets with intersection of size \(a\)
give cells of sizes \(a,k-a,k-a,a+d\). Exchangeability also gives
\begin{equation}\label{s3:bias}
 \mathbb P_A(I_{n,k}=j)
 =\frac{j\,\mathbb P(I_{n,k}=j)}{\mathbb E I_{n,k}}\quad(j\ge1).
\end{equation}
Indeed, summing \(\mathbb P(I_{n,k}=j,\ A\text{ invariant})\) over \(A\)
counts each permutation on the event \(I_{n,k}=j\) exactly \(j\) times.
Thus the prescribed-set law weights the unconditioned law by its
multiplicity. Dividing by \(j\) and summing over \(j\ge1\) gives
\begin{equation}\label{s3:existence}
 i_\theta(n,k)=\frac{h_kh_{n-k}}{h_n}\,\mathbb E_A(1/I_{n,k}).
\end{equation}

For a small cyclation check, enumerate the matchings by pairing the
least unpaired endpoint with each available partner and repeating on
the remaining endpoints. For a matching with \(c_j\) loops of length \(j\),
the count is \(I_{n,2}=c_2+\binom{c_1}{2}\).
Among the \(105\) matchings at \(n=4\), the values \(0,2,6\)
occur \(80,24,1\) times. Among the \(945\) matchings at \(n=5\),
the values \(0,1,2,4,10\) occur \(624,240,60,20,1\) times.
Dividing \eqref{s3:second} by \eqref{s3:first} gives the three
contributions indexed by \(s=2-a\) in Table~\ref{tab:cyclation-check}.
For equal sides the last term counts the automatic exchange to the
complement. With one more label on the other side this exchange
requires leaving out a fixed point there, whose expected number is
\(3/5\) in an Ewens permutation of size three at \(\theta=1/2\).

\begin{table}[ht]
\centering
\caption{Exact cyclation checks for \(k=2\). The three overlap
terms sum to \(\E(I_{n,2}^2)/\E I_{n,2}\).}
\label{tab:cyclation-check}
\begin{tabular}{@{}cccc ccc c@{}}
\toprule
\(n\) & \(d\) & \(i_{1/2}(n,2)\) & \(\E I_{n,2}\) & \(s=0\) & \(s=1\) & \(s=2\)
    & \(\E I_{n,2}^2/\E I_{n,2}\) \\
\midrule
4 & 0 & \(5/21\) & \(18/35\) & 1 & \(4/9\) & 1 & \(22/9\) \\
5 & 1 & \(107/315\) & \(10/21\) & 1 & \(2/5\) & \(3/5\) & 2 \\
\bottomrule
\end{tabular}
\end{table}

\subsection{The reciprocal identity and common sums}\label{s3:critical}

Superpose the two profiles, retaining distinguishable copies. Write
\(\mathcal M=(Z_j+Z'_j)_{j\ge1}\) for the merged profile. We use
\(J_N\) and \(T_N\) for the number and total length of copies through
length \(N\) in a single profile, and attach a superscript \((2)\) for
the corresponding quantities after two profiles have been merged. Thus
\[
 J_N^{(2)}=\sum_{j\le N}(Z_j+Z'_j),\qquad
 T_N^{(2)}=\sum_{j\le N}j(Z_j+Z'_j).
\]
Let \(R_N(t)\) count subsets of these copies of total length
\(t\), taking \(R_N(t)=0\) when \(t\) is unattainable. By Poisson
splitting, conditional on \(\mathcal M\), the copies of each length
can be assigned independently to the two sides with equal
probabilities. This realizes the two original profiles. Let \(S_N\)
be the total length on the first side. Write \(V_N\) for the number
of equal exchanges, including the empty exchange, using copies of
individual length at most \(N\). The total of such an exchange may
exceed \(N\). For example, with two copies of length two on each side
and \(N=2\), the exchange totals \(0,2,4\) contribute \(1,4,1\)
respectively. Hence \(V_2=6\), whereas a cutoff of two on the
exchange total gives only \(1+4=5\).

Let \(\Sigma_N^{(2)}=\{t:R_N(t)>0\}\), and write \(\Sigma_N\) for the
attainable totals through length \(N\) in a single Poisson profile.
Use \(\E_\rho\) for expectation over one profile with means \(\rho/j\),
and \(\E_{\theta,\theta}\) for expectation over two independent profiles
with means \(\theta/j\).

\begin{lemma}[Reciprocal exchange identity]\label{lem:reciprocal-exchange}
For a cutoff \(N\) on individual cycle lengths,
\begin{equation}\label{s3:fair}
 V_N=R_N(S_N).
\end{equation}
Conditional on the merged profile,
\begin{equation}\label{s3:reciprocal-support}
 \mathbb E(1/V_N\mid\mathcal M)
     =2^{-J_N^{(2)}}|\Sigma_N^{(2)}|
     \le (T_N^{(2)}+1)2^{-J_N^{(2)}}.
\end{equation}
These two assertions hold also at every almost surely finite cutoff
determined by \(\mathcal M\). For deterministic \(N\), averaging over
the merged profile gives
\begin{equation}\label{eq:exchange-intensity}
 \E_{\theta,\theta}(1/V_N)
   =\E_{2\theta}\!\left[2^{-J_N}|\Sigma_N|\right]
   =e^{-\theta H_N}\E_\theta|\Sigma_N|.
\end{equation}
\end{lemma}

\begin{proof}
Let \(A_1\) consist of the copies marked for the first side through
length \(N\). An equal exchange removes \(E_1\subseteq A_1\) and adds
\(E_2\subseteq A_1^c\) of the same total, where the complement is taken
among the merged copies through \(N\). Its image
\(B=(A_1\setminus E_1)\cup E_2\) has total \(S_N\). Conversely, every
such \(B\) determines \(E_1=A_1\setminus B\) and \(E_2=B\setminus A_1\),
proving \eqref{s3:fair}.

Given \(\mathcal M\), the side marks are independent and fair, so
\(S_N=t\) has probability \(2^{-J_N^{(2)}}R_N(t)\). The factor
\(R_N(t)\) cancels its reciprocal for every attainable \(t\). This
gives \eqref{s3:reciprocal-support}, since \(\Sigma_N^{(2)}\) is
contained in \(\{0,\ldots,T_N^{(2)}\}\). If the cutoff is determined
by \(\mathcal M\), conditioning on \(\mathcal M\) fixes that cutoff and
leaves the marks fair.

For a deterministic cutoff, the merged counts have means \(2\theta/j\).
A count of \(a\) copies of length \(j\) contributes the weighted mass
\[
 2^{-a}e^{-2\theta/j}\frac{(2\theta/j)^a}{a!}
 =e^{-\theta/j}\,e^{-\theta/j}\frac{(\theta/j)^a}{a!}.
\]
Multiplication over \(1\le j\le N\) gives the factor
\(e^{-\theta H_N}\) and changes the means to \(\theta/j\), proving
\eqref{eq:exchange-intensity}.
\end{proof}

The empty exchange gives \(V_N\ge1\). Every finite exchange appears as
soon as its largest copy is admitted, whence \(V_N\uparrow W_\infty\).
For the merged lengths \(1,1,2\), the polynomial
\((1+x)^2(1+x^2)=1+2x+2x^2+2x^3+x^4\) gives eight assignments and
five attainable totals. Thus the conditional mean of \(1/V_N\) is
\(5/8\). At \(N=1\) a single profile has \(1+Z_1\) attainable totals,
so \eqref{eq:exchange-intensity} gives \(e^{-\theta}(1+\theta)\).

Green and Sawhney \cite[Eq.~(1.5), Definition~6.1, and Lemma~6.2]{gs}
study the normalized support and bound it along walk excursions.
In a profile of intensity \(\rho/j\), the expected number of copies
in \((2^{v-1},2^v]\) tends to \(\rho\log 2\). For the merged pair this
limit equals one at \(\theta=\theta_c\), where the merged counts have
the intensities of the lower process of Green and Sawhney. For a single
profile the limit equals one at \(\theta=\lambda_0\).

A selection of total \(S_N+d\) corresponds to a toggle in which the
copies added from the second side exceed those removed from the
first side by total length \(d\). We will also use the count
\begin{equation}\label{eq:shifted-exchange}
 V_{N,d}=R_N(S_N)+R_N(S_N+d).
\end{equation}
For \(d=0\) the two terms agree. The first term always counts
the empty toggle.

\begin{theorem}[Common sums at the critical parameter]\label{s3:poisson-critical}
Let \(\theta\ge\theta_c=(2\log 2)^{-1}\), and take two independent
Poisson profiles with means \(\theta/j\). Almost surely their sumsets
have infinitely many common positive sums. The same holds
simultaneously for every fixed lower cutoff when both profiles are
restricted to lengths greater than that cutoff. In particular,
\(W_\infty=\infty\) almost surely.
\end{theorem}

\begin{proof}
Set \(\theta=\theta_c\). Since \(\mathbb E T_{2^\ell}^{(2)}=2^\ell/\log 2\),
Markov's inequality and the Borel--Cantelli lemma give
\(T_{2^\ell}^{(2)}\le 2^\ell \ell^2\) eventually almost surely. The block counts \(B_v=J_{2^v}^{(2)}-J_{2^{v-1}}^{(2)}\), \(v\ge1\), are
independent Poisson variables of means
\[
 \mu_v=\frac{H_{2^v}-H_{2^{v-1}}}{\log 2}=1+O(2^{-v}).
\]
Since \(1/x\) is strictly decreasing, the harmonic sum in the
numerator is less than
\(\int_{2^{v-1}}^{2^v} dx/x=\log 2\), so \(\mu_v<1\). Add an
independent Poisson variable \(E_v\) of mean \(1-\mu_v\) in each
block to obtain independent \(\operatorname{Pois}(1)\) variables
\(Y_v\). Then
\[
 J_{2^\ell}^{(2)}-\ell=J_1^{(2)}+\sum_{v=1}^\ell(Y_v-1)-\sum_{v=1}^\ell E_v .
\]
Since \(\sum_v\E E_v<\infty\) and each \(E_v\) is a
nonnegative integer, almost surely only finitely many \(E_v\) are nonzero.
Thus \(J_1^{(2)}-\sum_{v=1}^\ell E_v\) eventually stays constant. The centered increments \(Y_v-1\) have variance one.
By the law of the iterated logarithm, their partial sums exceed
\(\tfrac12\sqrt{2\ell\log\log \ell}\) infinitely often almost surely.
This grows faster than \(3\log_2 \ell\), even after adding the
finite random constant. Hence \(J_{2^\ell}^{(2)}-\ell\ge3\log_2 \ell\)
infinitely often almost surely, and along these indices
\((T_{2^\ell}^{(2)}+1)2^{-J_{2^\ell}^{(2)}}\le \ell^{-1}+2^{-\ell}\ell^{-3}\).
This assertion depends only on the merged profile and so holds
under its marginal law. With \(1/\infty=0\),
monotonicity of \(V_N\) and \eqref{s3:reciprocal-support} give
\[
 0\le \E(1/W_\infty\mid\mathcal M)
 \le \E(1/V_N\mid\mathcal M)
 \le (T_N^{(2)}+1)2^{-J_N^{(2)}}
\]
almost surely, simultaneously at every dyadic cutoff by countability.
Taking the infimum gives \(W_\infty=\infty\) almost surely.
For each fixed \(a\), the merged count \(J_a^{(2)}\) is almost surely
finite, so the term \(U_aU'_a\) in \eqref{s3:W} is finite. This holds
simultaneously for all \(a\). Hence there are infinitely many distinct
common positive totals.

For an integer lower cutoff \(L_0\ge0\), the same argument applies with
\(J_N^{(2)}-J_{L_0}^{(2)}\) and \(T_N^{(2)}-T_{L_0}^{(2)}\) in the bound, taking \(J_0^{(2)}=T_0^{(2)}=0\).
For \(N>L_0\) the new upper bound is at most \(2^{J_{L_0}^{(2)}}\) times the
old one. Countability gives the conclusion simultaneously for all
integer cutoffs \(L_0\). A real cutoff is equivalent to its integer part.
For \(\theta>\theta_c\), add independent Poisson counts of
means \((\theta-\theta_c)/j\) to both profiles. This can only enlarge
the sumsets and preserves all common sums after every lower cutoff.
\end{proof}

\subsection{The conditional law and marked coefficients}

\begin{theorem}\label{s3:main}
Let \(1\le k\le n/2\), put \(d=n-2k\), and condition on a prescribed
\(k\)-set \(A\) being invariant. Use \(\lambda_0=1/\log 2\),
\(\theta_c=\lambda_0/2\), and \(D(\theta)\) as in
\eqref{eq:rate-function}. If \(0<\theta<\theta_c\), set
\begin{equation}\label{s3:beta}
 \beta=
 \begin{cases}
  1-2\theta,&0<\theta\le\lambda_0/4,\\
  D(\theta),&\lambda_0/4<\theta<\theta_c,
 \end{cases}
 \qquad \eta=\frac{\beta}{1+\beta}.
\end{equation}
Each \(W_d\), for \(d\in\{0,1,2,\ldots\}\cup\{\infty\}\), is almost
surely finite. With \(c_d=\mathbb E(1/W_d)\),
\begin{align}
 d_{\mathrm{TV}}\bigl(\mathcal L_{\mathbb P_A}(I_{n,k}),\mathcal L(W_d)\bigr)
  &=O_\theta(k^{-\eta}),\label{s3:TV}\\
 i_\theta(n,k)&=\frac{h_kh_{n-k}}{h_n}
       \bigl(c_d+O_\theta(k^{-\eta})\bigr),\label{s3:asymptotic}
\end{align}
uniformly in \(1\le k\le n/2\). In particular the approximating law
uses the current \(d=n-2k\), even when \(d\) oscillates with \(k\).
The constants satisfy \(c_0=c_\infty/2\)
and \(0<c_\infty-c_d=O_\theta((1+d)^{-\beta})\) for every finite \(d\).
They are bounded below by a positive constant depending only on
\(\theta\). In particular, \(i_\theta(n,k)\asymp_\theta k^{\theta-1}\)
uniformly in this range.

If \(\theta\ge\theta_c\), then \(W_\infty=\infty\) almost surely, and
\[
 i_\theta(n,k)=o_\theta(h_kh_{n-k}/h_n)
 \qquad(k\to\infty)
\]
uniformly over \(n\ge2k\).
\end{theorem}

Since \(D'(\theta)=2-\lambda_0/\theta\), \(D\) has its unique minimum
\(D(\theta_c)=0\). The two expressions for \(\beta\) agree in value and
derivative at \(\lambda_0/4\), with derivative \(-2\).
For \(\theta=\theta_c\), Proposition~\ref{s3:critical-rate} gives the
order \((1+\log k)^{-1/2}\) for the first-moment ratio.
Proposition~\ref{prop:large-theta} gives the fixed-set probabilities
for \(\theta>\theta_c\).
To control large exchanges in the subcritical range, we mark short
cycles and bound a fractional moment of each overlap count.

For integers \(L\ge1\) and \(t\ge0\), use the coefficients
\(b_t(v,L)\) to give each cycle of length at most \(L\) a factor
\(v\ge0\). Put
\[
 B(z;v,L)=(1-z)^{-\theta}
       \exp\left(\theta(v-1)\sum_{j\le L}\frac{z^j}{j}\right),
 \qquad b_t(v,L)=[z^t]B(z;v,L).
\]
Write \(\theta_j=\theta v\) for \(j\le L\) and \(\theta_j=\theta\)
otherwise. Since \(\log B(z;v,L)=\sum_{j\ge1}\theta_jz^j/j\),
all the coefficients \(b_t(v,L)\) are nonnegative, including
\(b_t(0,L)=g_{t,L}\). Logarithmic differentiation gives,
for \(t\ge1\),
\[
 t b_t(v,L)=\sum_{j=1}^{t}\theta_j b_{t-j}(v,L)
 \le \theta\max(1,v)\sum_{i=0}^{t-1}b_i(v,L).
\]
The last sum is at most \(e B(e^{-1/t};v,L)\), since
\(1\le e^{1-i/t}\) for \(0\le i<t\). Moreover,
\[
 \sum_{j\le L}\frac{e^{-j/t}}{j}
   =\log(1+\min(t,L))+O(1).
\]
For \(j\le\min(L,t)\), replacing \(e^{-j/t}\) by one costs at
most \(1/t\) per term in the sum. If \(L>t\), bound the sum over
\(t<j\le L\) by \(t^{-1}\sum_{j>t}e^{-j/t}\).
The error in this harmonic-sum estimate is bounded after
multiplication by \(\theta(v-1)\), including at \(v=0\).
It follows, with \(t=0\) immediate, that for \(v\ge0\)
\begin{equation}\label{s3:coef}
 b_t(v,L)\ll_{\theta,v}(t+1)^{\theta-1}
                    (1+\min(t,L))^{\theta(v-1)} .
\end{equation}

\begin{lemma}[Marked fixed-set counts]\label{lem:marked-fixed-sets}
Let \(K_{\le k}\) count the cycles of length at most \(k\) in an
unconditioned Ewens permutation. For \(y,v,\theta>0\) and
\(1\le k\le n/2\),
\begin{align}
 \mathbb E y^{K_{\le k}}
 &= \frac{b_n(y,k)}{h_n}\ll_{\theta,y} k^{\theta(y-1)}, \notag\\
 \mathbb E(I_{n,k}v^{K_{\le k}})
 &= \frac{h_k(v\theta)b_{n-k}(v,k)}{h_n}
    \ll_{\theta,v} k^{2\theta v-\theta-1}.
 \label{s3:fixed-mark}
\end{align}
\end{lemma}

\begin{proof}
Marking the cycles through length \(k\) gives the first identity.
For the second, sum over prescribed \(k\)-sets and split at each
invariant set. Every cycle on the first side has length at most \(k\),
so its coefficient is \(h_k(v\theta)\). On the other side the
coefficient is \(b_{n-k}(v,k)\). The two bounds follow from
\eqref{s3:coef} and \eqref{s2:gamma}, since \(n-k\ge k\) and \(n-k\asymp n\).
\end{proof}

\subsection{Marked estimates for the overlap ranges}\label{sec:overlap}

At each scale in Table~\ref{tab:overlap}, we use one estimate for the number of competing sets with a mark on short cycles,
and a second estimate for the mark alone. Their constants are uniform in \(k,n,d\)
and in the scale. Dependence on the exponents is indicated below.

Let \(K^{\mathrm{tot}}_{\le L}\) count cycles of length at most \(L\) on the two sides
together. Independence gives
\(\mathbb E_A y^{K^{\mathrm{tot}}_{\le L}}=b_k(y,L)b_{k+d}(y,L)/(h_kh_{k+d})\).
If \(m\ge1\), \(L\le4m\), and \(y>1\), then \eqref{s3:coef} gives
\begin{equation}\label{s3:short-mark}
 \E_A y^{K^{\mathrm{tot}}_{\le L}}\ll_{\theta,y} m^{2\theta(y-1)}.
\end{equation}

Write \(N_{u,w}\) for the number of pairs of invariant subsets of
sizes \(u\) in \(A\) and \(w\) in \(A^c\). For
\(0\le u\le\min(k,L)\) and \(0\le w\le\min(k+d,L)\), every
cycle in a selected part is marked by \(v\). Only cycles of length at
most \(L\) are marked in the residual parts, which gives
\begin{equation}\label{s3:split-mark}
 \mathbb E_A(v^{K^{\mathrm{tot}}_{\le L}}N_{u,w})
 =h_u(v\theta)h_w(v\theta)
   \frac{b_{k-u}(v,L)b_{k+d-w}(v,L)}{h_kh_{k+d}}.
\end{equation}
The rows of Table~\ref{tab:overlap} partition the overlap indices
\((s,a)\), where \(a=k-s\), with \(s\ge1\) and \((a,d)\ne(0,0)\).
The omitted indices give \(A\) and, when \(d=0\), \(A^c\).
Rows~2 and~3 use the \(a\) labels kept in \(A\) and the \(a+d\)
omitted from \(A^c\), since these are the small parts near the
complement. When \(d>k\), a part of size \(a+d\) is too large for
this purpose, so row~4 uses the direct exchange. The first two
rows subdivide their index ranges at dyadic scales \(m\).
Here a dyadic scale is \(m=2^j\), \(j\ge0\). The ranges in the
first two rows include the displayed restrictions on \(s\) and \(a\).
Row~3 includes \(a=0\) for \(0<d\le k\), with \(d\le a+d<2d\).
At most one of rows~3 and 4 contributes an additional scale.

Let \(C\) count the competing invariant \(k\)-sets in one row
at the indicated scale.

\begin{table}[!ht]
\centering
\caption{Overlap ranges and their cycle-length cutoffs. The parts
of sizes \(u,w\) are exchanged in rows 1 and 4. In rows 2 and 3,
\(u\) is retained on the left and \(w\) is omitted on the right.}
\label{tab:overlap}
\[
\begin{array}{c@{\quad}l@{\qquad}c@{\qquad}c@{\qquad}c}
 &\text{Range}&\text{Scale }m&(u,w)&L\\
1&1\le s\le k/2,\ m\le s<2m&\text{dyadic}&(s,s)&4m\\
2&s> k/2,\ a\ge\max(1,d),\ m\le a<2m&\text{dyadic}&(a,a+d)&4m\\
3&s> k/2,\ 0\le a<d\le k&d&(a,a+d)&2d\\
4&s> k/2,\ d>k&k&(s,s)&2k
\end{array}
\]
\end{table}

For \(0<v\le1\), the marked count on each nonempty row and scale satisfies
\begin{equation}\label{s3:overlap-mark}
 \mathbb E_A(v^{K^{\mathrm{tot}}_{\le L}}C)\ll_{\theta,v}m^{4\theta v-2\theta-1}.
\end{equation}
In row 1 the residual sizes \(k-s\) and \(k+d-s\) are comparable to
their original side sizes and are at least \(m\). The ratio in
\eqref{s3:split-mark} is \(O_{\theta,v}(m^{2\theta(v-1)})\) by
\eqref{s3:coef}. Summing the two selected-part coefficients gives
\(O_{\theta,v}(m^{2v\theta-1})\). In row 2 both residual sizes equal
\(k-a>k/2\). The original side sizes are \(k\) and \(k+d\le 3k/2\)
since \(d\le a<k/2\), and both selected sizes are comparable to \(m\). The same bounds
apply. In row 3 the residual sizes again equal \(k-a>k/2\ge d/2\)
and \(k+d\le2k\), so the residual ratio has the same bound at \(m=d\).
The row has \(0\le a<\min(d,k/2)\). Extending the selected-size sum
to \(0\le a<d\) gives
\[
 \sum_{a<d}h_a(v\theta)h_{a+d}(v\theta)
 \ll_{\theta,v}d^{v\theta-1}\sum_{a<d}h_a(v\theta)
 \ll_{\theta,v}d^{2v\theta-1}.
\]
In row 4 the residual sizes are \(a=k-s<k/2\) and \(a+d\).
Since \(a\le L=2k\), \(b_a(v,L)=h_a(v\theta)\). The second residual
size is comparable to \(k+d\) and at least \(k\). Thus in
\(b_{a+d}(v,2k)/h_{k+d}\) the unmarked factor
\((a+d+1)^{\theta-1}\) is comparable to \(h_{k+d}\), and the mark
contributes \(O_{\theta,v}(k^{\theta(v-1)})\). With \(s\asymp k\)
and \(h_k^{-1}\asymp_\theta k^{1-\theta}\), the sum in
\eqref{s3:split-mark} is bounded by a constant times
\[
 k^{2v\theta-2}k^{1-\theta}k^{\theta(v-1)}
          \sum_{a<k/2}h_a(v\theta)
 \ll_{\theta,v} k^{4v\theta-2\theta-1}.
\]
This proves \eqref{s3:overlap-mark}.

Put \(\Phi_\theta(r)=2\theta(2^r-1)-r\). For \(0<r<1\), take
\(v=2^{r-1}\in(1/2,1)\) and \(y=v^{-r/(1-r)}=2^r\in(1,2)\).
Since
\(v^r y^{1-r}=1\), the product of the two marked factors, raised
to powers \(r\) and \(1-r\), is \(C^r\). By H\"older's inequality,
\begin{equation}\label{s3:fractional-shell}
 \mathbb E_A C^r
 \le \bigl(\mathbb E_A(v^{K^{\mathrm{tot}}_{\le L}}C)\bigr)^r
      \bigl(\mathbb E_A y^{K^{\mathrm{tot}}_{\le L}}\bigr)^{1-r}
 \ll_{\theta,r}m^{\Phi_\theta(r)}.
\end{equation}
Indeed, \(4\theta v-2\theta=2\theta(y-1)\), so the exponents in
\eqref{s3:overlap-mark} and \eqref{s3:short-mark} combine to give
\(r[2\theta(y-1)-1]+(1-r)2\theta(y-1)=\Phi_\theta(r)\).
The last bound holds at \(r=1\) by \eqref{s3:overlap-mark} at \(v=1\).
Since \(\mathbb P_A(C>0)\le\mathbb E_A C^r\), minimizing the exponent
over \(0<r\le1\) gives
\begin{equation}\label{s3:overlap-tail}
 \mathbb P_A(C>0)\ll_\theta m^{-\beta}
 \qquad(0<\theta<\theta_c).
\end{equation}
The derivative is \(\Phi_\theta'(r)=2\theta 2^r\log 2-1\).
When \(\theta\le\lambda_0/4\), its zero is at or beyond one and the
constrained minimum occurs at \(r=1\). Otherwise the minimizer is
\(r=\log_2(\lambda_0/(2\theta))\). The minimum is \(-\beta\)
in \eqref{s3:beta}. At \(\theta=1/2\) the number of merged Poisson copies through length
\(4m\) has mean \(2\theta H_{4m}=\log m+O(1)\). The optimizing
cycle-count level is \(\lambda_0\log m\). The corresponding
rate is
\(D(1/2)=1-\lambda_0+\lambda_0\log\lambda_0=D(1)\), the exponent
\(1-(1+\log\log 2)/\log 2\) of \cite{efg,gs} at \(\theta=1\). Both rates compare a cycle intensity of one with the
level at which \(2^{\lambda_0\log m}=m\). Here
\(\beta=0.086071\ldots\) bounds the probability of a large exchange.
Balancing this tail estimate with the coupling error below gives
\(\eta=\beta/(1+\beta)=0.079250\ldots\) in \eqref{s3:TV}. Summing
\eqref{s3:overlap-tail} over the dyadic rows and the one extra scale,
when present, bounds the probability of any contribution at scales
at least \(M\) by \(O_\theta(M^{-\beta})\).

\subsection{The subcritical law}\label{s3:subcritical}

Fix a finite imbalance \(d\) and take \(n=2k+d\). For fixed positive
\(s\), the exchange of total \(s\) is eventually in row~1 of
Table~\ref{tab:overlap}. For fixed \(a\) with \((a,d)\ne(0,0)\), the
shifted exchange of totals \(a,a+d\) is eventually in row~2 if
\(a\ge\max(1,d)\) and in row~3 if \(a<d\). Assign the corresponding
scale to each Poisson term. The empty exchange and, when \(d=0\),
the automatic complement require no tail estimate.

Let \(\mathcal T_d(M)\) be the event that at least one Poisson term
with assigned scale at least \(M\) is positive. For \(d=\infty\)
retain only row~1. We claim that
\begin{equation*}
 \sup_{d\in\{0,1,\ldots\}\cup\{\infty\}}
          \Pp(\mathcal T_d(M))\ll_\theta M^{-\beta},\qquad M\ge1 .
\end{equation*}
Indeed, take any finite collection of the indicated terms and let
\(L\) be their largest selection total. In finite degree the sum of
\eqref{s3:overlap-tail} over scales at least \(M\) bounds the
probability of a contribution by \(O_\theta(M^{-\beta})\).
For \(k\ge2L\) sufficiently large, all terms in the collection fall
in their indicated rows. The two Feller couplings through \(L\)
transfer this bound to the Poisson profiles with error \(O_\theta(L/k)\).
Let \(k\to\infty\) first, then exhaust the terms by finite collections.
The finite-degree bound is independent of both the collection and
\(d\), proving the claim. For the unshifted terms alone the same proof
uses row~1 at, for example, \(n=2k\).

Consequently the contributing scales are almost surely bounded. Each
bounded range of scales contains only finitely many indices, and each
term is finite almost surely. This proves that every \(W_d\) is
almost surely finite. When \(d\ge1\), the shifted terms have scale
\(d\) if \(a<d\), and scale greater than \(a/2\ge d/2\) otherwise. Hence
\begin{equation}\label{s3:d-tail}
 \mathbb P(W_d\ne W_\infty)\ll_\theta(1+d)^{-\beta}.
\end{equation}

For \(k\ge4\) and an integer \(1\le M\le k/4\), let
\(I_{n,k}^{(M)}\) count the sets with \(s\le M\) or \(a+d\le M\).
These ranges are disjoint. Their Poisson counterparts are
\[
 W_\infty^{(M)}=\sum_{s=0}^M U_sU'_s,\qquad
 W_d^{(M)}=W_\infty^{(M)}+\sum_{a=0}^{M-d}U_aU'_{a+d},
\]
with an empty second sum for \(d>M\). Here \(M\) bounds each selection's total length, while \(N\) in
\(V_{N,d}\) bounds each individual copy length.

A discarded exchange in row~1 has \(s>M\) and therefore scale
\(m>M/2\). A discarded shifted exchange has \(a+d>M\). In row~2,
\(a+d\le2a<4m\) gives \(m>M/4\). In row~3, \(a+d<2d=2m\) gives
\(m>M/2\). In particular, if \(d>M\), the entire shifted Poisson
sum is omitted. Its terms with \(a<d\) have scale \(d>M\), and
those with \(a\ge d\) have scale greater than \(a/2\ge d/2>M/2\).
For the finite count, row~4 uses a direct exchange of scale
\(k\ge4M\). Summing the finite-degree tail
bounds and using the preceding Poisson tail bound gives
\begin{equation}\label{s3:trunc-tail}
 \mathbb P_A(I_{n,k}^{(M)}\ne I_{n,k})
   +\mathbb P(W_d^{(M)}\ne W_d)\ll_\theta M^{-\beta}.
\end{equation}
Each retained count is the same function of cycle counts through
length \(M\). Since \(M\le k/4\le (k+d)/4\),
\eqref{s2:feller} applies on both sides. Under the independent
couplings the retained counts differ with probability at most
\(O_\theta(M/k+M/(k+d))=O_\theta(M/k)\). The error is therefore
\(O_\theta(M^{-\beta}+M/k)\). Take \(M=\lfloor k^{1/(1+\beta)}\rfloor\), which satisfies
\(1\le M\le k/4\) for sufficiently large \(k\). Both
\(M^{-\beta}\) and \(M/k\) are \(O_\theta(k^{-\eta})\), giving
\eqref{s3:TV}. Bounded \(k\) are absorbed into the implied constant. Equation
\eqref{s3:existence} and the bounded test function \(j\mapsto1/j\)
give \eqref{s3:asymptotic}. The same test function and
\eqref{s3:trunc-tail} give
\[
 0\le \E(1/W_d^{(M)})-c_d\ll_\theta M^{-\beta}.
\]

Since \(W_0=2W_\infty\), we have \(c_0=c_\infty/2\). For \(d\ge1\)
the event \(Z'_d\ge1\) makes the shifted sum positive and has positive
probability, so \(c_d<c_\infty\). Since \(W_d\ge W_\infty\ge1\),
\[
 0\le W_\infty^{-1}-W_d^{-1}\le \ind_{\{W_d\ne W_\infty\}}.
\]
Together with \eqref{s3:d-tail} this gives the bound on
\(c_\infty-c_d\). Each \(c_d>0\) by finiteness, and convergence to
\(c_\infty>0\) gives a uniform lower
bound. By \eqref{s2:gamma}, \(h_kh_{n-k}/h_n\asymp_\theta k^{\theta-1}\).
For bounded \(k\), \eqref{s3:bias} identifies
\(\mathbb E_A I_{n,k}\) with the ratio of \eqref{s3:second} to
\eqref{s3:first}, which is bounded uniformly in \(d\). Indeed,
\(0<\theta<\theta_c<1\), so \eqref{s2:gamma} gives
\[
 \frac{h_{a+d}}{h_{k+d}}
 \ll_\theta\left(\frac{k+d+1}{a+d+1}\right)^{1-\theta}
 \le (k+1)^{1-\theta}\qquad(0\le a\le k).
\]
Thus \(\mathbb E_A(1/I_{n,k})\ge1/\mathbb E_A I_{n,k}\) is
bounded below uniformly in \(d\) for these \(k\) as well.

\subsection{Existence at and above the threshold}\label{s3:supercritical}

For each fixed \(M\), the count of exchanges of size at most \(M\)
is a lower bound for \(I_{n,k}\) under \(\mathbb P_A\). For
\(k\ge2M\), coupling the two profiles through length \(M\) as
above gives \(W_\infty^{(M)}\) with error \(O_\theta(M/k)\)
uniformly over \(n\ge2k\). Consequently
\[
 \limsup_{k\to\infty}\sup_{n\ge2k}\mathbb E_A(1/I_{n,k})
       \le \mathbb E(1/W_\infty^{(M)}).
\]
We have taken \(k\to\infty\) with \(M\) fixed, after taking
the supremum over \(n\ge2k\). Now \(M\to\infty\) gives zero on the
right by bounded convergence and Theorem~\ref{s3:poisson-critical}.
Equation \eqref{s3:existence} proves the little-oh assertion.

The same distinction appears in the expected number of attainable
totals of a single Poisson profile. Since \(V_N\uparrow W_\infty\) and
\(0<1/V_N\le1\), bounded convergence in \eqref{eq:exchange-intensity} gives
\[
 \lim_{N\to\infty} e^{-\theta H_N}\E_\theta|\Sigma_N|
   =
   \begin{cases}
      c_\infty,&0<\theta<\theta_c,\\
      0,&\theta\ge\theta_c.
   \end{cases}
\]
In particular, below \(\theta_c\) the expected support is asymptotic
to \(c_\infty e^{\theta H_N}\).

\subsection{Existence-conditioned laws and graph applications}

For a probability law \(\nu\) on the positive integers, define
\[
 \mathcal R(\nu)\{j\}=\frac{\nu\{j\}/j}{\sum_{\ell\ge1}\nu\{\ell\}/\ell}.
\]
This transformation reverses the size weighting in \eqref{s3:bias}.

\begin{corollary}\label{s3:exist-cor}
Let \(d=n-2k\) and \(0<\theta<\theta_c\). The law of \(I_{n,k}\)
conditional on \(I_{n,k}>0\) is within \(O_\theta(k^{-\eta})\) in
total variation of \(\mathcal R(\mathcal L(W_d))\), uniformly in
\(1\le k\le n/2\).
For finite \(d>0\) the approximating law has an atom at one, and
for \(d=0\) it has an atom at two. The law \(\mathcal R(\mathcal L(W_\infty))\)
also has an atom at one. These atom masses are bounded below by a
positive constant depending only on \(\theta\).
\end{corollary}

\begin{proof}
Equations \eqref{s3:bias} and \eqref{s3:existence} express the law in
question as \(\mathcal R(\mathcal L_{\mathbb P_A}(I_{n,k}))\).
Multiplication by \(1/j\) contracts the sum of the absolute mass
differences. The two normalizing constants are
\(\mathbb E_A(1/I_{n,k})\) and \(c_d\), both bounded below uniformly
by Theorem~\ref{s3:main} and its proof. Dividing by them preserves
the order of the total variation error.

Take \(M\) large enough that the omitted Poisson contribution in
\eqref{s3:trunc-tail} has probability at most \(1/2\), uniformly in
\(d\). The event \(Z_j=Z'_j=0\) for \(j\le M\) has probability
\(e^{-2\theta H_M}\). The counts beyond \(M\) are independent of those through \(M\).
Setting the latter to zero can only remove omitted contributions.
Thus conditioning on their absence can only decrease the probability
of any omitted contribution. On this event the retained count is
two if \(d=0\) and one otherwise. With conditional probability at least \(1/2\) it is also
the full count. The same argument with just the unshifted terms
applies to \(W_\infty\). After reversing the size weighting, the
atom has mass at least \(e^{-2\theta H_M}/4\), since the relevant
normalizing constant is at most one.
\end{proof}

\begin{proposition}\label{prop:constant-support}
If \(0<\theta<\theta_c\), then for each finite \(d\ge0\),
\begin{equation}\label{eq:constant-support}
 c_d=\tfrac12 c_\infty+
     \tfrac12\,\E\!\left(\frac{\ind_{\{W_d=W_\infty\}}}{W_\infty}\right).
\end{equation}
Moreover, \(c_d-c_0\) is bounded below by a positive constant
depending only on \(\theta\), uniformly over \(d\ge1\).
\end{proposition}

\begin{proof}
By \eqref{s3:fair} and \eqref{eq:shifted-exchange},
the two terms in \(V_{N,d}\) are \(R_N(S_N)\) and \(R_N(S_N+d)\).
Given the merged profile,
\[
 \mathbb E(1/V_{N,d}\mid\mathcal M)
 =2^{-J_N^{(2)}}\sum_{t\in\mathbb Z}
       \frac{R_N(t)}{R_N(t)+R_N(t+d)}
 =2^{-J_N^{(2)}-1}|\Sigma_N^{(2)}\cup(\Sigma_N^{(2)}-d)|,
\]
where a zero denominator contributes zero. For the second equality
use \(R_N(t)=R_N(T_N^{(2)}-t)\). If \(u=R_N(t)\) and \(v=R_N(t+d)\),
reflection \(t\mapsto T_N^{(2)}-d-t\) interchanges \(u\) and \(v\).
For \(u+v>0\) the two reflected summands have sum
\(u/(u+v)+v/(u+v)=1\). A fixed point with positive denominator contributes \(1/2\).
The indices with \(u+v>0\) form \(\Sigma_N^{(2)}\cup(\Sigma_N^{(2)}-d)\).
Its cardinality is \(|\Sigma_N^{(2)}|+|\Sigma_N^{(2)}\setminus(\Sigma_N^{(2)}-d)|\),
since translation preserves cardinality. Divided by \(2^{J_N^{(2)}}\),
the two terms give, respectively, the conditional mean of \(1/V_N\)
and that of \(\ind_{\{R_N(S_N+d)=0\}}/V_N\).
As \(N\) increases, \(R_N(S_N+d)\) increases to the shifted sum
in \eqref{s3:W}, since it counts shifted selections whose individual
cycle lengths are at most \(N\). Consequently its zero indicator
decreases to \(\ind_{\{W_d=W_\infty\}}\). Bounded convergence in the
reciprocal identity gives \eqref{eq:constant-support}. If \(d>0\),
the event \(W_d=1\) lies in \(\{W_d=W_\infty\}\). Its probability is
bounded below uniformly in \(d\) by the proof of
Corollary~\ref{s3:exist-cor}.
\end{proof}

For cyclations, Corollary~\ref{cor:cyclation-consequences} gives
\(\sqrt{k}\,i_{1/2}(n,k)\to c_\infty/\sqrt{\pi(1-\alpha)}\) when
\(k\to\infty\), \(k/n\to\alpha\), and \(d\to\infty\).
For fixed \(d\) the limit is \(\sqrt{2/\pi}\,c_d\).
At equal halves complementation pairs each invariant set with its
complement. Thus \(I_{2k,k}/2\) counts unordered partitions into two
unions of loops, and \(I_{2k,k}=2\) means that the partition is unique.

Take \(\theta=1/2\) and put
\[
 h_m^{[3]}=[z^m]\exp\left(\frac12\sum_{j\ge3}\frac{z^j}{j}\right).
\]
Write \(I_{n,k}^{[3]}\) for the number of unions of components
having total size \(k\) in a uniformly chosen simple \(2\)-regular
graph on the labeled vertex set \([n]\). Define \(W_d^{[3]}\) by
\eqref{s3:W} with \(Z_1=Z'_1=Z_2=Z'_2=0\) and the other means
unchanged, and put \(c_d^{[3]}=\mathbb E(1/W_d^{[3]})\).

\begin{corollary}\label{s3:simple-cor}
Let \(3\le k\le n/2\), put \(d=n-2k\), and fix a vertex \(k\)-set
\(A\). As \(k\to\infty\), uniformly over \(n\ge2k\),
\[
 d_{\mathrm{TV}}\bigl(\mathcal L(I_{n,k}^{[3]}\mid
 A\text{ is a union of components}),\mathcal L(W_d^{[3]})\bigr)
   =O(k^{-\eta}),
\]
where \(\eta\) is from \eqref{s3:beta} at \(\theta=1/2\). In particular,
\[
 \mathbb P(I_{n,k}^{[3]}>0)=
 \frac{h_k^{[3]}h_{n-k}^{[3]}}{h_n^{[3]}}
        \bigl(c_d^{[3]}+O(k^{-\eta})\bigr),
 \qquad c_0^{[3]}=\tfrac12c_\infty^{[3]}.
\]
The conditional law given \(I_{n,k}^{[3]}>0\) is within
\(O(k^{-\eta})\) in total variation of
\(\mathcal R(\mathcal L(W_d^{[3]}))\).
\end{corollary}

\begin{proof}
Each simple 2-regular graph arises from \(2^n\) endpoint pairings, by
assigning the two edges incident to each vertex to its endpoints.
The matching profile conditioned on simplicity is thus the Ewens
profile at \(\theta=1/2\) with lengths one and two absent. Given the
prescribed union \(A\), simplicity is the joint event \(E\) that both
independent side profiles avoid these lengths. Let \(F\) be the
analogous Poisson event. In the coupling used for \eqref{s3:TV},
the choice \(M=\lfloor k^{1/(1+\beta)}\rfloor\) satisfies
\(2\le M\le k/4\) for sufficiently large \(k\), so the required
window is contained in the coupling on both sides. The indicators of
\(E,F\) agree whenever the coupled profiles agree through \(M\).
Together with \eqref{s3:trunc-tail}, this shows that the
multiplicities and the indicators agree jointly off an event of
probability \(O(k^{-\eta})\). The probability of \(F\) is
\(\exp[-2(1/2+1/4)]=e^{-3/2}\). The two joint probabilities of any
test event on the counts with \(E,F\), as well as the two
normalizing probabilities, therefore differ by at most the
coupling error. Dividing by the normalizing probabilities, which
are bounded away from zero, gives the conditional estimate. The
Poisson count under \(F\) has law \(\mathcal L(W_d^{[3]})\).

The restriction \(k\ge3\) ensures that both sides can be completed
in the simple model, since each can be a single cycle. Splitting
at \(A\) with the restricted weights gives the exact first
moment \(h_k^{[3]}h_{n-k}^{[3]}/h_n^{[3]}\). The size-bias identity
\eqref{s3:bias} holds by exchangeability. Setting the short counts to zero gives \(W_d^{[3]}\le W_d\), hence
\(c_d^{[3]}\ge c_d\) uniformly in \(d\). Applying the conditional
estimate to \(j\mapsto1/j\) bounds the corresponding finite
normalizing constants away from zero for all sufficiently large \(k\).
Reversing the size weighting now gives the existence probability
and the law conditional on existence. The identity at \(d=0\) follows
from the two identical sums. The coupling through length two also gives
\(h_m^{[3]}/h_m=e^{-3/4}+O(m^{-1})\), as used in
Corollary~\ref{cor:cyclation-consequences}.

\end{proof}

Given that a prescribed half \(A\) is a union of components, the
graphs on \(A\) and \(A^c\) are independent uniform simple
\(2\)-regular graphs on \(k\) vertices. If \(A_{i,s}\) counts the
\(s\)-vertex unions of components in graph \(i\), the count of
invariant halves is
\[
 T=\sum_{s=0}^{k} A_{1,s}A_{2,s}.
\]
Indeed, removing a union of size \(s\) from \(A\) and adding one of
the same size from \(A^c\) gives each invariant half exactly once.
D'Ascenzo and Michielan \cite[Section~3.1 and Appendix~A]{das}
bound \(\E T\) by \(O(\log k)\) and prove tightness of the number of
distinct common sizes in Lemma~A.1.

Their estimates also give qualitative convergence of \(T\) at equal
halves. To see this, take \(R=\lfloor\sqrt{k}\rfloor\). In the
unrestricted Ewens profiles their Lemma~A.3 excludes common sizes
in \([R,k-R]\) with probability tending to one. For fixed \(M<R\),
coupling the counts through length \(R\) by \eqref{s2:feller}
bounds the probability of a common size in \((M,R)\) by the
probability of a Poisson common sum greater than \(M\), plus
\(O(R/k)\). Their Lemma~A.2 shows that the latter Poisson intersection
is almost surely finite, so its tail tends to zero. The event of
simplicity on both sides has probability tending to \(e^{-3/2}>0\),
and hence the same exclusion holds in the simple graphs after letting
\(k\to\infty\) and then \(M\to\infty\). By complementation, on this
event \(T=2\sum_{s=0}^{M}A_{1,s}A_{2,s}\) for \(2M<k\). The
fixed-window counts converge to the independent Poisson profiles
with lengths one and two omitted, giving the limiting law
\(\mathcal L(W_0^{[3]})\). Corollary~\ref{s3:simple-cor} supplies
the rate in total variation and treats all imbalances \(d=n-2k\)
in the same estimate.

\section{Invariable generation and critical probabilities}\label{sec:generation}

An intransitive group preserves a proper subset, and conjugation can
align any sets of the same size. This leaves two questions: how often
do the samples have a common invariant size, and when does invariable
transitivity force generation of \(S_n\)? We first answer these in the
strict ranges. At the first discontinuity, long-cycle completions
transfer the critical Poisson estimates to finite degree. The same walk estimates will give the probability that a single
permutation has an invariant set of prescribed size.

\subsection{Invariable generation in the strict ranges}\label{s3:discontinuity}

Recall that \(p_{0,m}(\theta)\) is the probability that \(m\)
independent Poisson sumsets meet only at zero.

\begin{theorem}[The generation probability]\label{thm:strict-generation}
Fix an integer \(m\ge2\) and suppose \(0<\theta<(1-1/m)/\log 2\).
Then \(p_{0,m}(\theta)>0\), and the probability that \(m\)
independent Ewens permutations of parameter \(\theta\) invariably
generate \(S_n\) tends to \((1-2^{-m})p_{0,m}(\theta)\) as \(n\to\infty\).
\end{theorem}

\begin{proposition}[The strict converse]\label{thm:strict-converse}
Let \(\theta>0\) and let \(m\ge2\) be an integer such that
\[
                   m(1-\theta\log 2)<1.
\]
The intersection of \(m\) independent Poisson sumsets of parameter
\(\theta\) contains infinitely many positive integers almost surely.
The same assertion holds simultaneously after all lengths through any
fixed cutoff are omitted from every profile. The probability that
\(m\) independent Ewens permutations of parameter \(\theta\)
invariably generate \(S_n\) tends to zero.
\end{proposition}

We first prove Theorem~\ref{thm:strict-generation}, using a tail
bound for common sizes and an estimate for transitive overgroups.
The proof of Proposition~\ref{thm:strict-converse} follows. Its
Fourier estimate counts attainable differences, which are turned
into common positive sums by adding one copy to each profile.

\begin{lemma}[The tail of the common sizes]\label{lem:common-size-tail}
Suppose \(m\ge2\) and \(m(1-\theta\log 2)>1\). There is a
\(b_1=b_1(\theta,m)>0\) such that for \(M\ge1\),
\[
\begin{split}
 \Pp(\text{a common Poisson sum greater than }M)&=O_{\theta,m}(M^{-b_1}),\\
 \sup_{n\ge2}\Pp(\text{a common invariant size in }(M,n/2])
      &=O_{\theta,m}(M^{-b_1}).
\end{split}
\]
The \(m\)-fold Poisson intersection is almost surely finite, and
\(p_{0,m}(\theta)>0\).
\end{lemma}

\begin{proof}
In this proof the implied constants may depend on \(\theta,m\);
the auxiliary numbers \(\lambda,y\) below are fixed in terms of them.
Let \(P_k=\sum_{j\le k}Z_j\) count the copies through length \(k\)
in one Poisson profile. Weighting by \(2^{-P_k}\) changes the means
through length \(k\) to \(\theta/(2j)\), with normalizing factor
\(e^{-\theta H_k/2}\). Thus
\[
 \E(U_k2^{-P_k})
   =e^{-\theta H_k/2}\E_{\theta/2}U_k
   =e^{-\theta H_k/2}h_k(\theta/2)\ll_\theta k^{-1}.
\]
For the middle expectation, the Poisson generating function gives
\(\E_{\theta/2}\prod_{j\le k}(1+u^j)^{Z_j}
=\exp((\theta/2)\sum_{j\le k}u^j/j)\), whose coefficient of \(u^k\)
is \(h_k(\theta/2)\). For comparison, the same mark at \(\theta=1/2\) appears in
the proof of D'Ascenzo and Michielan \cite[Lemma~A.2]{das}.
On the finite side,
Lemma~\ref{lem:marked-fixed-sets} gives
\(\E(I_{n,k}2^{-K_{\le k}})\ll_\theta k^{-1}\) for
\(1\le k\le n/2\).

Choose \(\theta<\lambda<(1-1/m)/\log 2\) and put
\(\delta_\lambda=1-\lambda\log 2>1/m\). On
\(\{U_k>0,\ P_k\le\lambda\log k\}\) the marked count is at least
\(2^{-\lambda\log k}\). The corresponding statement holds for
\(I_{n,k}\) and \(K_{\le k}\), so each restricted probability is
\(O(k^{-\delta_\lambda})\). Independence gives a bound of
\(O(k^{-m\delta_\lambda})\) for a common size \(k\) when all
\(m\) counts satisfy the restriction.

To control violations, choose \(y>1\) so close to one that
\(b_0=\lambda\log y-\theta(y-1)>0\). Since \(k>M\ge1\),
the dyadic block containing \(k\) has index \(j\ge1\).
If \(2^j\le k<2^{j+1}\) and
\(P_k+1>\lambda\log k\), then \(P_{2^{j+1}}+1>\lambda j\log 2\).
The latter event has probability \(O(2^{-jb_0})\) by the Poisson
exponential moment. The pathwise bound \(K_{\le k}\le P_k+1\)
in the Feller coupling transfers the same dyadic bound to the
finite profiles. Sum the probability of violations over the profiles and
dyadic scales, and sum the \(m\)-sample probabilities of the
restricted events over \(k>M\). This gives both bounds with
\(b_1=\min(b_0,m\delta_\lambda-1)>0\).

In particular the Poisson intersection is almost surely finite.
Choose \(M\) so that having no common Poisson sum greater than \(M\)
has positive probability. Deleting every copy through length \(M\)
from the first profile can only increase this probability and
rules out common sums in \([1,M]\). The profile after deletion
has the law of the first profile conditional on an event of
probability \(e^{-\theta H_M}\), namely having no copies through
length \(M\). This proves \(p_{0,m}(\theta)>0\).
\end{proof}

\begin{lemma}[Reduction to few blocks]\label{lem:ewens-overgroup-reduction}
Fix \(0<\theta\le\lambda_0=1/\log 2\) and put \(R=\exp((\log n)^{3/5})\).
Except on an event of probability \(O_\theta((\log n)^{-2\theta/5})\), an
Ewens permutation \(\pi_n\) has the following properties. Every primitive
overgroup of \(\pi_n\) contains \(A_n\). The permutation preserves no
system of \(r>R\) equal-sized blocks of size at least two.
\end{lemma}

\begin{proof}
We use the block-system argument of Brito, Fowler, Junge, and
Levy \cite[Theorem~1.3]{brito}. Take \(n\) large enough that
\(R<n^{1/4}\) and \(n^{1/2}\le n/2\). The implied constant
covers the remaining degrees \(n\ge2\).
The choice \(\log R=(\log n)^{3/5}\) serves two purposes.
With exponent \(a\) in place of \(3/5\), the mean in the
large-prime-factor estimate below is \(\gg_\theta(\log n)^{2a-1}\).
For the two-sample lower bound in Theorem~\ref{thm:first-discontinuity},
the probability at \(\theta_c\) that neither sample
has a cycle of prime length in \((R,n^{1/4}]\) is \(O((\log n)^{-2\theta_c(1-a)})\).
We need \(a>1/2\) for the former and
\(a<1-(\log 2)/2\) for the latter exponent to exceed \(1/2\).

\emph{Primitive overgroups.}
We first find a prime \(R<p\le n^{1/4}\) for which a power of \(\pi_n\)
is a \(p\)-cycle. This cycle also rules out systems of small blocks.
The Poisson counts at the prime lengths in this range have
total mean \(\theta(\tfrac25\log\log n+O(1))\) by Mertens'
estimate \cite[Theorem~4.12]{apostol}. The probability of no such
cycle is \(O_\theta((\log n)^{-2\theta/5})\), using the Poisson
probability of zero and the \(O_\theta(n^{-3/4})\) error from
\eqref{s2:feller}.
For a prime \(p\) in the range, the expected number of ordered
pairs consisting of a \(p\)-cycle and a distinct cycle of
length divisible by \(p\) is
\[
 \sum_{\substack{a\ge1\\(a+1)p\le n}}
     \frac{\theta^2}{ap^2}\frac{h_{n-(a+1)p}}{h_n}
 \ll_\theta\frac{\log n}{p^2}+\frac{n^{-\theta}}p.
\]
The expectation follows by taking two cycle derivatives in
the profile generating function \eqref{eq:marked}. When \(a=1\),
the derivatives give the second factorial moment of the number
of \(p\)-cycles. For \(ap\le n/2\) and sufficiently large \(n\), we have
\(n-(a+1)p\ge n/2-n^{1/4}\ge n/4\). Thus \eqref{s2:gamma}
bounds the coefficient ratio uniformly. For \(ap>n/2\), use
\(1/(ap^2)\le 2/(np)\) and sum \(h_{n-(a+1)p}\) over an
arithmetic progression of step \(p\). For \(0<\theta<1\),
the term with smallest index \(n-(a+1)p<p\) is at most one
since \(h_j(\theta)\le1\) for \(j\ge0\).
Equation~\eqref{s2:gamma} bounds the remainder by a constant times
\(p^{\theta-1}\sum_{j\le n/p}j^{\theta-1}=O_\theta(n^\theta/p)\).
For \(\theta\ge1\), there are \(O(n/p)\) terms of size
\(O_\theta(n^{\theta-1})\), giving the same bound for the full sum.
In both cases the sum is \(O_\theta(1+n^\theta/p)\).
Since \(np h_n\asymp_\theta n^\theta p\), the contribution from
\(ap>n/2\) is \(O_\theta(n^{-\theta}/p+p^{-2})\).
The \(p^{-2}\) term is absorbed into \((\log n)/p^2\).
Since \(\sum_{p>R}p^{-2}\ll R^{-1}\) and
\(\sum_{p\le n^{1/4}}p^{-1}\ll\log\log n\) by Mertens' estimate,
the probability of any such pair is
\(O_\theta((\log n)/R+n^{-\theta}\log\log n)\).
In the absence of such a pair, the
selected \(p\)-cycle is the only cycle of length divisible by \(p\).
Let \(L_p\) be the least common multiple of the other cycle lengths.
Since \(p\nmid L_p\), the power \(g=\pi_n^{L_p}\) acts as a \(p\)-cycle
on the selected support and fixes every other point. Since \(p\le n^{1/4}\le n-3\)
for sufficiently large \(n\), Jordan's theorem \cite[Theorem~3.3E]{dm}
gives that every primitive overgroup contains \(A_n\).

\emph{Systems with more than \(R\) blocks.}
Suppose \(\pi_n\) preserves \(r\) blocks of size \(s=n/r\), where
\(r,s\ge2\). Any cycle of length \(j\) projects onto an orbit of \(t\)
blocks. Hence \(t\mid j\), and the cycle visits each of the \(t\) blocks
\(j/t\le s\) times.

For \(R<r\le n/R\), we use a cycle length divisible by a
prime \(q>n/R\). To see that such a cycle occurs with
probability tending to one, use the Poisson profile through
\(B=\lfloor n/\log n\rfloor\). The primes in
\((n/R,n/\sqrt R]\) exceed \(\sqrt n\) for large \(n\), so
their sets of multiples through \(B\) are disjoint. The mean
number of copies at those multiples is
\[
 \theta\sum_{n/R<q\le n/\sqrt R}
    \frac{H_{\lfloor B/q\rfloor}}q
 \gg_\theta\frac{(\log R)^2}{\log n}\to\infty,
\]
where the sum is over primes. Since
\(\lfloor B/q\rfloor\ge\sqrt R/(2\log n)\), the harmonic factor
is \(\gg\log R\). Also \(\log q\le\log n\), so
\(\sum_{n/R<q\le n/\sqrt R}1/q\gg (\log R)/\log n\) by
\(\sum_{q\le x}(\log q)/q=\log x+O(1)\) \cite[Theorem~4.10]{apostol}.
The chance of no such Poisson copy is
\(O_\theta(\exp(-c(\log n)^{1/5}))\), for some \(c=c(\theta)>0\).
For large \(n\), \(B\le n/2\), so the coupling error in \eqref{s2:feller} is
\(O_\theta(B/n)=O_\theta((\log n)^{-1})\). This is
\(o_\theta((\log n)^{-2\theta/5})\) since \(2\theta/5\le 2\lambda_0/5<1\).
For the resulting cycle length \(j\), the orbit
of blocks has length \(t\le r\le n/R<q\), and each block holds
\(j/t\le s=n/r<n/R<q\) points of the cycle. Neither factor of
\(j=t(j/t)\) is divisible by \(q\), which contradicts \(q\mid j\).

Finally, if \(r>n/R\), then \(s<R<p\). The \(p\)-cycle
\(g\in\langle\pi_n\rangle\) preserves every block system preserved by
\(\pi_n\). Its action on the blocks has order \(1\) or \(p\).
If it has order \(p\), it moves at least \(p\) blocks. This would
move every point in those blocks, for a total of at least \(ps>p\)
points. Thus the \(p\)-cycle fixes all blocks setwise. Its support
must then lie in one block, which is impossible when \(s<p\).
\end{proof}

\begin{lemma}[Transitive overgroups below \(\lambda_0\)]\label{lem:ewens-overgroups}
Let \(0<\theta<\lambda_0\) and \(n\ge2\). Except on an event of
probability \(O_\theta((\log n)^{-2\theta/5})\), every transitive
subgroup containing an Ewens permutation of parameter \(\theta\)
contains \(A_n\).
\end{lemma}

\begin{proposition}[Transitive overgroups at \(\lambda_0\)]
\label{prop:endpoint-overgroups}
Let \(\pi_n\) have the Ewens law of parameter \(\lambda_0\), and let
\(F_n\) be the event that some transitive subgroup containing \(\pi_n\)
does not contain \(A_n\). As \(n\to\infty\),
\[
 \Pp(F_n)=
 \begin{cases}
 i_{\lambda_0}(n,n/2)+O((\log n)^{-2\lambda_0/5}),&n\text{ even},\\
 O((\log n)^{-2\lambda_0/5}),&n\text{ odd}.
 \end{cases}
\]
In particular \(\Pp(F_n)\to0\), and along even degrees it has order
\((\log n)^{-1/2}\).
\end{proposition}

\begin{proof}[Proof of Lemma~\ref{lem:ewens-overgroups} and Proposition~\ref{prop:endpoint-overgroups}]
Write \(F_n\) for the event in Proposition~\ref{prop:endpoint-overgroups}
at any fixed \(0<\theta\le\lambda_0\). Throughout this proof the
constants may depend on \(\theta\). By Lemma~\ref{lem:ewens-overgroup-reduction},
apart from an event of probability \(O_\theta((\log n)^{-2\theta/5})\),
any transitive overgroup missing \(A_n\) preserves a system of
\(2\le u\le R\) blocks of size \(v=n/u\).
If the induced permutation of the blocks has orbit lengths
\(t_1,\ldots,t_\ell\), then \(\sum_i t_i=u\).
The union of the blocks in orbit \(i\) is an invariant set of size
\(t_iv\). Each of its cycles visits every block of the orbit the
same number of times, so its length is divisible by \(t_i\).

For a given tuple \(\boldsymbol t=(t_1,\ldots,t_\ell)\), let
\(X_{\boldsymbol t}\) count ordered partitions of \([n]\) into
invariant sets of sizes \(t_iv\) satisfying these divisibility conditions.
Marking each cycle by \(1/\ell\) gives
\begin{equation}\label{eq:overgroup-orbit-mark}
 \mathbb E_\theta[X_{\boldsymbol t}\ell^{-K_n}]
    =\frac{1}{h_n(\theta)}
       \prod_{i=1}^{\ell} h_v\!\left(\frac{\theta}{\ell t_i}\right).
\end{equation}
Indeed, on part \(i\) the marked generating function is
\(\exp((\theta/\ell)\sum_{j\ge1}z^{t_i j}/(t_i j))\).
Its coefficient of \(z^{t_i v}\) is the indicated factor,
and the multinomial count of ordered partitions cancels the factorials
of their sizes.

Fix \(\epsilon>0\), to be chosen below. As \(K_n\) is a sum of independent
Bernoulli variables with total mean \(\theta\log n+O_\theta(1)\), its
exponential moment gives \(K_n\le(\theta+\epsilon)\log n\) outside
an event of probability \(O_{\theta,\epsilon}(n^{-c})\), with
\(c=c(\theta,\epsilon)>0\). On this count event,
\(X_{\boldsymbol t}>0\) implies
\(X_{\boldsymbol t}\ell^{-K_n}\ge n^{-(\theta+\epsilon)\log\ell}\).
For \(0<a\le\theta\) and \(v\ge1\), we have the uniform bound
\[
 h_v(a)=\frac av\prod_{j=1}^{v-1}(1+a/j)
       \le\frac av\exp(aH_{v-1})
       \le \theta e^\theta v^{a-1},
\]
where \(H_{v-1}\le 1+\log v\). Substitution of \(v=n/u\) in
\eqref{eq:overgroup-orbit-mark} introduces the power
\(u^{\ell-(\theta/\ell)\sum_i t_i^{-1}}\), which is at most \(R^\ell\)
since \(1\le u\le R\). Let \(E_n=\{K_n\le(\theta+\epsilon)\log n\}\). It follows that
\[
 \mathbb P_\theta(X_{\boldsymbol t}>0, E_n)\ll_\theta (C_\theta R)^\ell n^{e(\boldsymbol t)},\qquad
 e(\boldsymbol t)=1-\theta-\ell
   +\frac{\theta}{\ell}\sum_i\frac1{t_i}
   +(\theta+\epsilon)\log\ell.
\]
For \(\ell=1\) this exponent is \(-\theta(1-1/u)\le-\theta/2\).
For \(\ell\ge2\) it is at most \(1-\ell+(\theta+\epsilon)\log\ell\).

Suppose first that \(\theta<\lambda_0\), and choose \(\epsilon\) so that
\((\theta+\epsilon)\log 2<1\). Since \(\log\ell/(\ell-1)\) decreases for
\(\ell\ge2\), we obtain \(e(\boldsymbol t)\le-c_1\ell\) for some
\(c_1=c_1(\theta)>0\) and every \(\ell\ge2\). For each \(\ell\) there
are at most \(R^\ell\) tuples to consider. Since \(R=n^{o(1)}\),
\[
 \sum_{\boldsymbol t}\mathbb P_\theta(X_{\boldsymbol t}>0,E_n)
    \ll_\theta R^2 n^{-\theta/2}
       +\sum_{\ell=2}^{\lfloor R\rfloor}
          (C_\theta R^2 n^{-c_1})^\ell
   \ll_\theta n^{-c_2}
\]
for some \(c_2=c_2(\theta)>0\). This proves
Lemma~\ref{lem:ewens-overgroups}.

At \(\theta=\lambda_0\), choose \(\epsilon>0\) so that
\((\theta+\epsilon)\log3<2\) and \(\epsilon\log2<\theta/4\).
The same argument bounds \(e(\boldsymbol t)\le-c_1\ell\) for
\(\ell\ge3\). If \(\ell=2\) and \(u\ge3\), at least one of the two
orbit lengths is at least two, so \(1/t_1+1/t_2\le3/2\). Using
\(\theta\log2=1\) gives \(e(\boldsymbol t)\le-\theta/4+\epsilon\log2<0\).
Summing as before, all these cases together with \(\ell=1\) have
probability \(O(n^{-c_3})\) on the count event, for some \(c_3>0\).
The remaining tuple is \((1,1)\), which consists of two invariant
halves and requires \(n\) even. Conversely, the unordered pair formed by an invariant half and its
complement has transitive stabilizer \(S_{n/2}\wr S_2\), whose index
\(\tfrac12\binom{n}{n/2}\ge3\) for even \(n\ge4\) excludes \(A_n\).
This proves the asserted formula for \(\mathbb P(F_n)\).
By Proposition~\ref{prop:large-theta} below, the invariant-half
probability at \(\lambda_0\) has order \((\log n)^{-1/2}\).
Since \(2\lambda_0/5>1/2\), this gives the stated order in even degrees.
\end{proof}

When \(\theta>\lambda_0\), Proposition~\ref{prop:large-theta} instead
gives a positive lower bound for the invariant-half probability in
even degrees. Thus the transitive-overgroup assertion of
\cite[Theorem~1.3]{brito} holds through \(\lambda_0\) and fails above it.
Proposition~\ref{prop:overgroup-counterexample} gives an explicit
bound at \(\theta=4\).

\paragraph{Signs and short counts.}
For each fixed integer \(M\ge1\), the sign of an Ewens permutation is asymptotically
fair and independent of its cycle counts through length \(M\).
Couple the profiles of degrees \(n\) and \(t\) using the same
Feller sequence, where \(2M\le t\le n\). By \eqref{s2:feller}
their short counts disagree with probability \(O_\theta(M/t)\).
The degree-\(t\) counts are functions of the first \(t\) Bernoulli
variables. The total cycle count at degree \(n\) is
\(K_n=\sum_{j=1}^n\xi_j\), and \(\operatorname{sgn}(\pi_n)=(-1)^{n-K_n}\).
The absolute value of the conditional mean of the sign, given the
first \(t\) variables, is therefore
\[
 \prod_{j=t+1}^n\left|1-\frac{2\theta}{\theta+j-1}\right|.
\]
For \(j\le1+\theta\) the absolute value of each factor is at most one.
For \(j>1+\theta\) it is at most \(\exp(-2\theta/(\theta+j-1))\). The harmonic series in the
exponent diverges, so for fixed \(t\) the product tends to zero as
\(n\to\infty\). For any test
function of the short counts bounded by one in absolute value, the
absolute expectation of its product with the sign is at most this
product plus \(O_\theta(M/t)\), uniformly over the test function.
Here the error comes from replacing the degree-\(n\) counts by those
at degree \(t\), which are measurable with respect to
\(\xi_1,\ldots,\xi_t\). The limit superior as \(n\to\infty\) is at most
\(O_\theta(M/t)\) for every fixed \(t\ge2M\). Sending \(t\to\infty\)
with \(M\) fixed gives asymptotic fairness and independence through
this cutoff.

Let \(E_{n,M}\) be the event that the \(m\) samples have no common
invariant size in \([1,M]\). Write \(p_m^{(M)}(\theta)\) for the
probability that the \(m\) Poisson sumsets have no common element
in \([1,M]\). To find the probability of \(E_{n,M}\) with
every sample even, multiply \(\ind_{E_{n,M}}\) by
\[
 2^{-m}\prod_{r=1}^m\bigl(1+\operatorname{sgn}(\pi_n^{(r)})\bigr)
\]
and expand. In each nonconstant term choose a sample whose sign occurs
and condition on the other \(m-1\) samples. The factor multiplying
the remaining sign is a test function of that sample's short profile,
bounded in absolute value by one. Each such term tends to zero by
the uniform bound above. The constant term is \(2^{-m}\Pp(E_{n,M})\),
and \eqref{s2:feller} gives \(\Pp(E_{n,M})\to p_m^{(M)}(\theta)\).
It follows that
\begin{equation}\label{eq:generation-short-sign}
 \Pp(E_{n,M},\ \text{at least one odd sample})
 \longrightarrow (1-2^{-m})p_m^{(M)}(\theta).
\end{equation}
Both conditions are necessary for invariable generation when \(n\ge2M\).
As \(M\to\infty\), the probabilities \(p_m^{(M)}(\theta)\) decrease
to \(p_{0,m}(\theta)\). Thus, for every fixed \(\theta>0\) and \(m\ge2\),
\begin{equation}\label{eq:generation-limsup}
 \limsup_{n\to\infty}G_{n,m}(\theta)\le (1-2^{-m})p_{0,m}(\theta).
\end{equation}

\paragraph{The transitivity obstruction.}
A common invariant size in \([1,n/2]\) precludes invariable
generation, since the invariant sets of that size can be aligned
by conjugation. In its absence every choice of conjugates generates
a transitive subgroup. Otherwise a smallest orbit would have size
at most \(n/2\) and would give a common invariant size. The overgroup
property in Lemma~\ref{lem:ewens-overgroups}
is invariant under conjugation. If at least one sample has
this property, every transitive group generated by conjugates
contains \(A_n\). An odd sample then forces the group to be \(S_n\).

\begin{proof}[Proof of Theorem~\ref{thm:strict-generation}]
For \(n\ge2M\), the probability on the left of
\eqref{eq:generation-short-sign} differs from \(G_{n,m}(\theta)\)
by at most the probability of a common size in \((M,n/2]\) plus the
probability that the first sample has an exceptional overgroup.
Lemmas~\ref{lem:common-size-tail} and~\ref{lem:ewens-overgroups} bound their
sum by \(O_{\theta,m}(M^{-b_1})+O_\theta((\log n)^{-2\theta/5})\).
Apply \eqref{eq:generation-short-sign} at fixed \(M\) and then send
\(M\) to infinity. Positivity of \(p_{0,m}(\theta)\) follows from
Lemma~\ref{lem:common-size-tail}.
\end{proof}

\begin{proof}[Proof of Proposition~\ref{thm:strict-converse}]
Throughout the proof implied constants may depend on \(\theta,m\)
and on the auxiliary parameters chosen in terms of them.
We use the difference-set argument of Eberhard, Ford, and Green
\cite[Section~3]{efg-invariable} in the form developed by Brito,
Fowler, Junge, and Levy \cite[Section~5]{brito}. Let \(Z^{(1)},\ldots,Z^{(m)}\) be the independent Poisson profiles, with \(Z_j^{(i)}\) of mean \(\theta/j\). For an interval \(I\), write \(\mathcal S_i(I)\) for the set of totals that can be obtained from the copies of lengths in \(I\) in profile \(i\). The strict inequality lets us choose \(\gamma\in(0,1)\) and \(\delta_F\in(0,1)\) so that
\[
                    a_F=\frac{m}{m-1}\theta\gamma(1-\delta_F)\log 2>1.
\]
Take an integer \(k\ge2\) and let \(I_k=(k^{1-\gamma},k]\). For a real cutoff \(1\le u\le k\), put \(Y_i(u)=\sum_{u<j\le k}Z_j^{(i)}\).

\paragraph{A simultaneous lower count bound.}
There is a constant \(A_F\), independent of \(k\), such that the event \(G_k\) given by
\[
                   Y_i(u)\ge \theta(1-\delta_F)\log(k/u)-A_F
                   \qquad(1\le i\le m,\ 1\le u\le k)
\]
has probability at least \(3/4\). To see this, set \(u_b=ke^{-b}\) for integers \(0\le b\le\lfloor\log k\rfloor\). The mean of \(Y_i(u_b)\) is \(\theta b+O_\theta(1)\). Choose \(v>0\) small enough that \(1-e^{-v}>v(1-\delta_F)\). The Poisson exponential moment bounds
\[
 \Pp\bigl(Y_i(u_b)<\theta(1-\delta_F)(b+1)-A_F\bigr)
       \ll_{\theta,v} e^{-vA_F}
                 e^{-\theta(1-e^{-v}-v(1-\delta_F))b}.
\]
Sum over \(b\) and \(i\) and take \(A_F\) large. For general \(u\) use \(b=\lfloor\log(k/u)\rfloor\). Then \(u\le u_b\), so \(Y_i(u)\ge Y_i(u_b)\), and \(\log(k/u)<b+1\).

\paragraph{The difference set.}
For \(\boldsymbol t=(t_1,\ldots,t_m)\) on the zero-sum torus
\[
      \mathbb T_0^m
       =\{\boldsymbol t\in(\mathbb R/\mathbb Z)^m:\ t_1+\cdots+t_m=0\}
\]
define
\[
 F_k(\boldsymbol t)
  =\prod_{i=1}^m\prod_{j\in I_k\cap\mathbb Z}
              \left(\frac{1+e^{2\pi\mathrm ijt_i}}2\right)^{Z_j^{(i)}}.
\]
Conditional on the copies in \(I_k\), include each independently with probability \(1/2\). The Fourier coefficients of \(F_k\) in the coordinates \(t_1,\ldots,t_{m-1}\) are the probabilities of the resulting difference vectors \((s_1-s_m,\ldots,s_{m-1}-s_m)\). Hence their support is the set
\[
        \mathcal D_k=\{(s_1-s_m,\ldots,s_{m-1}-s_m):
                                     s_i\in\mathcal S_i(I_k)\}.
\]
Parseval and Cauchy--Schwarz give
\[
              |\mathcal D_k|
                   \ge \left(\int_{\mathbb T_0^m}|F_k|^2\right)^{-1},
\]
where the integral uses normalized Haar measure.

Let \(\|t\|\) denote distance to the nearest integer, put \(x_i=\max(1,k\|t_i\|)\), and take \(u_i=k/\min(k^\gamma,x_i)\). The copies counted by \(Y_i(u_i)\) all have lengths in \(I_k\). The Poisson generating function gives
\[
 \E\left(|F_k(\boldsymbol t)|^2\prod_{i=1}^m2^{Y_i(u_i)}\right)
 \le
 \exp\left(\theta\sum_{i=1}^m\sum_{u_i<j\le k}
                                \frac{\cos(2\pi jt_i)}j\right)\ll_{\theta,m}1.
\]
A marked copy contributes \(1+\cos(2\pi jt_i)\). The factor from an
unmarked copy is \((1+\cos(2\pi jt_i))/2\in[0,1]\). The unmarked
terms therefore contribute at most one to the expectation. If \(x_i=1\),
the marked sum is empty. Otherwise \(u_i\ge1/\|t_i\|\), and the
geometric-series bound \(O(1/\|t_i\|)\) for partial sums of
\(e^{2\pi\mathrm ijt_i}\) gives \(O(1/(u_i\|t_i\|))\) for the cosine
sum by summation by parts.

On \(G_k\) the mark contributes a factor at least a constant times \(\prod_i(k/u_i)^{\theta(1-\delta_F)\log 2}\). Since \(1\le x_i\le\max(1,k/2)\le k\), we have \(k/u_i\ge x_i^\gamma\). If \(x_{i_*}\) is largest, it is at least the geometric mean of the other \(x_i\). Consequently
\[
      \E\bigl(\ind_{G_k}|F_k(\boldsymbol t)|^2\bigr)
                \ll_{\theta,m,\gamma,\delta_F}\prod_{i\ne i_*}x_i^{-a_F}.
\]
For each possible \(i_*\), parameterize the torus by the other \(m-1\) coordinates. Each of their integrals is bounded by
\[
              \int_{\mathbb R/\mathbb Z}\max(1,k\|t\|)^{-a_F}\,dt
                    \ll_{a_F} k^{-1}.
\]
Thus \(\E(\ind_{G_k}\int_{\mathbb T_0^m}|F_k|^2)\ll_{\theta,m} k^{1-m}\). Markov's inequality and \(\Pp(G_k)\ge3/4\) give \(|\mathcal D_k|\gg_{\theta,m} k^{m-1}\) with probability at least \(5/8\), after choosing the constant in this last bound sufficiently small. Each profile has expected total length at most \(\theta k\) in \(I_k\).
Choose an integer \(C_F\ge\max(2,8m\theta)\). By Markov's inequality,
the probability that the combined total exceeds \(C_F k\) is at most
\(1/8\). On its complement every difference belongs to
\([-C_F k,C_F k]^{m-1}\). Consequently
\[
         |\mathcal D_k|\gg_{\theta,m} k^{m-1},\qquad
         \mathcal D_k\subset[-C_Fk,C_Fk]^{m-1}
\]
hold together with probability at least \(1/2\).

\paragraph{Positive sums and transfer.}
Condition on the profiles through \(k\) on the event in the last
display. Require exactly one new copy in \([2C_Fk,3C_Fk]\) for
profile \(m\) and exactly one in \([C_Fk,4C_Fk]\) for each other
profile. The intervals contain only their integer lengths. The new
counts are independent of the profiles through \(k\), and the
requirements have joint probability bounded below by a positive
constant for all large \(k\). Given the length \(j_m\) of the new copy in profile \(m\), a vector \((d_1,\ldots,d_{m-1})\in\mathcal D_k\) specifies admissible lengths \(j_i=j_m-d_i\) in the other profiles. In each profile the conditional law of the unique length is proportional to \(1/j\) on the indicated interval. Thus, given \(j_m\), the probability of the tuple \((j_1,\ldots,j_{m-1})\) under these conditions is \(\asymp_{C_F,m} k^{1-m}\). Distinct difference vectors prescribe disjoint events for these lengths. It follows that with probability bounded below independently of \(k\), every profile has a selection with the same total \(s_m+j_m\), which is positive since each selection uses its new copy.

A success is determined by the lengths in \((k^{1-\gamma},4C_Fk]\). We used \(G_k\) to bound its probability and do not impose \(G_k\) in its definition. With \(C_F\) fixed, choose a sequence of \(k\)'s with disjoint intervals. The success events are independent, so the second Borel--Cantelli lemma gives infinitely many successes almost surely. The common totals tend to infinity along the successes since they are at least \(2C_Fk\). Since the lower endpoints of the intervals also tend to infinity, the conclusion holds simultaneously after every fixed lower cutoff.

For finite permutations, let \(E_M\) be the event that some
common total belongs to \(\{1,\ldots,M\}\). These events increase
to an event of probability one. Given \(\epsilon>0\), choose a
fixed \(M\) with \(\Pp(E_M)>1-\epsilon\). Since all lengths are
positive, \(E_M\) depends only on the counts of lengths at most
\(M\). For \(n\ge2M\), the Feller coupling \eqref{s2:feller}
transfers \(E_M\) to the \(m\) Ewens profiles with error
\(O_{\theta,m}(M/n)\). The common total is at most \(M<n\), so
the invariant sets are proper. Conjugates can align them and
prevent generation. Let \(n\to\infty\), then \(\epsilon\downarrow0\).
\end{proof}

\subsection{Exchanges at the critical parameter}\label{sec:critical-walk}

In this subsection and in Proposition~\ref{s3:critical-rate}, we
absorb the fixed parameter \(\theta_c=1/(2\log 2)\) into the implied constants.
The expected count in each dyadic block of the
merged profile is one up to a summable error. At the first maximum
of the centered count, the reciprocal identity bounds the mean
of \(V_N^{-\omega}\). For the lower bounds we keep the walk
below a barrier. The probability of making any specified net
change with a given choice of copies is controlled by the
largest block used. A copy from block \(v\) must hit one
prescribed length, an event of probability \(O(2^{-v})\).

\begin{lemma}[Poisson insertion]\label{lem:poisson-insertion}
Let \(Y=(Y_1,\ldots,Y_b)\) have independent \(\operatorname{Pois}(1)\)
coordinates, and let \(e_v\) be the \(v\)-th unit vector. For each
nonnegative function \(f\) on \(\mathbb Z_{\ge0}^b\),
\[
      \E[Y_v f(Y)]=\E[f(Y+e_v)].
\]
If the increments \(Y_j-1\) are read in any fixed order, insertion at
coordinate \(v\) raises every partial sum from that step onward by one.
In the reversed prefix \(Y_i-1,\ldots,Y_1-1\) the shift starts
at step \(i+1-v\) for \(v\le i\). The increments after time \(i\) are
unchanged.
\end{lemma}

\begin{proof}
Condition on the other coordinates and use
\(y\Pp(Y_v=y)=\Pp(Y_v=y-1)\) for \(y\ge1\). The assertion about
partial sums follows by adding the increments in the indicated order.
\end{proof}

\begin{lemma}[Critical walk estimates]\label{lem:critical-walk}
Let \(L_i^*=\sum_{v=1}^i(Y_v-1)\), where \(L_0^*=0\) and the
\(Y_v\sim\operatorname{Pois}(1)\) are independent. For integers
\(a,\ell\ge0\),
\begin{equation}\label{s3:critical-walk-survival}
 \Pp\!\left(\max_{0\le i\le \ell}L_i^*\le a\right)\ll
      \frac{a+1}{\sqrt{\ell+1}},
 \qquad
 \Pp\!\left(\min_{0\le i\le \ell}L_i^*\ge-a\right)\ll
      \frac{a+1}{\sqrt{\ell+1}}.
\end{equation}
For \(i,a\ge0\) the two bridge probabilities satisfy
\begin{equation}\label{s3:critical-walk-bridge}
 \Pp\!\left(\min_{0\le j\le i}L_j^*\ge0,\ L_i^*=a\right)
 +\Pp\!\left(\max_{0\le j\le i}L_j^*\le0,\ L_i^*=-a\right)
 \ll \frac{a+1}{(i+1)^{3/2}}.
\end{equation}
Fix \(0<\omega\le1\). For the exchange count \(V_N\) of two independent
Poisson profiles at \(\theta_c\),
\begin{equation}\label{s3:critical-poisson-upper}
 \E V_{2^\ell}^{-\omega}\ll_\omega(\ell+1)^{-1/2}.
\end{equation}
Finally, put \(G=\{\max_{0\le i\le \ell}L_i^*\le0\}\). Then
\begin{equation}\label{s3:critical-weighted-survival}
 \E\!\left[\ind_G\sum_{v=1}^{\ell}2^{L_v^*}\right]\ll(\ell+1)^{-1/2}
 \asymp \Pp(G),
 \qquad \E[Y_v\ind_G]\le \Pp(G)\quad(1\le v\le \ell).
\end{equation}
\end{lemma}

\begin{proof}
\emph{Survival and bridges.}
Let \(\tau_r\) be the first visit of \(L^*\) to \(-r\), for
an integer \(r\ge1\). Since the increments are at least \(-1\), the
law of the iterated logarithm gives \(\tau_r<\infty\) almost surely.
For \(N\ge1\), the ballot identity is
\[
 \Pp(\tau_r=N)=\frac rN\Pp(L^*_N=-r).
\]
To see this, extend a fixed word of \(N\) increments of sum \(-r\)
periodically in both directions, so its partial sums \(s_j\) satisfy
\(s_{j+N}=s_j-r\). A cyclic word ending at position \(j\) first
reaches its final level at its last step precisely when \(s_j\) is
strictly below the preceding \(N\) heights. Earlier periods have
larger heights, so this is a strict record minimum of the extended
path. Its running minimum drops by \(r\) in each period, one unit
at each strict record. Thus exactly \(r\) of the \(N\) cyclic
starting positions work. Conditional on the sum, the increments
are exchangeable, which gives the identity. When \(r>N\) both
probabilities vanish.

Stirling's formula gives \(\sup_x\Pp(L^*_N=x)\ll N^{-1/2}\) and
\(\Pp(L^*_N=-1)\asymp N^{-1/2}\), since \(N+L^*_N\) is Poisson
of mean \(N\). Summing the ballot identity over \(N>\ell\) gives
\[
 \Pp\!\left(\min_{0\le j\le\ell}L^*_j\ge-a\right)
   =\Pp(\tau_{a+1}>\ell)\ll (a+1)(\ell+1)^{-1/2},
\]
with a matching lower bound at \(a=0\).

For the other direction, reverse the first \(i\) increments of
a path ending at \(-a\). The original path stays at most zero
precisely when the reversed path stays at least \(-a\).
Appending an independent increment equal to \(-1\), an event of
probability \(e^{-1}\), then gives a first visit to \(-a-1\).
Conversely, a first visit to \(-a-1\) at time \(i+1\) must make
its last step from \(-a\). Hence for \(0\le a\le i\),
\begin{equation}\label{eq:critical-ballot}
 \Pp\!\left(\max_{0\le j\le i}L^*_j\le0,\ L^*_i=-a\right)
 =e\Pp(\tau_{a+1}=i+1)
 =\frac{e(a+1)}{i+1}\Pp\bigl(\operatorname{Pois}(i+1)=i-a\bigr).
\end{equation}
At \(a=0\), reversal gives the same probability for a nonnegative
bridge ending at zero. Summing over \(a\), with \(N=i+1\) and
\(X\sim\operatorname{Pois}(N)\), gives \(e\E(N-X)_+/N=e\Pp(X=N-1)\).
The last equality follows by summing \(j\Pp(X=j)=N\Pp(X=j-1)\)
for \(1\le j\le N-1\). Together with the preceding estimate this
gives
\begin{equation}\label{eq:critical-zero-survival}
 \Pp\!\left(\max_{0\le j\le i}L^*_j\le0\right)
 \asymp\Pp\!\left(\min_{0\le j\le i}L^*_j\ge0\right)
 \asymp(i+1)^{-1/2}.
\end{equation}

Let \(D\) be the first time \(L^*\) is positive. The successive
waiting times between strict record highs are independent copies of
\(D\), and each new record increases the height by at least one.
Thus the first crossing of \(a\) occurs no later than the sum of
\(a+1\) such waiting times. If their sum exceeds \(\ell\ge1\),
their sum after truncating each at \(\ell\) is at least \(\ell\).
Markov's inequality and \eqref{eq:critical-zero-survival} give
\[
 \Pp\!\left(\max_{0\le j\le\ell}L^*_j\le a\right)
 \le\frac{a+1}{\ell}\E\min(D,\ell)
 =\frac{a+1}{\ell}\sum_{t=0}^{\ell-1}\Pp(D>t)
 \ll\frac{a+1}{\sqrt{\ell+1}}.
\]
This proves \eqref{s3:critical-walk-survival}.

Equation~\eqref{eq:critical-ballot} and the Poisson point-mass bound
give the nonpositive bridge estimate in \eqref{s3:critical-walk-bridge}.
For a nonnegative bridge ending at \(a\), take the first and last
\(r=\lfloor i/3\rfloor\) increments. The first block must stay
nonnegative. The last block, read in reverse, must have maximum
at most \(a\), since its partial sums are \(a-L^*_{i-j}\).
Given these two blocks, the middle \(i-2r\) increments must have
one specified sum. By independence, \eqref{s3:critical-walk-survival},
and the Poisson point-mass bound, the probability is
\(O((r+1)^{-1/2}(a+1)(r+1)^{-1/2}(i-2r)^{-1/2})\) for \(i\ge3\).
This gives \eqref{s3:critical-walk-bridge}, with the remaining
values \(i\le2\) covered by enlarging the constant.

\emph{The first maximum.}
Put \(B_v=J_{2^v}^{(2)}-J_{2^{v-1}}^{(2)}\) for \(1\le v\le \ell\)
and \(L_i^\circ=\sum_{v=1}^i(B_v-1)\), with \(L_0^\circ=0\).
Let \(\tau\) be the first time at which \(L_i^\circ\) reaches its
maximum over \(0\le i\le \ell\). As noted after \eqref{s3:reciprocal-support}, the cutoff may depend on the
merged profile. In particular \(\tau\) is determined before the side
marks are drawn, so those marks are still independent and fair at \(2^\tau\). By monotonicity
of \(V_N\) and conditional Jensen's inequality for the concave
function \(x^\omega\) (\(0<\omega\le1\)),
\[
 \mathbb E V_{2^\ell}^{-\omega}\le
 \mathbb E\bigl[\bigl((T_{2^\tau}^{(2)}+1)2^{-J_{2^\tau}^{(2)}}\bigr)^\omega\bigr].
\]
Since \(T_{2^i}^{(2)}\le J_1^{(2)}+\sum_{v=1}^i2^vB_v\) and
\(J_{2^i}^{(2)}=J_1^{(2)}+i+L_i^\circ\),
\[
 \bigl((T_{2^\tau}^{(2)}+1)2^{-J_{2^\tau}^{(2)}}\bigr)^\omega
 \le 2^{-\omega L_\tau^\circ}
       \left(1+\sum_{v=1}^{\tau}2^{v-\tau}B_v\right).
\]
Indeed, \((j+1)2^{-j}\le1\) for every integer \(j\ge0\).
Taking \(j=J_1^{(2)}\) removes the length-one contribution, and
\((1+x)^\omega\le1+x\) for \(x\ge0\).

The independent \(B_v\) have Poisson means \(\mu_v<1\), with
\(\sum_{v\ge1}(1-\mu_v)<\infty\), as in the proof of
Theorem~\ref{s3:poisson-critical}. Their joint law through \(\ell\) is bounded
above by a constant times that of the \(Y_v\), since at
\(x_1,\ldots,x_\ell\) the ratio is
\(\exp(\sum_{v=1}^\ell(1-\mu_v))\prod_{v=1}^\ell\mu_v^{x_v}\).
Since \(\mu_v<1\), the product is at most one, and the exponential
factor is bounded independently of \(\ell\). It suffices to bound the
expectation on the right with
\(B_v\) replaced by \(Y_v\) and \(\tau\) the corresponding first
maximum of \(L_i^*\).

On \(\{\tau=0\}\) the walk stays at most zero, which by
\eqref{s3:critical-walk-survival} has probability \(O((\ell+1)^{-1/2})\).
If \(\tau=i\ge1\), put \(a=L_i^*\) and write
\(S'_t=L_i^*-L_{i-t}^*\) for \(0\le t\le i\).
Then \(S'_t\ge1\) for \(1\le t\le i\) because this is the
first maximum, and \(S'_i=a\ge1\). The reversed prefix has the
same law as the original prefix and is independent of the
suffix, which stays nonpositive relative to \(L_i^*\). The
bridge estimate \eqref{s3:critical-walk-bridge} for the reversed
prefix and the survival estimate \eqref{s3:critical-walk-survival}
at height zero for the suffix give
\[
  \Pp(\tau=i,L_i^*=a)
  \ll \frac{a+1}{(i+1)^{3/2}(\ell-i+1)^{1/2}}.
\]
The same bound holds with \(Y_v\) inserted in the expectation for
\(v\le i\). By Lemma~\ref{lem:poisson-insertion}, strict positivity of
the reversed prefix after insertion requires a nonnegative path ending
at \(a-1\) before insertion. Its probability is
\(O(a(i+1)^{-3/2})\) by \eqref{s3:critical-walk-bridge}.
The suffix increments are unchanged.

Summing first over \(v\) with weights \(2^{v-i}\), and then over
\(a\) with weights \(2^{-\omega a}\), gives a bound
\(O_\omega((i+1)^{-3/2}(\ell-i+1)^{-1/2})\) for the contribution with
\(\tau=i\). Splitting the sum over \(i\) at \(\ell/2\) gives
\eqref{s3:critical-poisson-upper}.

\emph{Weighted survival.}
For \(1\le v\le \ell\) and \(a\ge0\), first use the nonpositive
bridge bound up to time \(v\), and then use the upper barrier bound
from position \(-a\) for the remaining \(\ell-v\) steps. This gives
\[
 \Pp(G,L_v^*=-a)
 \ll\frac{(a+1)^2}{(v+1)^{3/2}(\ell-v+1)^{1/2}}.
\]
Sum with weights \(2^{-a}\) and split the sum over \(v\) at \(\ell/2\)
to obtain the first bound of \eqref{s3:critical-weighted-survival}.
The comparison with \(\Pp(G)\) follows from
\eqref{eq:critical-zero-survival}. Since increasing any coordinate can
only remove an outcome from \(G\), Lemma~\ref{lem:poisson-insertion}
also gives \(\E[Y_v\ind_G]\le\Pp(G)\).
\end{proof}

For the lower bounds we delete the first few lengths and use the
same block count for exact and shifted exchanges.

\begin{lemma}[Exchanges after deletion of short lengths]\label{lem:critical-toggle}
In two independent critical Poisson profiles let \(N=2^\ell\) for
integers \(0\le b<\ell\), and let \(\mathcal H_b\) be the event that there are no
copies through length \(2^b\). For \(b<v\le \ell\), add independent
\(E_v\sim\operatorname{Pois}(1-\mu_v)\) to the merged block counts
\(B_v\), as in the proof of Theorem~\ref{s3:poisson-critical}, and put \(Y_v=B_v+E_v\).
Write
\[
 L_i^{(b)}=\sum_{u=1}^i (Y_{b+u}-1),\qquad
 G_b=\left\{\max_{0\le i\le \ell-b}L_i^{(b)}\le0\right\}.
\]
For any integer \(t\), let \(D_N(t)=R_N(S_N+t)-\ind_{\{t=0\}}\),
the number of nonempty toggles with net change \(t\). Then
\[
 \E[D_N(t)\ind_{G_b}\mid \mathcal H_b]\ll 2^{-b}\Pp(G_b),
 \qquad
 \Pp(G_b,\ T_N^{(2)}>x\mid \mathcal H_b)\le \frac{2N}{x}\Pp(G_b)
 \quad(x>0).
\]
The constants are uniform in \(b,\ell,t\).
\end{lemma}

\begin{proof}
Condition on \(\mathcal H_b\) and on the sigma-field \(\mathcal F_b\) generated
by the block counts, the additions, and the side marks of the copies
indexed before their lengths are sampled. The event \(G_b\) is
\(\mathcal F_b\)-measurable. Given these data, lengths of copies in
block \(v\) are independent draws with probabilities proportional to
\(1/j\) on \(2^{v-1}<j\le2^v\). The reciprocal sum on the block is at least \(1/2\).
As each \(j>2^{v-1}\), the largest point mass is less than
\(4\cdot2^{-v}\).

A toggle whose largest occupied block is \(b+i\) selects a
subset of the copies in blocks \(b+1,\ldots,b+i\), including at
least one from block \(b+i\). There are at most
\(2^{\sum_{u=1}^iB_{b+u}}\) choices. For each choice, reveal all
selected lengths except for one in the last block. The coefficient
of the unrevealed length in the net change is \(1\) or \(-1\),
so a net change of \(t\) fixes its value. The conditional expected
number of such toggles with net change \(t\) is therefore at most
a constant times
\[
  2^{-b-i}2^{\sum_{u=1}^iB_{b+u}}
  \le 2^{-b}2^{L_i^{(b)}}.
\]
Sum over \(i\) and apply \eqref{s3:critical-weighted-survival} to the
restarted walk. The conditioning on \(\mathcal H_b\) leaves that walk unchanged.
This proves the exchange bound. On \(\mathcal H_b\) we have
\(T_N^{(2)}\le\sum_{i=1}^{\ell-b}2^{b+i}Y_{b+i}\). The insertion estimate in
\eqref{s3:critical-weighted-survival} bounds the mean of this sum on
\(G_b\) by \(2N\Pp(G_b)\). Markov's inequality gives the length bound.
\end{proof}

\subsection{Long-cycle completions and the first discontinuity}\label{sec:critical-two-sample}

\begin{lemma}[Long-cycle completion]\label{lem:long-completion}
Fix \(m\ge2\), \(0<\theta<\lambda_0\), and \(0<\epsilon\le1/16\).
For \(\epsilon n/2<N\le\epsilon n\), let \(B_N\) be the event
that \(m\) independent Poisson profiles through length \(N\) have
no common positive sum and combined total length at most \(n/4\).
For all sufficiently large \(n\),
\begin{equation}\label{eq:generation-completion}
 G_{n,m}(\theta)\ge c_{\theta,m,\epsilon}\Pp(B_N)
     -O_{\theta,m}\bigl((\log n)^{-2m\theta/5}\bigr),
\end{equation}
where \(c_{\theta,m,\epsilon}>0\).
\end{lemma}

\begin{proof}
Fix a tuple of short profiles in \(B_N\). Their combined length is at
most \(n/4\), so each individual total is at most \(n/4\). Complete every sample
except the second by adding a single cycle of length \(n-u\),
where \(u\) is its short total. The second sample is completed by
one or two cycles. Changing the number of completing cycles by one
reverses the sign, so choose the option that makes this sample odd.
If its short total is \(u'\), the two-cycle option uses lengths
\(a\) and \(n-u'-a\), with \(n/3\le a<3n/8\). The inequality
\[
 n-u'-a>n-n/4-3n/8=3n/8>a\ge n/3>n/4>N
\]
gives distinct lengths above the cutoff. There are
\(\asymp n\) choices of \(a\), each giving a different completed
profile. Their combined coefficient weight
\(\sum_a\theta^2/[a(n-u'-a)]\) has order \(n^{-1}\), as does the
weight \(\theta/(n-u)\) of a one-cycle completion.

By \eqref{eq:profile}, the ratio of the Ewens probability of a
completion to the Poisson probability of its short profile is its
completion weight multiplied by \(e^{\theta H_N}/h_n\). The latter
factor has order \(n\), since \(N\asymp_\epsilon n\).
The short profiles can be recovered from each completion by
discarding all cycles longer than \(N\). Summing over distinct tuples
gives a completed event of probability at least
\(c_{\theta,m,\epsilon}\Pp(B_N)\). A common size at most \(n/2\)
would use only short cycles in the first sample, since its completing
cycle has length at least \(3n/4\). It would therefore be at most
\(n/4\), too small to use a completing cycle in any sample. This
contradicts the defining condition of \(B_N\).

Write \(C\) for this completion event and \(F_i\) for the event that
sample \(i\) belongs to a transitive subgroup missing \(A_n\).
On \(C\), every choice of conjugates generates a transitive subgroup
containing an odd permutation. If some \(F_i\) fails, that subgroup
contains \(A_n\) and hence equals \(S_n\). Consequently
\[
  G_{n,m}(\theta)\ge \Pp(C)-\Pp\!\left(\bigcap_{i=1}^m F_i\right).
\]
The \(F_i\) are independent under the original product Ewens law.
Lemma~\ref{lem:ewens-overgroups} thus gives
\(\Pp(\bigcap_iF_i)=\prod_i\Pp(F_i)\ll_{\theta,m}
(\log n)^{-2m\theta/5}\), which proves \eqref{eq:generation-completion}.

\end{proof}

\begin{theorem}[The first discontinuity]\label{thm:first-discontinuity}
At \(\theta_c=(2\log 2)^{-1}\), one has \(s_{\theta_c}=m_{\theta_c}=3\).
The probability that two independent Ewens permutations with parameter
\(\theta_c\) invariably generate \(S_n\) has order
\((\log n)^{-1/2}\) for all sufficiently large \(n\).
\end{theorem}

\begin{proof}
Theorem~\ref{s3:poisson-critical} gives \(s_{\theta_c}\ge3\). Since
\(3(1-\theta_c\log2)=3/2>1\), Theorem~\ref{thm:strict-generation}
gives \(s_{\theta_c}=3\) and \(m_{\theta_c}\le3\). The two-sample
upper bound below gives the remaining inequality \(m_{\theta_c}\ge3\).

\emph{Upper bound.} Take \(N=2^{\lfloor\frac12\log_2 n\rfloor}\)
in the two Poisson profiles and couple them to the Ewens profiles
through length \(N\). For \(n\ge4\), \(N\le\sqrt n\le n/2\),
so \eqref{s2:feller} gives error \(O(n^{-1/2})\).
On agreement, each equal exchange gives a total attained by both
Ewens samples. Invariable generation excludes every total strictly
between \(0\) and \(n\). Only the empty selections attain \(0\),
and the only possible selections attaining \(n\) use every cycle in
both samples. Thus \(V_N\le2\), and
\[
 G_{n,2}(\theta_c)\le \Pp(V_N\le2)+O(n^{-1/2})
     \le 2\E(1/V_N)+O(n^{-1/2}).
\]
Equation~\eqref{s3:critical-poisson-upper} gives the upper bound.

\emph{Lower bound.} Fix an integer \(b\ge1\) and a constant
\(0<\epsilon\le1/16\), to be chosen in that order. For \(n\) sufficiently
large in terms of \(b\) and \(\epsilon\), choose \(N=2^\ell\) with
\(\epsilon n/2<N\le\epsilon n\) and \(\ell>b\).
Let \(\mathcal H_b\) be the event that the two Poisson profiles have no copies
through length \(2^b\). Its probability is a positive constant and it
is independent of the counts in the remaining blocks. Add independent
Poisson counts in blocks \(b+1,\ldots,\ell\) and use the restarted walk
\(L_i^{(b)}\) and its survival event \(G_b\) from
Lemma~\ref{lem:critical-toggle}. In particular,
\(\Pp(G_b)\asymp(\ell-b+1)^{-1/2}\).

Apply Lemma~\ref{lem:critical-toggle} with \(t=0\) and \(x=n/4\).
Since \(D_N(0)=V_N-1\), its two bounds give
\[
 \Pp\!\left(G_b,\ \{V_N>1\}\cup\{T_N^{(2)}>n/4\}\mid\mathcal H_b\right)
 \le (C2^{-b}+8\epsilon)\Pp(G_b).
\]
Choose \(b\) large and then \(\epsilon\) small so that the
parenthesis is at most \(1/2\). This gives probability \(\gg(\log n)^{-1/2}\) for
\(\mathcal H_b\cap G_b\cap\{V_N=1,\ T_N^{(2)}\le n/4\}\).
After the additions are integrated out, the probability of short
profiles with no common positive sum and total length at most \(n/4\)
is at least this large.

Apply Lemma~\ref{lem:long-completion} with \(m=2\).
The overgroup error has exponent
\(4\theta_c/5=2/(5\log 2)>1/2\), so it is
\(o((\log n)^{-1/2})\). Equation~\eqref{eq:generation-completion}
therefore gives the lower bound.
\end{proof}

The same decay describes the least positive common total \(H_{\min}\)
of two independent critical Poisson profiles. For sufficiently large \(x\), take
dyadic \(M,N\) with \(M\le\sqrt{x}<2M\) and \(x\le N<2x\).
If \(V_N=1\), no positive total at most \(x\) is common, since all
copies used to attain such a total have length at most \(N\).
Choose \(b\) as in the lower bound above. Conditional on \(\mathcal H_b\),
Lemma~\ref{lem:critical-toggle} shows that at most half the probability
of \(G_b\) is lost by requiring \(V_N=1\). Since \(\Pp(\mathcal H_b)>0\)
is fixed, \(\Pp(V_N=1)\gg(\log N)^{-1/2}\).
Conversely, if \(V_M>1\) and \(T_M^{(2)}\le x\), a nonempty equal
exchange gives a common total at most \(x\). Hence
\[
  \Pp(V_N=1)\le \Pp(H_{\min}>x)
  \le \E(1/V_M)+\Pp(T_M^{(2)}>x)\ll(\log x)^{-1/2},
\]
by \eqref{s3:critical-poisson-upper} and \(\E T_M^{(2)}=2\theta_c M\).
Thus \(\Pp(H_{\min}>x)\asymp(\log x)^{-1/2}\) as \(x\to\infty\).

\begin{corollary}[Equality of the two thresholds]\label{cor:threshold-equivalence}
For each \(\theta>0\) and each fixed \(m\ge2\),
\[
  \liminf_{n\to\infty}G_{n,m}(\theta)>0
  \quad\Longleftrightarrow\quad p_{0,m}(\theta)>0 .
\]
In particular, \(s_\theta=m_\theta\) for every \(\theta>0\).
\end{corollary}

\begin{proof}
Equation~\eqref{eq:generation-limsup} gives necessity. When
\(\theta\ge\lambda_0\), Proposition~\ref{thm:strict-converse}
rules out positivity in both models. Suppose \(0<\theta<\lambda_0\)
and \(p_{0,m}(\theta)>0\). The probability of no common positive
sum in the Poisson profiles through length \(N\) is at least
\(p_{0,m}(\theta)\), and their combined total length has mean
\(m\theta N\). Choose \(\epsilon\le1/16\) so that
\(4m\theta\epsilon\le p_{0,m}(\theta)/2\), and take
\(N=\lfloor\epsilon n\rfloor\). For \(n\ge 2/\epsilon\) we have
\(\epsilon n/2<N\le\epsilon n\). Markov's inequality gives
\(\Pp(B_N)\ge p_{0,m}(\theta)/2\) for the event in
Lemma~\ref{lem:long-completion}. Thus
\eqref{eq:generation-completion} gives a positive limit inferior.
One profile almost surely has positive sums, because the sum
of its Poisson means diverges, and one permutation cannot
generate \(S_n\) for \(n\ge3\). This proves equality of the thresholds.
\end{proof}

For \(0<\theta<\lambda_0\), put \(x_\theta=(1-\theta\log 2)^{-1}\)
and \(g_\theta=\lceil x_\theta\rceil\). The two strict-range results
give
\[
 \begin{cases}
 s_\theta=m_\theta=g_\theta,&x_\theta\notin\mathbb Z,\\
 g_\theta\le s_\theta=m_\theta\le g_\theta+1,&x_\theta\in\mathbb Z.
 \end{cases}
\]
At \(x_\theta=2\), Theorem~\ref{thm:first-discontinuity} gives the
upper of the two choices. At higher integer values of \(x_\theta\)
there remains one common threshold to determine. For
\(\theta\ge\lambda_0\), both thresholds are infinite.

\begin{remark}[A summable perturbation]\label{rem:summable-perturbation}
Suppose the finite-dimensional cycle counts of a sequence of
permutation laws converge to independent counts \((Q_j)_{j\ge1}\)
satisfying
\[
 \sum_{j\ge1}\TV\bigl(\Law(Q_j),
                   \operatorname{Pois}(\theta_c/j)\bigr)<\infty.
\]
Then with probability tending to one, two independent samples
have conjugates that generate an intransitive group. Indeed,
independent maximal couplings and the Borel--Cantelli lemma
make each limiting profile agree with its Poisson counterpart
above an almost surely finite length. The simultaneous cutoff
assertion in Theorem~\ref{s3:poisson-critical} gives infinitely many
common positive sums for the limiting profiles. Given \(\epsilon>0\),
choose \(M\) so that one of those sums lies in \([1,M]\) with
probability at least \(1-\epsilon\). This event depends only on
counts through \(M\), so its probability converges under the finite
laws. For \(n>M\), conjugation aligns the two proper invariant sets.
Now let \(\epsilon\downarrow0\).

For a logarithmic structure of parameter \(\theta\), the uniform logarithmic
condition of Arratia, Barbour, and Tavar\'e
\cite[Eqs.~(1.11)--(1.12) and Theorem~3.1]{abt2000} gives fixed-window
convergence to the defining counts \(Q_j\). The summability condition of
Brito et al.\ \cite[Eq.~(1.2)]{brito} at the \emph{same} parameter gives,
with the constants \(c_r\) of \eqref{eq:strong-logarithmic-condition},
\[
 \TV\bigl(\Law(Q_j),\operatorname{Pois}(\theta/j)\bigr)
     \le \Bigl(c_1+\sum_{r\ge2}c_r\Bigr)\frac{e(j)}j
          + O_\theta(j^{-2}).
\]
Indeed, this bounds the discrepancy at one and the sum of the two
probabilities of two or more copies. The discrepancy at zero is
bounded by their sum. The comparison applies here with \(\theta=\theta_c\).
For a different \(\theta\), the two Poisson probabilities of zero already
differ by \(\lvert\theta-\theta_c\rvert/j+O_\theta(j^{-2})\).
Under the printed summability condition alone, fixed-window convergence
requires a separate argument, as Proposition~\ref{prop:lacunary-logarithmic} shows.

\end{remark}

\begin{proof}[Proof of Theorem~A]
Theorem~\ref{thm:strict-generation} gives the positive limit, including
the factor contributed by the signs. Proposition~\ref{thm:strict-converse}
gives the zero limit and the almost-sure assertion after each cutoff in
the opposite range. In particular, no finite number of samples suffices
when \(\theta\ge\lambda_0\). Theorem~\ref{s3:poisson-critical} gives
the assertion about two Poisson profiles at equality and above it.
Theorem~\ref{thm:first-discontinuity} gives the two-sample rate as
\(n\to\infty\), and three samples at \(\theta_c\) lie in the positive
strict range. Finally Corollary~\ref{cor:threshold-equivalence} covers
the equality of thresholds at every parameter, including the higher jumps.
The rate concerns sufficiently large degrees. The all-degree convention
and the exceptional degree six are treated in
Proposition~\ref{prop:all-degree-threshold}.
\end{proof}

\subsubsection*{Critical reciprocal moments}
\addcontentsline{toc}{subsubsection}{Critical reciprocal moments}\label{sec:critical-reciprocal}

The same survival estimate applies when a prescribed \(k\)-set is
invariant. The imbalance \(d=n-2k\) then allows a second kind of
exchange, with shift \(d\).

\begin{proposition}[Critical reciprocal moments]\label{s3:critical-rate}
At \(\theta=\theta_c=(2\log 2)^{-1}\), for each fixed
\(\omega>0\),
\[
 \mathbb E_A I_{n,k}^{-\omega}\asymp_\omega(1+\log k)^{-1/2}
\]
uniformly for \(1\le k\le n/2\). In particular,
\(i_{\theta_c}(n,k)/\mathbb E I_{n,k}\asymp(1+\log k)^{-1/2}\).
\end{proposition}

\begin{proof}
\emph{The upper bound.}
It suffices to consider \(0<\omega\le1\), since for \(\omega>1\)
we have \(I_{n,k}^{-\omega}\le I_{n,k}^{-1}\) under \(\Pp_A\).
For \(k\ge4\), let \(\ell=\lfloor\tfrac12\log_2 k\rfloor\), so
\(1\le \ell\) and \(2^\ell\le\sqrt k\le k/2\). Couple both conditioned
side profiles through length \(2^\ell\) to their Poisson counterparts,
with error \(O(2^\ell/k)\) by \eqref{s2:feller}, uniformly in \(n\ge2k\).
On agreement, each Poisson copy through length \(2^\ell\) can be
matched to a cycle on its own side. Every exchange counted by
\(V_{2^\ell}\) then gives a distinct invariant \(k\)-set. Both
reciprocal powers lie in \([0,1]\), so
\(\mathbb E_A I_{n,k}^{-\omega}\le\mathbb E V_{2^\ell}^{-\omega}
+O(2^\ell/k)\). The expectation is \(O_\omega((\log k)^{-1/2})\) by
\eqref{s3:critical-poisson-upper}, and the coupling error is \(O(k^{-1/2})\).
For \(1\le k\le3\), use \(\mathbb E_A I_{n,k}^{-\omega}\le1\).

\emph{The lower bound.}
Take \(0<\epsilon\le1/16\), and let \(N=2^\ell\) satisfy
\(\epsilon k/2<N\le\epsilon k\), with \(k\) large enough that \(\ell\ge1\).
Use the construction of Lemma~\ref{lem:critical-toggle} with \(b=0\).
Thus \(\mathcal H_0\) is the event \(J_1^{(2)}=0\), of fixed positive probability.
The survival event \(G_0\) is independent of \(\mathcal H_0\) and has probability
\(\asymp(\ell+1)^{-1/2}\). Put
\[
 E=\mathcal H_0\cap G_0\cap\{T_N^{(2)}\le k/4\}.
\]
By the length bound,
\(\Pp(E)\ge \Pp(\mathcal H_0)\Pp(G_0)(1-8N/k)\).
Since \(8N/k\le8\epsilon\le1/2\) and \(E\subseteq G_0\), we have
\(\Pp(E)\asymp(\ell+1)^{-1/2}\).

The empty toggle contributes exactly \(1+\ind_{\{d=0\}}\) to
\(V_{N,d}\). From \eqref{eq:shifted-exchange},
\[
 V_{N,d}=1+\ind_{\{d=0\}}+D_N(0)+D_N(d).
\]
The exchange bound of Lemma~\ref{lem:critical-toggle} gives
\(\E[V_{N,d}\ind_E]\ll\Pp(E)\) uniformly in \(d\). Jensen's inequality
for \(x\mapsto x^{-\omega}\), conditional on \(E\), now gives
\begin{equation}\label{s3:critical-poisson-lower}
 \E\!\left[V_{N,d}^{-\omega}\ind_E\right]
 \gg_\omega(\ell+1)^{-1/2}.
\end{equation}

\emph{Transfer to finite profiles.}
Let \(\mathcal C\) be the pair of short Poisson profiles through
length \(N\), and write \(u,u'\) for their total lengths.
Integrate out the additions by putting \(w(\mathcal C)=
\Pp(E\mid\mathcal C)\). Then \(0\le w\le1\), \(w=0\) unless
\(u+u'\le k/4\), and
\[
 \E[V_{N,d}^{-\omega}w(\mathcal C)]
     =\E[V_{N,d}^{-\omega}\ind_E],
\]
since the exchange count is a function of \(\mathcal C\).

Fix \(C_0\ge2\) large enough that \(k\le c_{\theta_c}(k+d)/2\)
whenever \(k+d>C_0k\), where \(c_{\theta_c}\) is the constant in
Lemma~\ref{lem:short-profile-comparison}. Now fix a profile pair in the
support of \(w\). For a side with size \(n_0\in[k,C_0k]\) and short
total \(u_0\), complete the short cycles by one cycle of length
\(n_0-u_0\ge3k/4>N\). This leaves the profile through \(N\) unchanged.
By \eqref{eq:profile} the ratio of the
probability of the completion to the Poisson probability of
the short profile is
\[
 \frac{\theta_c e^{\theta_c H_N}}{(n_0-u_0)h_{n_0}}.
\]
This is bounded below uniformly, since \(N\asymp k\) and
\(h_{n_0}\asymp n_0^{\theta_c-1}\).

If \(k+d\le C_0k\), use the one-cycle completion on both sides.
An invariant \(k\)-set then uses exactly one long cycle. The short
copies alone have total length at most \(k/4\), and the long
cycles together have length greater than \(k\). The left long
cycle leaves a short total of \(u=S_N\). The right long cycle
leaves \(k-(k+d-u')=u'-d\). Choosing this total is equivalent,
by complementation among the short copies, to choosing
\((u+u')-(u'-d)=S_N+d\). When \(u'<d\), we have \(u'-d<0\) and \(S_N+d>u+u'\),
so the shifted count is zero. Thus \(I_{n,k}=V_{N,d}\).

If \(k+d>C_0k\), then \(d>k\), since \(C_0\ge2\).
Use the one-cycle completion on the left.
On the right require the specified short profile through \(N\)
and zero counts at lengths \(N<j\le k\). For sufficiently large \(k\), the right-hand degree \(k+d\ge k\)
exceeds the threshold in Lemma~\ref{lem:short-profile-comparison}.
We also have \(u'\le(k+d)/2\) and \(k\le c_{\theta_c}(k+d)/2\) by the
choice of \(C_0\). The lemma at cutoff \(k\) therefore bounds the
probability of this event below by a fixed multiple of its Poisson
probability through \(k\). The probability
of the additional zero counts is \(e^{-\theta_c(H_k-H_N)}\),
bounded away from zero. A \(k\)-set must use the left long cycle, since the short
cycles have total length at most \(k/4\) and all remaining right
cycles exceed \(k\). Thus \(I_{n,k}=V_N\). Since \(T_N^{(2)}\le k/4\) and \(d>k\),
\(S_N+d>T_N^{(2)}\), the shifted count is zero, and \(I_{n,k}=V_{N,d}\).

In either range, distinct short profile pairs specify disjoint
events and the count on each event is \(V_{N,d}\). The probability
of each event is at least a fixed multiple of the Poisson
probability of its pair. Weighting by \(w(\mathcal C)\), summing,
and applying
\eqref{s3:critical-poisson-lower} gives
\[
 \E_A I_{n,k}^{-\omega}
 \gg \E[V_{N,d}^{-\omega}w(\mathcal C)]
 \gg_\omega (\log k)^{-1/2}.
\]
For each bounded \(k\), equations \eqref{s3:first}--\eqref{s3:bias}
and \(h_j\asymp(j+1)^{\theta_c-1}\) give
\(\mathbb E_A I_{n,k}\ll_k1\) uniformly in \(n\ge2k\).
Jensen's inequality extends the lower bound, in the stated form, to all \(k\ge1\).
Equation~\eqref{s3:existence} gives the assertion for the ratio
of the existence probability to the first moment.
\end{proof}

\subsection{Poisson sumsets and fixed-set probabilities}
\label{sec:single-profile}

In a single Poisson profile of parameter \(\theta\), let \(J_L\) be the
number of copies of lengths at most \(L\), \(T_L\) their total length,
and \(\Sigma_L\) the set of totals attainable from these copies.
A short selection of total \(s\) can be completed to size \(k\) by a
cycle of length \(k-s\). Thus distinct values of \(s\) give distinct
opportunities for completion, even when a given \(s\) has many
representations. We first make this observation uniform in the degree.

\begin{lemma}[Completing a short selection to size \(k\)]\label{lem:single-profile-completion}
Fix \(\theta>0\). There is an \(\epsilon_\theta>0\) such that for
all sufficiently large \(k\), all \(n\ge2k\), and all
\(1\le L\le\epsilon_\theta k\), the following holds. Given an
Ewens profile through \(L\) with total length \(u\le k/8\), write
\(\Sigma_L\) for its attainable totals. The conditional probability
of an invariant \(k\)-set is \(\gg_\theta |\Sigma_L|/k\).
\end{lemma}

\begin{proof}
Choose \(\epsilon_\theta\le\min(1/8,c_\theta/2)\), where \(c_\theta\)
is from Lemma~\ref{lem:short-profile-comparison}. After we condition
on the short profile, the weights of the remaining cycles are
normalized by \(g_{n-u,L}\).
The candidate lengths
\(\mathcal C=\{k-s:s\in\Sigma_L\}\) are distinct and lie in
\([7k/8,k]\). They exceed \(L\), and any one of them completes
a short selection to an invariant \(k\)-set.

Let \(X_{\mathcal C}\) count the cycles of lengths in \(\mathcal C\) in the
remaining Ewens profile. Its conditional first moment is
\[
 \E_\theta[X_{\mathcal C}\mid\text{short profile}]
 =\sum_{\ell\in\mathcal C}
  \frac{\theta}{\ell}
  \frac{g_{n-u-\ell,L}}{g_{n-u,L}}
 \asymp_\theta \frac{|\mathcal C|}{k}.
\]
Since \(u\le k/8\le n/16\) and \(\ell\le k\le n/2\), we have \(n-u-\ell\ge7n/16\). The lower bound from
Lemma~\ref{lem:short-profile-comparison} and the upper bound
\eqref{s3:coef} at mark zero are comparable for both coefficients.
The cutoffs satisfy the comparison hypothesis by the choice of
\(\epsilon_\theta\). For \(n\le4k\), each candidate cycle has length at least \(7k/8\),
so \(X_{\mathcal C}\le\lfloor 4k/(7k/8)\rfloor=4\). If \(n>4k\),
the analogous conditional factorial second moment is
\[
 \E_\theta[X_{\mathcal C}(X_{\mathcal C}-1)\mid\text{short profile}]
 =\sum_{\ell,\ell'\in\mathcal C}
  \frac{\theta^2}{\ell\ell'}
  \frac{g_{n-u-\ell-\ell',L}}{g_{n-u,L}}
 \ll_\theta (|\mathcal C|/k)^2 .
\]
Here \(u+\ell+\ell'\le17k/8<17n/32\), so \(n-u-\ell-\ell'>15n/32\) and the same coefficient bounds
apply. The diagonal uses two distinct copies. The factorial moment bound and
\(1\le|\mathcal C|\le u+1\le k\) give
\(\E_\theta[X_{\mathcal C}^2\mid\text{short profile}]
\ll_\theta |\mathcal C|/k\) when \(n>4k\).
For \(n\le4k\) the same bound follows from \(X_{\mathcal C}\le4\)
and the first moment. Cauchy--Schwarz now gives the stated conditional
probability, since \(|\mathcal C|=|\Sigma_L|\).
\end{proof}

The following comparison at \(\theta=1\) is due to Eberhard, Ford,
and Green \cite[Proposition~1.3]{efg}.

\begin{proposition}[Comparison with a Poisson sumset]\label{cor:finite-support-comparison}
For every fixed \(\theta>0\), uniformly over \(1\le k\le n/2\),
\[
 i_\theta(n,k)\asymp_\theta \frac{\E_\theta|\Sigma_k|}{k}.
\]
\end{proposition}

\begin{proof}
\emph{Upper bound.}
We compare the largest cycle on each side of a split, as in the proof
of \cite[Proposition~1.3]{efg}. A permutation with an
invariant \(k\)-set has a split of its cycles into parts of sizes
\(k\) and \(n-k\). If either part consists of fixed points, there are at least \(k\)
fixed points. By \eqref{s3:coef} at cutoff \(1\) and mark \(2\),
this has probability \(O_\theta(2^{-k})\).
Otherwise each part contains a cycle of length at least two. Let
\(\ell\) be the smaller of the maximum cycle lengths on the two sides.
Then \(2\le\ell\le k\), so \(k\ge2\) and \(\ell\in(L,2L]\) for some
dyadic \(L<k\). Both parts contain a cycle longer than \(L\), and one
of them has no cycle longer than \(2L\).

For the capped part put
\[
 a_{t,L}=[z^t]\exp\left(\theta\sum_{L<j\le2L}\frac{z^j}{j}\right).
\]
The generating function at \(z=e^{1/L}\) is bounded in terms of
\(\theta\). Logarithmic differentiation gives
\(t a_{t,L}=\theta\sum_{L<j\le\min(2L,t)}a_{t-j,L}\).
Multiplying by \(e^{t/L}\), using \(e^{j/L}\le e^2\), and
summing the weighted coefficients gives
\(t a_{t,L}e^{t/L}\ll_\theta1\).
For the uncapped part use \eqref{s3:coef} at mark zero. For
\(t>L\) the resulting bounds are
\[
 g_{t,L}\ll_\theta t^{\theta-1}L^{-\theta},
 \qquad
 a_{t,L}\ll_\theta t^{\theta-1}L^{-\theta}e^{-t/L},
\]
where the second uses \(t^{-1}\le t^{\theta-1}L^{-\theta}\) for
\(t>L\) and \(\theta>0\).

Fix a short profile through \(L\) of total length \(u\).
For an attainable short total \(s\in\Sigma_L\) on the first
side, the long parts have sizes \(t=k-s\) and
\(t'=n-k-u+s\). For a long profile \(\mathbf m=(m_j)_{j>L}\)
write \(w(\mathbf m)=\prod_{j>L}(\theta/j)^{m_j}/m_j!\).
For each split \(\mathbf m^{(1)}+\mathbf m^{(2)}=\mathbf m\),
\[
 w(\mathbf m^{(1)})w(\mathbf m^{(2)})
   =w(\mathbf m)\prod_{j>L}\binom{m_j}{m_j^{(1)}}
   \ge w(\mathbf m).
\]
Thus the weight of the long profiles admitting a witnessing split
at this cutoff and short total is at most \(a_{t,L}g_{t',L}+g_{t,L}a_{t',L}\),
provided \(t,t'>L\).

Divide by \(h_n\) and convert the short weights to Poisson
probabilities by multiplying by \(e^{\theta H_L}\).
Since \(n-k\ge n/2\), the resulting bound per attainable total
is \(O_\theta(k^{\theta-1}L^{-\theta})\) times
\[
 \left(\frac{t}{k}\right)^{\theta-1}
 \left(\frac{t'}{n-k}\right)^{\theta-1}
 (e^{-t/L}+e^{-t'/L}).
\]
For \(u\le k/2\), both ratios are at least \(1/2\) and each
exponential is at most \(e^{-k/(2L)}\). If \(u>k/2\), the
product of the two powers is at most one for \(\theta\ge1\), since
\(t\le k\) and \(t'\le n-k\). For \(0<\theta<1\), use
\(k\le t+u\), \(n-k\le t'+u\), and \(t,t'>L\) to bound the product
by \((1+u/L)^{2(1-\theta)}\ll_\theta e^{u/(2L)}\).
As \(u>k/2\) gives \(e^{u/(2L)}\le e^{u/L-k/(4L)}\), the display is
\(O_\theta(e^{u/L-k/(4L)})\) for every \(u\).

Weighting the short profile by \(e^{T_L/L}\) changes the Poisson
mean at \(j\le L\) from \(\theta/j\) to \((\theta/j)e^{j/L}\).
The normalizing factor is \(e^{\Lambda_L}\), where
\[
 \Lambda_L=\theta\sum_{j\le L}\frac{e^{j/L}-1}{j}
 \le\theta(e-1).
\]
The tilted profile is obtained by adding an independent family whose
total count \(D_L\) has law \(\operatorname{Pois}(\Lambda_L)\).
Each added copy at most doubles the support, and
\(\E 2^{D_L}=e^{\Lambda_L}\). Hence
\begin{equation}\label{eq:tilted-poisson-support}
 \E_\theta[e^{T_L/L}|\Sigma_L|]\le e^{2\Lambda_L}\E_\theta|\Sigma_L|.
\end{equation}
Summing first over \(s\in\Sigma_L\), then over the dyadic cutoffs \(L<k\), gives
\[
 i_\theta(n,k)\ll_\theta 2^{-k}+\frac1k
 \sum_{\substack{L<k\\L\ {\rm dyadic}}}
 (k/L)^\theta e^{-k/(4L)}\E_\theta|\Sigma_L|.
\]
We have \(\Sigma_L\subseteq\Sigma_k\), and the sum of
\((k/L)^\theta e^{-k/(4L)}\) over these cutoffs is bounded in
terms of \(\theta\). Since \(0\in\Sigma_k\), this also absorbs the
term \(2^{-k}\).

\emph{Lower bound.}
Choose \(0<\epsilon\le\min(\epsilon_\theta,c_\theta/2)\) so small that
\(e^{-1/(8\epsilon)+2\theta(e-1)}\le1/2\), and take
\(L=\lfloor\epsilon k\rfloor\) for sufficiently large \(k\). By
\eqref{eq:tilted-poisson-support},
\[
 \E_\theta[|\Sigma_L|\ind_{\{T_L>k/8\}}]
 \le e^{-k/(8L)+2\Lambda_L}\E_\theta|\Sigma_L|
 \le\tfrac12\E_\theta|\Sigma_L|.
\]
For each short profile with \(T_L\le k/8\), Lemma~\ref{lem:short-profile-comparison}
bounds its Ewens probability below by a constant times its Poisson probability.
Lemma~\ref{lem:single-profile-completion} then bounds the conditional
probability of fixing a \(k\)-set below by a constant times
\(|\Sigma_L|/k\). Summing over these profiles gives
\(i_\theta(n,k)\gg_\theta k^{-1}\E_\theta|\Sigma_L|\).
The number of copies of lengths in \((L,k]\) is independent of the
profile through \(L\) and is Poisson with mean \(\theta(H_k-H_L)\). Each further copy
at most doubles the support, so
\[
 \E_\theta|\Sigma_k|
 \le e^{\theta(H_k-H_L)}\E_\theta|\Sigma_L|
 \ll_\theta\E_\theta|\Sigma_L|,
\]
where the last inequality uses \(L\asymp_\theta k\).
For each remaining \(k\), the probability of a \(k\)-cycle tends to
\(1-e^{-\theta/k}>0\) as \(n\to\infty\) by \eqref{s2:feller}. In each
of the finitely many smaller degrees \(n\ge2k\) the identity has
positive Ewens probability. Adjusting the constant gives the lower
bound for these \(k\) as well.
\end{proof}

For \(\theta>\theta_c\) it remains to estimate \(\E_\theta|\Sigma_k|\).
The mean number of copies in a dyadic block tends to \(\theta\log 2\).
For \(\theta_c<\theta<\lambda_0\), changing the block means to one
accounts for the power of \(k\). The logarithmic factor comes from a
centered walk that stays nonnegative and returns to zero. At
\(\theta=\lambda_0\) the endpoint is free, and the probability of
staying nonnegative accounts for the exponent \(-1/2\).
For the corresponding transition for continuum subset sums, see
\cite[Section~22]{ArratiaScale}.

\begin{proposition}[Fixed-set probabilities above the exchange threshold]\label{prop:large-theta}
Uniformly over \(1\le k\le n/2\),
\begin{equation}\label{eq:single-profile-intermediate}
 i_\theta(n,k)\asymp_\theta
 k^{\theta-1-D(\theta)}(1+\log k)^{-3/2}
 \qquad(\theta_c<\theta<\lambda_0).
\end{equation}
At and above \(\lambda_0\) one has
\[
 i_{\lambda_0}(n,k)\asymp (1+\log k)^{-1/2},
 \qquad
 i_\theta(n,k)\asymp_\theta 1\quad(\theta>\lambda_0).
\]
\end{proposition}

The dyadic block means at
\(\theta=\lambda_0\) coincide with those of the merged pair at
\(\theta=\theta_c\). Write \(L=2^b\),
where \(b\ge0\), and set \(B_v=J_{2^v}-J_{2^{v-1}}\) for
\(1\le v\le b\). When \(\theta_c<\theta\le\lambda_0\), put
\[
 \rho=\frac{\theta}{\lambda_0},\qquad
 \delta_\theta=\lambda_0(\rho-1-\log\rho)=D(\theta)+1-\theta .
\]
The means of the \(B_v\) are \(\rho\mu_v\), where
\(\mu_v=\lambda_0(H_{2^v}-H_{2^{v-1}})<1\) and
\(1-\mu_v=O(2^{-v})\). Under the reference law \(\Pref\), with expectation \(\E_{\!*}\), the
block counts are independent \(Y_v\sim\operatorname{Pois}(1)\).
The law of \(J_1\) and the conditional length laws are the same
under \(\Pref\) and \(\mathbb P_\theta\). Put
\(L_i^*=\sum_{v=1}^i(Y_v-1)\), with \(L_0^*=0\). The likelihood ratio is
\begin{equation}\label{eq:single-profile-likelihood}
 \exp\!\left(\sum_{v=1}^b(1-\rho\mu_v)\right)
       \prod_{v=1}^b(\rho\mu_v)^{Y_v}
 =L^{-\delta_\theta}\rho^{L_b^*}
   \exp\!\left(\rho\sum_{v=1}^b(1-\mu_v)\right)
       \prod_{v=1}^b\mu_v^{Y_v}.
\end{equation}
Here \(\sum_{v=1}^bY_v=b+L_b^*\) and
\(e^{b(1-\rho+\log\rho)}=L^{-\delta_\theta}\).

\begin{lemma}[Expected Poisson sumset size]\label{lem:single-support-upper}
For \(\theta_c<\theta\le\lambda_0\) and \(L=2^b\),
\begin{equation}\label{eq:single-support-upper}
 \E_\theta|\Sigma_L|\ll_\theta L^{1-\delta_\theta}
 \begin{cases}
 (1+\log L)^{-3/2},&\theta_c<\theta<\lambda_0,\\
 (1+\log L)^{-1/2},&\theta=\lambda_0 .
 \end{cases}
\end{equation}
\end{lemma}

\begin{proof}
Put \(S_i=\sum_{v=1}^i(B_v-1)\). The support through \(L\) is a
union of at most \(2^{\sum_{v>i}B_v}\) translates of \(\Sigma_{2^i}\),
as in \cite[Lemma~3.1]{efg}. Since \(|\Sigma_{2^i}|\le T_{2^i}+1\),
for \(0\le i\le b\),
\begin{equation}\label{eq:single-support-split}
 \frac{|\Sigma_L|}{L}
 \le \frac{T_{2^i}+1}{2^i}\,2^{S_b-S_i}
 \le 2^{S_b-S_i}
       \left(2^{-i}(J_1+1)+\sum_{v=1}^i2^{v-i}B_v\right).
\end{equation}

Let \(\tau\) be the first maximizer of the block walk over
\(\{0,\ldots,b\}\), and apply \eqref{eq:single-support-split} at
\(i=\tau\). Under \(\Pref\), on \(\{\tau=i\}\) write
\(u=L_i^*\) and \(a=L_i^*-L_b^*\). For \(i>0\) the reversed prefix
has partial sums \(L_i^*-L_{i-j}^*\ge1\) for \(1\le j\le i\) and
ends at \(u\ge1\). The independent suffix has partial sums
\(L_{i+j}^*-L_i^*\le0\) for \(1\le j\le b-i\) and ends at \(-a\).
At \(i=0\) the prefix is empty and \(u=0\). The bridge estimate
\eqref{s3:critical-walk-bridge} gives
\[
 \Pref(\tau=i,L_i^*=u,L_b^*=u-a)
 \ll \frac{(u+1)(a+1)}{(i+1)^{3/2}(b-i+1)^{3/2}} .
\]
The same bound holds when the indicator of this event is weighted
by \(Y_v\), \(1\le v\le i\). By Lemma~\ref{lem:poisson-insertion}, the
reversed prefix before insertion must stay nonnegative and end at
\(u-1\). Its bridge bound and the original suffix bound apply.

In \eqref{eq:single-profile-likelihood} the product of the
\(\mu_v^{Y_v}\) is at most one and the exponential factor is bounded
uniformly in \(b\). The remaining weight \(\rho^{L_b^*}\) times
the factor \(2^{-a}\) in \eqref{eq:single-support-split} is
\(\rho^u(2\rho)^{-a}\). When \(1/2<\rho<1\), both factors decay
geometrically in their respective variables. The sum over \(u,a,v\)
at fixed \(i\) is bounded by a constant depending only on \(\theta\),
since \(\sum_{v\le i}2^{v-i}\le2\) and \(J_1\) is independent of the
block counts, with mean \(\theta\). The sum of the two bridge factors
over \(i\) is \(O((b+1)^{-3/2})\) by splitting at \(b/2\).
At \(\rho=1\) the height has weight one. Leave the prefix endpoint
unspecified and use the survival bound \eqref{s3:critical-walk-survival},
also after insertion. The sum over \(i\) is then
\[
 \sum_{i=0}^b (i+1)^{-1/2}(b-i+1)^{-3/2}
 \ll (b+1)^{-1/2}.
\]
Restoring \(L^{-\delta_\theta}\) gives the claim.
\end{proof}

\begin{lemma}[Collision probability]\label{lem:single-collision}
For \(L=2^b\), \(b\ge1\), omit the copies of length one and condition on the
counts \(y_v\) in the blocks \((2^{v-1},2^v]\), \(1\le v\le b\).
Let \(R_i=\sum_{v=b-i+1}^b(y_v-1)\) for \(0\le i\le b\).
Given the lengths, let \(q_L\) be the probability that two independent
fair selections of these copies have the same total. Then
\begin{equation}\label{eq:single-collision-suffix}
 \E[q_L\mid(y_v)]\ll L^{-1}\sum_{i=0}^b2^{-R_i}.
\end{equation}
\end{lemma}

\begin{proof}
The lengths of the \(y_v\) copies in block \(v\) are independent
draws with probabilities proportional to \(1/j\) on that block.
The probability that the selections agree on blocks above \(v\) is
\(2^{-\sum_{w>v}y_w}\). On the event that \(v\) is their last block of
disagreement, condition on the selections and expose all lengths
except that of one differing copy in block \(v\). Its coefficient
in the difference of the totals is \(1\) or \(-1\), so equality
fixes its length. The largest point mass in that block is
\(O(2^{-v})\).
The contribution of identical selections is
\(2^{-\sum_vy_v}=L^{-1}2^{-R_b}\). For differing selections
put \(i=b-v\), so \(2^{-v}2^{-\sum_{w>v}y_w}=L^{-1}2^{-R_i}\).
Summing over the last block of disagreement gives
\eqref{eq:single-collision-suffix}; see also \cite[Lemma~4.1]{efg}.
\end{proof}

\begin{lemma}[Poisson profiles with many attainable totals]\label{lem:single-support-lower}
For each fixed \(\theta>\theta_c\) there is a \(C=C(\theta)\) such
that at dyadic \(L=2^b\), \(b\ge1\),
\[
 \Pp_\theta\{T_L\leq CL,\ \lvert\Sigma_L\rvert\geq L/C\}
 \gg_\theta
 \begin{cases}
 L^{-\delta_\theta}(b+1)^{-3/2},&\theta_c<\theta<\lambda_0,\\
 (b+1)^{-1/2},&\theta=\lambda_0,\\
 1,&\theta>\lambda_0 .
 \end{cases}
\]
\end{lemma}

\begin{proof}
If the attainable totals of a fair selection have masses \(p_s\),
then \(q_L=\sum_s p_s^2\). Cauchy--Schwarz gives
\(\lvert\Sigma_L\rvert\ge q_L^{-1}\). We use \(T_L\) and \(\Sigma_L\) for lengths
\(2,\ldots,L\) throughout the proof and require no copies of length
one at the end.

\emph{Nonnegative walks.}
Suppose \(\theta_c<\theta\le\lambda_0\), and under \(\Pref\)
set \(y_v=Y_v\). The reversed walk \(R_i\) has the same law as \(L_i^*\).
For \(\theta=\lambda_0\) take the unpinned event
\(G=\{\min_{0\le i\le b}R_i\ge0\}\).
By \eqref{eq:critical-zero-survival},
\(\Pref(G)\asymp(b+1)^{-1/2}\).

For \(\theta_c<\theta<\lambda_0\), pin the endpoint at zero by taking
\(G=\{\min_{0\le i\le b}R_i\ge0,\ R_b=0\}\).
Reversing a bridge from zero to zero interchanges its nonnegative
and nonpositive versions. Equation~\eqref{eq:critical-ballot} at
\(a=0\) therefore gives
\[
 \Pref(G)=\frac{e}{b+1}\,\mathbb P(\operatorname{Pois}(b+1)=b)
 \asymp(b+1)^{-3/2}
\]
by Stirling's formula.

For the unpinned walk, a nonnegative bridge to \(a\) at time \(i\)
and a suffix that stays at least \(-a\) relative to its start give
\[
 \Pref(G,R_i=a)\ll
 \frac{(a+1)^2}{(i+1)^{3/2}(b-i+1)^{1/2}}\qquad(a\ge0).
\]
For the walk required to return to zero, the reversed suffix has
partial sums \(R_b-R_{b-j}\le0\) for \(1\le j\le b-i\)
and ends at \(R_b-R_i=-a\). The second bridge bound in
\eqref{s3:critical-walk-bridge} replaces the power \(1/2\) in
the denominator by \(3/2\).
Summing against \(2^{-a}\) and splitting the sum over \(i\) at
\(b/2\) gives
\(\E_{\!*}[q_L\ind_G]\ll L^{-1}\Pref(G)\) in both cases.

For \(1\le v\le b\), Lemma~\ref{lem:poisson-insertion} gives
\(\E_{\!*}[Y_v\ind_G]\ll\Pref(G)\). In the unpinned case the walk
before insertion must stay at or above \(-1\), so
\eqref{s3:critical-walk-survival} applies. In the pinned case it must
also end at \(-1\). Prepending an independent increment equal to
\(+1\), which has probability \(e^{-1}/2\), yields a nonnegative
bridge to zero of length \(b+1\). Its probability is
\(O((b+1)^{-3/2})\) by \eqref{s3:critical-walk-bridge}.
Consequently
\(\E_{\!*}[T_L\ind_G]\le\sum_v2^v\E_{\!*}[Y_v\ind_G]
\ll L\Pref(G)\). For a sufficiently large constant \(C\),
\[
 E=G\cap\{q_L\le C/L,\ T_L\le CL\}
 \qquad\text{satisfies}\qquad
 \Pref(E)\gg\Pref(G),\quad
 \E_{\!*}[Y_v\mid E]\ll1.
\]

In \eqref{eq:single-profile-likelihood} the factor
\(\rho^{L_b^*}\) equals one on \(E\). This follows from
\(\rho=1\) in the unpinned case and \(R_b=L_b^*=0\) in the
pinned case. Also
\(\sum_v-\log\mu_v<\infty\). The conditional laws of the lengths
agree under the two measures, so conditional Jensen gives
\[
 \Pp_\theta(E)\ge L^{-\delta_\theta}\Pref(E)
       \exp\!\left(\sum_{v=1}^b
                (\log\mu_v)\E_{\!*}[Y_v\mid E]\right)
 \gg L^{-\delta_\theta}\Pref(G).
\]

\emph{Positive drift.}
Now fix \(\theta>\lambda_0\) and use the Poisson block counts \(y_v=B_v\).
Write \(\mu=\theta\log 2>1\). By the harmonic sum estimate, the number
of copies in the last \(i\) blocks is Poisson of mean
\(\theta(H_{2^b}-H_{2^{b-i}})=\mu i+O_\theta(1)\), uniformly in
\(0\leq i\leq b\). Choose \(s\in(0,1)\) so small that
\[
 \rho_s=2^s\exp\{\mu(2^{-s}-1)\}<1 .
\]
Indeed the derivative of \(\log\rho_s\) at zero is
\((1-\mu)\log2<0\). The Poisson exponential moment gives
\(\E 2^{-sR_i}\ll_\theta\rho_s^i\), and subadditivity gives
\[
 \E\!\left(\sum_{i=0}^b2^{-R_i}\right)^s
 \leq\sum_{i=0}^b\E 2^{-sR_i}\ll_\theta1 .
\]
Thus a block-count event of probability at least \(3/4\) has
\(\sum_i2^{-R_i}\leq C_\theta\). On this event the conditional mean of
\(q_L\) is \(O_\theta(1/L)\) by \eqref{eq:single-collision-suffix}.
Also \(\E T_L\leq\theta L\). Markov's inequality first
conditionally for \(q_L\) and then for \(T_L\) shows that, for a
sufficiently large \(C=C(\theta)\), the inequalities \(q_L\leq C/L\)
and \(T_L\leq CL\) hold together with probability at least \(1/2\).

Requiring zero copies of length one costs the independent factor
\(e^{-\theta}\) and leaves the two statistics unchanged. This proves
the bound for the full profile.
\end{proof}

\begin{proof}[Proof of Proposition~\ref{prop:large-theta}]
At dyadic \(L\), multiply the probability in
Lemma~\ref{lem:single-support-lower} by \(L/C\) to get a lower
bound for \(\E_\theta|\Sigma_L|\). For \(\theta_c<\theta\le\lambda_0\)
this has the same order as the upper bound in
Lemma~\ref{lem:single-support-upper}. For \(\theta>\lambda_0\) the lower
bound has order \(L\), while \(|\Sigma_L|\le1+T_L\) and
\(\E_\theta T_L=\theta L\). Since \(|\Sigma_L|\) increases with \(L\),
the estimates at the adjacent powers of two give the same orders at
every integer \(k\ge2\). Proposition~\ref{cor:finite-support-comparison}
now gives the three stated probabilities, using
\(-\delta_\theta=\theta-1-D(\theta)\) in \eqref{eq:single-profile-intermediate}.
At \(k=1\) the comparison also applies, since \(\E_\theta|\Sigma_1|=1+\theta\).
\end{proof}

\begin{proof}[Proof of Theorem~B]
For \(0<\theta<\theta_c\), Theorem~\ref{s3:main} gives the total
variation estimate and the approximation to \(i_\theta(n,k)\).
The identity \(W_0=2W_\infty\) gives \(c_0=c_\infty/2\), and
\eqref{s3:d-tail} gives \(c_d\to c_\infty\) by boundedness of the
reciprocals. At \(\theta_c\), take \(\omega=1\) in
Proposition~\ref{s3:critical-rate} and use \eqref{s3:existence}.
Proposition~\ref{prop:large-theta} proves the remaining estimates.
\end{proof}

\section{Comparison with the earlier generation bounds}\label{sec:earlier-generation}

Brito, Fowler, Junge, and Levy \cite{brito} introduced the Poisson
sumset threshold and the generation threshold under Ewens sampling.
Their staircase \cite[Eq.~(1.4)]{brito}, denoted there by \(h(\alpha)\), is
\[
 g_\theta=\begin{cases}
 \lceil(1-\theta\log 2)^{-1}\rceil,&0<\theta<\lambda_0,\\
 \infty,&\theta\ge\lambda_0.
 \end{cases}
\]
Table~\ref{tab:earlier-results} gives the conclusions for the Poisson
model, the Ewens law, and the more general conditioned laws in
\cite{brito}.

\begin{table}[htbp]
\centering\small
\caption{Scope of the threshold and generation results in \cite{brito}.}
\label{tab:earlier-results}
\renewcommand{\arraystretch}{1.12}
\begin{tabular}{@{}p{.28\linewidth}@{\hspace{.025\linewidth}}p{.685\linewidth}@{}}
\toprule
Location in \cite{brito} & Conclusion\\
\midrule
\raggedright Theorem~1.1, Eq.~(1.2), \mbox{pp.\ 855--856}
 & A finite change to the defining counts can force a common fixed point.
 The printed hypotheses also allow failure of fixed-length convergence
 (Propositions~\ref{prop:strong-log-counterexamples}
 and~\ref{prop:lacunary-logarithmic}).\\[2pt]
\raggedright Eq.~(1.3), p.~856
 & For Ewens sampling, the infimum over all degrees gives
 \(m_\theta^{\mathrm{all}}=\max\{3,m_\theta\}\)
 (Proposition~\ref{prop:all-degree-threshold}).\\[2pt]
\raggedright Theorem~1.3, p.~857
 & The transitive-overgroup conclusion holds for \(\theta\le\lambda_0\).
 For \(\theta>\lambda_0\), an invariant half is an obstruction in even
 degrees (Lemma~\ref{lem:ewens-overgroups} and
 Propositions~\ref{prop:endpoint-overgroups} and~\ref{prop:overgroup-counterexample}).\\[2pt]
\raggedright Lemma~3.10, Corollary~3.13, \mbox{pp.\ 877--881}
 & The unconditional and odd-conditioned estimates for a fixed set size
 fail at \(\theta=1\). Lemma~\ref{lem:common-size-tail} bounds the
 common sizes needed for generation.\\[2pt]
\raggedright Propositions~4.2 and~5.1, p.~883
 & The threshold bounds hold. At a finite jump the summability argument
 uses \(g_\theta+1\) profiles. At continuity
 points below \(\lambda_0\), the intersection at \(g_\theta\) is
 almost surely finite (Lemma~\ref{lem:common-size-tail}).\\[2pt]
\raggedright Question~2, p.~860
 & Joint domination from above on every fixed window allows a common
 invariant size with probability tending to one, even with strictly
 smaller limiting Poisson means (Proposition~\ref{prop:short-domination}).\\
\bottomrule
\end{tabular}
\end{table}

\subsection{Poisson intersections}\label{subsec:earlier-poisson}

For \(0<\theta<\lambda_0\), Proposition~5.1 of \cite[p.~883]{brito}
asserts both \(s_\theta\ge g_\theta\) and almost-sure infinitude of the
\(g_\theta\)-fold intersection. For the threshold lower bound it
suffices to take \(m<g_\theta\), when
\(m(1-\theta\log2)<1\) and Proposition~\ref{thm:strict-converse}
applies. At a continuity point, however,
\(g_\theta(1-\theta\log2)>1\), and
Lemma~\ref{lem:common-size-tail} gives almost-sure finiteness at the
ceiling. Eberhard, Ford, and Green \cite[Corollary~2.5]{efg-invariable}
already proved this when \(\theta=1\) and \(g_1=4\).

At \(\theta=1/4\) the distinction follows from a first moment.
Here \(g_{1/4}=2\) and \(\E U_j=h_j(1/4)\). By independence and
Gauss's summation formula for the hypergeometric series,
\[
 \E W_\infty=\sum_{j\ge0} h_j(1/4)^2
 ={}_2F_1\!\left(\tfrac14,\tfrac14;1;1\right)
 =\frac{\Gamma(1/2)}{\Gamma(3/4)^2}.
\]
Every common positive sum contributes at least one selection pair
in addition to the empty pair. Thus there are almost surely finitely
many common sums, and
\(p_{0,2}(1/4)\ge 2-\E W_\infty=0.8196594\ldots\).

The two strict inequalities also explain the changes to the earlier
proofs. The Fourier argument in \cite[pp.~884--889]{brito} requires a
\(\beta<1\) with \(\beta\theta\log2>1-1/m\), hence
\(m(1-\theta\log2)<1\). In \cite[Eqs.~(5.8), (5.15)]{brito},
distance to the nearest integer gives \(\max(1,\lVert kt\rVert)=1\).
The exponent also changes from \(-1+\delta_2\) to \(-1-\delta_2\).
Proposition~\ref{thm:strict-converse} supplies the integrable factor
\(\max(1,k\lVert t\rVert)^{-1-\delta}\), \(\delta>0\), and uses a
positive completing length in each selection. For the upper bound in
\cite[Proposition~4.2]{brito}, the summability argument requires the
opposite inequality \(m(1-\theta\log2)>1\). At a finite jump
\(g_\theta=m\), the choice \(g_{\theta-2\epsilon}=m\) in that proof gives
exponent \(-1+m\epsilon\log2\) for small \(\epsilon>0\). Using
\(m+1\) profiles gives the stated bound \(s_\theta\le m+1\) by
Lemma~\ref{lem:common-size-tail}. Coupling from a slightly smaller
parameter, as on p.~890 of \cite{brito}, gives a finite intersection
there and leaves the intersection at the jump undetermined.
Theorem~\ref{s3:poisson-critical} decides the first jump. The higher
jumps remain open.

\subsection{Invariant sizes and transitive overgroups}\label{subsec:earlier-overgroups}

Brito et al.\ \cite[Theorem~1.3, p.~857]{brito} assert that, for
every fixed Ewens parameter, with probability \(1-o(1)\)
\textquotedblleft the only transitive subgroups containing \(\pi\) are
\(A_n\) and \(S_n\)\textquotedblright.
Lemma~\ref{lem:ewens-overgroups} proves the assertion for
\(\theta<\lambda_0\), and Proposition~\ref{prop:endpoint-overgroups}
treats equality. For \(\theta>\lambda_0\), Theorem~B gives an
invariant half with probability bounded away from zero in even degrees.
In degree \(2k\), the unordered pair consisting of the invariant
half and its complement has the transitive stabilizer \(S_k\wr S_2\). Eberhard, Ford, and Green \cite[Corollary~1.2]{efg} use
the same construction for uniform permutations. See also
\cite[Theorems~1.2 and~1.3]{EFK16}. At parameter four, two exact
moments give the obstruction directly.

\begin{proposition}[Transitive overgroups at parameter four]
\label{prop:overgroup-counterexample}
Let \(\pi_{2k}\) have the Ewens law of parameter \(4\).
For \(k\ge2\), the probability that \(\pi_{2k}\) lies in a
transitive subgroup that does not contain \(A_{2k}\) is at least
\[
 \frac{(k+1)^3}{\binom{2k+3}{3}},
\]
which tends to \(3/4\).
\end{proposition}

\begin{proof}
Put \(X=2^{-K_{2k}}I_{2k,k}\), and use a subscript on
the expectation to indicate the Ewens parameter. In the numerator
of the first moment the weight \(2^{-K_{2k}}\) changes the
parameter to two. In the second moment \(4^{-K_{2k}}\) changes it
to one. For an ordered pair of invariant halves meeting in
\(a\) points, the four membership cells have sizes
\(a,k-a,k-a,a\). Counting their labelings cancels their factorials
in the Ewens weights. With a single half there are two cells, each
of size \(k\). Thus
\[
 \E_4 X=\frac{h_k(2)^2}{h_{2k}(4)}
       =\frac{(k+1)^2}{h_{2k}(4)},\qquad
 \E_4 X^2=\frac{\sum_{a=0}^{k} h_a(1)^2h_{k-a}(1)^2}{h_{2k}(4)}
         =\frac{k+1}{h_{2k}(4)}.
\]
Since \(h_{2k}(4)=\binom{2k+3}{3}\), Cauchy--Schwarz gives
\[
 \Pp_4(I_{2k,k}>0)\ge\frac{(\E_4X)^2}{\E_4X^2}
 =\frac{(k+1)^3}{\binom{2k+3}{3}}\longrightarrow\frac34 .
\]
An invariant half and its complement form a system of two blocks.
Its stabilizer \(S_k\wr S_2\) is transitive and does not contain
\(A_{2k}\) for \(k\ge2\), since its index in \(S_{2k}\) is
\(\binom{2k}{k}/2\ge3\). At \(k=1\) the stabilizer is \(S_2\)
itself.
\end{proof}

\paragraph{A fixed size and a common size.}
Lemma~3.10 of \cite[pp.~877--878]{brito} asserts the unconditional
bound \(O_{\theta,\epsilon}(k^{-1+(\theta+\epsilon)\log2})\)
uniformly for \(k\le n/2\). The Poisson estimate it uses,
\cite[Lemma~4.1]{brito}, also imposes a restriction on the short
cycle counts. At \(\theta=1\), \(\epsilon=1/10\), the claimed
decay exponent is \(1-1.1\log2=0.237538\ldots\), whereas the
uniform fixed-set probability has order
\(k^{-D(1)}(1+\log k)^{-3/2}\) with \(D(1)=0.086071\ldots\)
by \cite[Theorem~1.1]{efg}. The same contradiction applies to the
odd-conditioned estimate in \cite[Corollary~3.13]{brito}. For fixed
\(k\), even and odd uniform permutations have the same limiting
short counts \cite[Lemma~2.6]{efg-invariable}, and one may then take
\(k\to\infty\).

There is a related dependence in \cite[Eq.~(3.36)]{brito}. The
terminal Feller spacing \(J_n^{\rm term}\), denoted \(J_n\) there,
and the infinite spacing sumset \(\mathcal L\) are read from the
same Bernoulli sequence. For \(n=4\), \(\theta=1\), the event
\(J_4^{\rm term}=1\) is \(\xi_4=1\). Given that event, any one
of \(\xi_2,\xi_3,\xi_5\) equal to one gives a spacing of length
one, since \(\xi_1=1\). Independence of these three Bernoulli
variables yields
\[
 \Pp(1\in\mathcal L\mid J_4^{\rm term}=1)
 \ge 1-(1-1/2)(1-1/3)(1-1/5)=11/15
 >1-e^{-1}=\Pp(1\in\mathcal L).
\]
Thus multiplying the marginals underestimates this term. For
generation, Lemma~\ref{lem:common-size-tail} instead bounds the
probability of a \emph{common} size above \(M\), also after
conditioning one sample to be odd. At the same cutoff, the sign
and short counts are asymptotically independent by the proof of
Theorem~\ref{thm:strict-generation}. Taking \(n\to\infty\) and
then \(M\to\infty\) gives the strict-range generation limit.

\subsection{Strong logarithmic laws and the degree convention}
\label{subsec:earlier-logarithmic}

Brito et al.\ start with independent nonnegative integer variables
\(Z_1,Z_2,\ldots\) and condition on \(\sum_{i=1}^n iZ_i=n\).
In their original symbols, the strong \(\alpha\)-logarithmic condition
\cite[Eq.~(1.2), p.~855, offprint p.~3]{brito} is
\begin{equation}\label{eq:strong-logarithmic-condition}
 \begin{aligned}
    |i\Pp[Z_i=1]-\alpha|&<e(i)c_1,\\
    i\Pp[Z_i=\ell]&\le e(i)c_\ell,\qquad \ell\ge2,
 \end{aligned}
\end{equation}
\textquotedblleft where \(\sum i^{-1}e(i)\) and \(\sum \ell c_\ell\) are both finite.\textquotedblright{}
For the rest of this subsection write \(Q_i=Z_i\) for these defining
counts and \(\theta=\alpha\). On conditioning, each permutation
with the resulting cycle type is given the same probability.

Theorem~1.1 of \cite[p.~856]{brito} asserts that the generation
threshold is \(g_\theta\) at continuity points and belongs to
\(\{g_\theta,g_\theta+1\}\) at finite jumps. Their justification
identifies the limiting counts of each fixed length with Poisson
variables of means \(\theta/i\). Under the uniform logarithmic
condition of Arratia, Barbour, and Tavar\'e \cite[Eqs.~(1.11)--(1.12)]{abt2000},
the error bound also satisfies \(e(i)\downarrow0\).
Their Theorem~3.1 approximates the counts through \(b_n=o(n)\)
in total variation by \(Q_1,\ldots,Q_{b_n}\). Their Theorem~3.2
compares the large components with those of an Ewens permutation.
The law of any particular \(Q_i\) can therefore affect the
generation threshold.

\begin{proposition}[Finite modifications of the defining counts]
\label{prop:strong-log-counterexamples}
For each \(\theta>0\) there are independent counts satisfying
\eqref{eq:strong-logarithmic-condition} and the uniform logarithmic
condition, with feasible conditioning in every degree \(n\ge2\),
for which the limit-inferior generation threshold is infinite. This
is possible with \(\Pp(Q_1=0)>0\). There are also such counts with
\(\Pp(Q_j=0)>0\) for all \(j\) for which the all-degree generation
threshold is infinite.
\end{proposition}

\begin{proof}
Take \(Q_1\) uniform on \(\{1,2\}\) and take
\(Q_j\sim\operatorname{Pois}(\theta/j)\) for \(j\ge2\), independently.
At the Poisson coordinates one has
\[
 |j\Pp(Q_j=1)-\theta|\le\theta^2/j,\qquad
 j\Pp(Q_j=a)\le\frac{\theta^{a}}{a!j}\quad(a\ge2).
\]
These bounds give \eqref{eq:strong-logarithmic-condition} with
\(e(j)=1/j\). Choose \(c_1\) strictly larger than
\(\max\{\theta^2,|1/2-\theta|\}\) and take
\(c_2\ge\max\{\theta^2/2,1/2\}\). The conditioning event has positive
probability in every degree \(n\ge2\). Take two fixed points
for \(n=2\), and a fixed point and an \((n-1)\)-cycle for \(n\ge3\).
Since every sample has a fixed point, conjugation can align a
common fixed point for any number of samples.

For an example with \(\Pp(Q_1=0)>0\), take \(Q_2\) uniform on
\(\{1,2\}\) and leave the other counts Poisson with means
\(\theta/j\). Choose \(c_1\) strictly larger than
\(\max\{\theta^2,2|1-\theta|\}\) and take
\(c_2\ge\max\{\theta^2/2,2\}\).
Conditioning is feasible by taking one \(2\)-cycle and \(n-2\)
fixed points. For every \(n\ge3\), the \(2\)-cycle supplies a
proper invariant set that conjugation can align across the
samples.

For the all-degree claim take \(Q_2=0\) almost surely and take the
other counts Poisson with means \(\theta/j\). Here
\(\Pp(Q_j=0)>0\) for every \(j\). Choose \(c_1\) strictly larger
than \(\max\{\theta^2,2\theta\}\). Fixed points alone make
conditioning feasible in every degree. In degree two the
resulting permutation is always the identity, so the infimum
over degrees is zero for every finite number of samples.
\end{proof}

The first two examples disprove the general threshold claimed in
\cite[Theorem~1.1]{brito} whenever \(g_\theta\) is finite, even
with the limit-inferior convention. The third exhibits a separate
obstruction due to the infimum over all degrees.

The sum \(\sum_j e(j)/j\) in \eqref{eq:strong-logarithmic-condition}
allows \(e(j)\) to be large at sparse lengths. The printed condition
gives a summable comparison of the independent defining counts with the
Poisson counts, as in Remark~\ref{rem:summable-perturbation}.
Conditioning on total size can amplify a discrepancy at a length close
to the degree, as the next example shows.

\begin{proposition}[Failure of fixed-length convergence]
\label{prop:lacunary-logarithmic}
For each \(0<\theta<1\) there are independent counts \(Q_j\) satisfying
\eqref{eq:strong-logarithmic-condition}, with \(Q_1\sim\operatorname{Pois}(\theta)\)
and feasible conditioning in every degree, whose conditioned profile
is a single cycle with probability tending to one along \(n=4^a\).
In particular, the number of fixed points does not converge in law to
\(Q_1\).
\end{proposition}

\begin{proof}
Let \(P_j\sim\operatorname{Pois}(\theta/j)\) and let \(D_j\) be Bernoulli
with mean \(b_j=j^{-(1-\theta)/2}\) at \(j=4^a\), \(a\ge1\), and zero
at all other lengths. Take all these variables independent and put
\(Q_j=P_j+D_j\). The probability of changing \(P_j\) by adding
\(D_j\) is \(b_j\), so
\[
 |j\Pp(Q_j=1)-\theta|\le \theta^2/j+jb_j.
\]
For \(r\ge2\), convolution of the two masses gives
\[
 j\Pp(Q_j=r)\le
 \frac{\theta^r}{r!j}+jb_j\frac{\theta^{r-1}}{(r-1)!}.
\]
These bounds give \eqref{eq:strong-logarithmic-condition} with
\(e(j)=j^{-1}+jb_j\), since \(\sum_j b_j<\infty\) and the bounds on the
mass at \(r\ge2\) allow factorially decreasing \(c_r\). The probability
that all Bernoulli variables vanish is \(\prod_j(1-b_j)>0\).
The Poisson fixed points alone give positive conditioning
probability in every degree.

For \(n=4^a\) the sum of all earlier Bernoulli lengths is less than
\(n/3\). If \(D_n=0\), conditioning on the earlier Bernoulli variables
makes the required Poisson total a value \(w\in[2n/3,n]\).
For \(w\le n\) its mass is
\(e^{-\theta H_n}h_w(\theta)=O_\theta(n^{-1})\) uniformly in this interval.
If \(D_n=1\), the conditioning event occurs exactly when all earlier
Bernoulli variables and all the Poisson counts through \(n\) vanish.
Its joint probability with \(D_n=1\) is
\[
 b_n e^{-\theta H_n}\prod_{j<n}(1-b_j)
   \asymp_\theta n^{-(1+\theta)/2}.
\]
Since \((1+\theta)/2<1\), the conditioning selects this event with
probability tending to one. It supplies a single \(n\)-cycle.
\end{proof}

For the Ewens law the effect of taking the infimum over all degrees
is exact.

\begin{proposition}[The all-degree threshold]
\label{prop:all-degree-threshold}
For the Ewens law of parameter \(\theta>0\), let
\(m_\theta^{\mathrm{all}}\) be the least number of samples for
which the infimum of the generation probabilities over \(n\ge2\)
is positive, as in \cite[Eq.~(1.3), p.~856]{brito}. Take the
value \(\infty\) if there is no such number. Then
\begin{equation}\label{eq:all-degree-threshold}
 m_\theta^{\mathrm{all}}=\max\{3,m_\theta\}.
\end{equation}
\end{proposition}

\begin{proof}
A positive infimum implies a positive limit inferior. No pair
invariably generates \(S_6\) \cite[Proposition~4.10(i)]{tracey}, so
\(m_\theta^{\mathrm{all}}\ge\max\{3,m_\theta\}\).

Suppose now that \(m_\theta<\infty\) and put \(m=\max\{3,m_\theta\}\).
Adding samples preserves invariable generation. Hence there are
\(c>0\) and \(n_0\ge3\) such that \(G_{n,m}(\theta)\ge c\) for
\(n\ge n_0\). It remains to check the individual smaller degrees.
In every degree \(n\ge3\), an \(n\)-cycle, an \((n-1)\)-cycle,
and a transposition invariably generate \(S_n\). Indeed, any
conjugates of the first two generate a transitive group in which
the stabilizer of the fixed point of the second is transitive on
the remaining points. The group is 2-transitive, so adjoining the
transposition supplies all transpositions. Each of these three cycle
types has positive Ewens probability. In degree two a transposition
suffices. Therefore \(G_{n,m}(\theta)>0\) for every \(n\ge2\), and
\(\min\{c,G_{2,m}(\theta),\ldots,G_{n_0-1,m}(\theta)\}>0\)
bounds the infimum from below.
\end{proof}

Consequently, for \(0<\theta<\theta_c\) the two thresholds are two
and three under the large-degree and all-degree conventions,
respectively. At \(\theta_c\) both equal three.

For general permutation laws, the following criterion separates the
contributions of the short counts, the signs, and the transitive
overgroups.

\begin{proposition}[Generation from short cycle counts]\label{prop:short-generation-criterion}
Let \(m\ge2\), let \(b_n\to\infty\) with \(b_n\le n/2\), and take
\(m\) independent samples from a law on \(S_n\). Let \((Q_j)_{j\ge1}\)
be independent nonnegative integer counts and write \(q\) for the probability
that the sumsets of \(m\) independent copies of this profile intersect
only at zero. Suppose the samples can be coupled with these \(m\)
profiles so that, with probability \(1-o(1)\), every sample count
through length \(b_n\) is at most its corresponding profile count.
Suppose also that
\begin{enumerate}
\item the probability of a common invariant size in \((b_n,n/2]\) is \(o(1)\),
\item given the counts of the \(m\) samples through \(b_n\), the conditional
probability of at least one odd sample is at least a constant \(c>0\)
outside an event of probability \(o(1)\), and
\item with probability \(1-o(1)\) every transitive subgroup containing the
first sample contains \(A_n\).
\end{enumerate}
Then the limit inferior of the invariable generation probabilities is at
least \(cq\). If the \(Q_j\) satisfy \eqref{eq:strong-logarithmic-condition}
with parameter \(\theta\), and \(m(1-\theta\log 2)>1\), then \(q>0\)
if and only if \(\Pp(Q_j=0)>0\) for every \(j\).
\end{proposition}

\begin{proof}
Let \(E_n\) be the event that there is no common invariant size in
\([1,b_n]\). On the coupling event, every such size would also belong
to each of the \(m\) comparison sumsets, so \(\Pp(E_n)\ge q-o(1)\).
The event \(E_n\) depends only on the short counts. By the sign
hypothesis, \(E_n\) occurs together with an odd sample with probability
at least \(cq-o(1)\). Exclude the two events of probability \(o(1)\)
in the first and third hypotheses. The absence of a common proper
invariant size forces any choice of conjugates to generate a transitive
group. This group contains \(A_n\) by the overgroup hypothesis, which
is invariant under conjugation. An odd sample then gives \(S_n\).

For the last assertion, if \(\Pp(Q_j=0)=0\) for some \(j\), then \(j\)
belongs to every sumset almost surely. Conversely, suppose every
\(\Pp(Q_j=0)>0\). As in Remark~\ref{rem:summable-perturbation},
\eqref{eq:strong-logarithmic-condition} bounds the total variation distance
from \(\Law(Q_j)\) to \(\operatorname{Pois}(\theta/j)\) by
\(O(e(j)/j+j^{-2})\). Independent maximal couplings therefore have
summable discrepancy probabilities.
Apply the coupling independently to the \(m\) profiles. The probability
of any discrepancy at \(j>M\) tends to zero with \(M\).
By Lemma~\ref{lem:common-size-tail}, the probability of a common Poisson
sum greater than \(M\) also tends to zero. The \(m\) profiles of
\(Q_j\)'s restricted to \(j>M\) therefore have no common positive
sum with probability tending to one. Fix \(M\) for which this probability
is positive. Independently, every count through \(M\) in all \(m\)
profiles vanishes with probability \(\prod_{j=1}^M\Pp(Q_j=0)^m>0\).
Together the two events give \(q>0\).
\end{proof}

When the counts of the finite permutations converge on each fixed window
to the defining counts \(Q_j\), a diagonal choice of \(b_n\) gives
a coupling in which the counts through \(b_n\) agree with probability
\(1-o(1)\). The other hypotheses of
Proposition~\ref{prop:short-generation-criterion} must hold for the
same \(b_n\). In particular the printed condition alone gives no
control of a common size comparable to \(n\).

\subsection{Domination of the short counts}\label{subsec:earlier-domination}

Question~2 of \cite[p.~860]{brito} uses the notation \(h^{-1}(2)\).
For the staircase in \cite[Eq.~(1.4)]{brito}, the set of positive
parameters at which \(h\) equals two is \((0,\theta_c]\). We take
\(h^{-1}(2)\) in that question to mean its right endpoint \(\theta_c\),
and asymptotic domination to mean joint stochastic domination of the
limiting cycle counts through each fixed length. For domination from below,
as in \cite[Corollary~1.2, pp.~856--857]{brito}, suppose each such vector
dominates independent Poisson coordinates of means \(\lambda_0/j\).
Proposition~\ref{thm:strict-converse} then gives failure of invariable
transitivity for any fixed number of independent samples. Indeed a common
positive sum in a fixed window is an increasing event in the counts,
and its probability in the Poisson model tends to one as the window grows.

Domination from above, as in Question~2, allows a common invariant
size of order \(n\). For \(n\ge4\), take a
permutation uniformly from the conjugacy class with cycle lengths
\(\lfloor n/2\rfloor\) and \(\lceil n/2\rceil\). Its counts through
each fixed length are eventually zero, and conjugation can align
an invariant \(\lfloor n/2\rfloor\)-set across all samples.

The obstruction persists for nondegenerate Poisson limits with
strict domination at each coordinate.

\begin{proposition}[Strict domination of the limiting counts]
\label{prop:short-domination}
Let \(Q_j\sim\operatorname{Pois}(\theta_c/(j+1))\) independently.
For each \(n\), condition \((Q_1,\ldots,Q_n)\) on total length \(n\)
and choose a permutation uniformly from the corresponding conjugacy
class. The limiting cycle counts through each fixed length are
the independent \(Q_j\), which are strictly stochastically dominated,
coordinate by coordinate, by the independent
\(\operatorname{Pois}(\theta_c/j)\) counts. With probability tending
to one, two independent samples have a common proper invariant size.
\end{proposition}

\begin{proof}
Every degree is feasible since \(Q_1\) has full support. Writing
\(x=\theta_c/(j+1)\) and using \(1-e^{-x}\le x\) gives
\[
 |j\Pp(Q_j=1)-\theta_c|\le\frac{\theta_c+\theta_c^2}{j},
 \qquad
 j\Pp(Q_j=a)\le\frac{\theta_c^a}{a!j}\quad(a\ge2).
\]
The uniform logarithmic condition holds with \(e(j)=1/j\),
and the limiting fixed-length counts are the \(Q_j\) by
\cite[Theorem~3.1]{abt2000}. Adding independent Poisson counts
of means \(\theta_c/[j(j+1)]\) gives the critical Poisson
profile, proving the claimed domination simultaneously in every
fixed collection of coordinates. Since \(\sum_j\theta_c/[j(j+1)]<\infty\),
the total variation discrepancies are summable. Remark~\ref{rem:summable-perturbation}
now gives the common proper invariant size with probability tending to one.
\end{proof}

Below the first discontinuity there is a positive statement for two
samples. If the comparison counts in
Proposition~\ref{prop:short-generation-criterion} are Poisson of means
\(\theta/j\) with \(\theta<\theta_c\), then its lower bound is
\(c\,p_{0,2}(\theta)>0\). The proposition requires control through a
growing window and of the common sizes above that window.

\section{Further questions}

For \(m\ge3\) and \(\theta=(1-1/m)/\log2\), is \(p_{0,m}(\theta)\)
positive? By Corollary~\ref{cor:threshold-equivalence} the answer
determines whether \(s_\theta=m_\theta=m\) or \(s_\theta=m_\theta=m+1\).
In fact, a zero answer would give infinitely many common sums after
every fixed lower cutoff. To see this, condition on the positive-probability
event that all \(m\) profiles have no copies through length \(M\).
If \(p_{0,m}(\theta)=0\), a common positive sum still exists almost surely
under this conditioning. The counts above \(M\) retain their original
law, and any such sum exceeds \(M\). Intersecting over the integer
cutoffs \(M\) proves the assertion. The difficulty at a higher jump
is detecting a positive total common to \emph{every} profile.

At \(\theta_c\), Proposition~\ref{s3:critical-rate} bounds
\(\sqrt{\log k}\,i_{\theta_c}(n,k)/\E I_{n,k}\) above and away from
zero for \(n\ge2k\ge4\). Does this ratio converge for each fixed
\(d\)? For \(d\to\infty\), can a limit depend on the growth of \(d\)?
At \(\lambda_0\), one may ask the same questions about
\(\sqrt{\log k}\,i_{\lambda_0}(n,k)\), or fix the limiting ratio
\(k/n\).

For \(\theta_c<\theta<\lambda_0\), Proposition~\ref{prop:large-theta}
gives the order with logarithmic factor \((1+\log k)^{-3/2}\).
The asymptotic equivalent for general intermediate parameters
remains to be determined. If \(n\to\infty\) first, \eqref{s2:feller} gives the probability that a
Poisson profile hits \(k\). At \(\theta=1\) its asymptotic has a smooth
positive periodic factor which Green and Sawhney \cite{gs} conjecture
is nonconstant. Do periodic factors occur at other parameters?

For \(0<\theta<\theta_c\), Proposition~\ref{prop:constant-support}
and Theorem~\ref{s3:main} give \(c_0<c_d<c_\infty\) for \(d\ge1\)
and \(c_d\to c_\infty\) as \(d\to\infty\). Is \(c_d\) strictly
increasing in \(d\)? The support expression in the proof involves
the overlap of a sumset with its translate. For the merged profile
through length two, condition on no copy of length one and exactly
one copy of length two. Then \(S_2\) is uniform on \(\{0,2\}\).
At shift one both reciprocal counts are one. At shift two they are
\(1/2\) when \(S_2=0\) and one when \(S_2=2\), so their mean is \(3/4\).
The mean is \(1/2\) at shift zero and one at shift three. Averaging over the merged profile could change this behavior.

\section{Disclosure of AI Use}

The ideas articulated in this paper are my own. However, due to laziness and a shortage of free time, I fed my fragmented notes to GPT-6.2 Astra (OpenAI) and asked it to synthesize everything into a completed manuscript. I take responsibility for the mathematical arguments, computations, and references in this paper.

\end{document}